\documentclass[3p,times]{elsarticle}
\usepackage{mathtools}
\usepackage{amssymb}
\usepackage{amsfonts}
\usepackage{dsfont}
\usepackage{amsthm}
\usepackage{cases}
\usepackage{bm}
\usepackage{geometry}
\usepackage{indentfirst}
\usepackage{tabularx}
\usepackage{booktabs}
\usepackage{graphicx}
\usepackage{graphicx}
\usepackage{ulem}
\usepackage{xspace}
\usepackage{xcolor}
\usepackage{aliascnt}

\usepackage{multirow}

\newcommand{\bea}{\begin{eqnarray}}
\newcommand{\eea}{\end{eqnarray}}
\newcommand{\lp}{\left (}
\newcommand{\rp}{\right )}
\newcommand{\bp}{\begin{pmatrix}}
\newcommand{\ep}{\end{pmatrix}}

\newcommand{\R}{\mathds{R}}

\newcommand{\dtdx}{\frac{\Delta t}{\Delta x}}

\newcommand{\dt}{{\Delta t}}
\newcommand{\dx}{{\Delta x}}
\newcommand{\dy}{{\Delta y}}
\newcommand{\dz}{{\Delta z}}
\newcommand{\eps}{\varepsilon}

\def\be{\begin{equation}}
\def\ee{\end{equation}}

\newcommand{\Frac}{\displaystyle\frac}

\newcommand{\Max}{\displaystyle\max}

\newcommand{\Q}{\mathbf{Q}}
\newcommand{\uu}{\bm{u}}
\newcommand{\ww}{\bm{w}}

\newcommand{\xx}{\bm{x}}
\newcommand{\bss}{\bm{\sigma}}

\newcommand{\U}[1]{{\bm{#1}}}

\newfont{\numerikEleven}{ecrm1000}
\newfont{\numerikTen}{cmss10}
\newfont{\numerikNine}{cmss9}
\newfont{\numerikEight}{cmss8}
\newfont{\numerikSeven}{cmss7}

\newcommand{\revised}[1]{\textcolor{black}{#1}}

\journal{Journal of Computational Physics}

\begin{document}
	
\begin{frontmatter}
\title{High order pressure-based semi-implicit IMEX schemes for the 3D Navier-Stokes equations at all Mach numbers}
\author[ferrara]{Walter Boscheri$^*$}
\ead{walter.boscheri@unife.it}
\cortext[cor1]{Corresponding author}

\author[ferrara]{Lorenzo Pareschi}
\ead{lorenzo.pareschi@unife.it}

\address[ferrara]{Department of Mathematics and Computer Science, University of
	Ferrara, Ferrara, Italy}


\begin{abstract}
This article aims at developing a high order pressure-based solver for the solution of the 3D compressible Navier-Stokes system at all Mach numbers. We propose a cell-centered discretization of the governing equations that splits the fluxes into a fast and a slow scale part, that are treated implicitly and explicitly, respectively. A novel semi-implicit discretization is proposed for the kinetic energy as well as the enthalpy fluxes in the energy equation, hence avoiding any need of iterative solvers. The implicit discretization yields an elliptic equation on the pressure that can be solved for both ideal gas and general equation of state (EOS). A nested Newton method is used to solve the mildly nonlinear system for the pressure in case of nonlinear EOS. High order in time is granted by implicit-explicit (IMEX) time stepping, whereas a novel CWENO technique efficiently
implemented in a dimension-by-dimension manner is developed for achieving high order in space for the discretization of explicit convective and viscous fluxes. A quadrature-free finite volume solver is then derived for the high order approximation of numerical fluxes. Central schemes with no dissipation of suitable order of accuracy are finally employed for the numerical approximation of the implicit terms. Consequently, the CFL-type stability condition on the maximum admissible time step is based only on the fluid velocity and not on the sound speed, so that the novel schemes work uniformly for all Mach numbers. Convergence and robustness of the proposed method are assessed through a wide set of benchmark problems involving low and high Mach number regimes, as well as inviscid and viscous flows.  
\end{abstract}

\begin{keyword}	
All Mach number flow solver \sep
Asymptotic preserving methods \sep
Semi-implicit IMEX schemes \sep
3D compressible Euler and Navier-Stokes equations \sep
Quadrature-free WENO \sep
General equation of state (EOS)
\end{keyword}
\end{frontmatter}

\section{Introduction} \label{sec:intro}
The unsteady compressible Navier-Stokes equations constitute a mathematical model for the simulation of a wide set of applications in fluid mechanics, that involve aerospace and
mechanical engineering as well as environmental engineering \cite{Vidal,Gresho1990,BrachetDNS,Ala2,Oriol2015,Patankar}. Atmospheric flows, geophysical flows in oceans, rivers and lakes can be described relying on the Navier-Stokes model, which is also used in industrial applications such as the design of  wind or water turbines, aircraft engines and cars. The governing equations are based on the physical principle of conservation and they can be derived from the conservation of mass, momentum and total energy. The compressible Navier-Stokes equations already embed several simplified sub-systems, such as the compressible Euler equations in the case of inviscid flows or the incompressible Navier-Stokes equations, that can be retrieved in the zero Mach number limit. The Mach number, which is the ratio between the fluid velocity and the sound speed, describes the regime of the fluid under consideration. High Mach number situations are typically encountered in industrial engineering, whereas geophysical phenomena mostly involve low Mach number flows. 


The numerical methods developed for the solution of high and low Mach number problems are quite different, because of the nature of the governing equations. For the high Mach number case explicit upwind finite difference and Godunov-type finite volume methods are very popular \cite{Godunov1959, HLL1983, LW1960, LeVequeBook, Munz1994, OS1982}. On the other hand, in the incompressible regime the elliptic behavior of the pressure introduces a very severe restriction on the maximum admissible time step for low Mach number flows. Indeed, the CFL-type stability condition for explicit methods depends also on the sound speed which becomes predominant in the zero Mach limit. Furthermore, in \cite{Dellacherie1} the effect of numerical viscosity on the slow waves introduced by upwind-type schemes is proven to degrade the accuracy. As a consequence, implicit strategies for time discretization have been proposed in order to avoid the acoustic CFL restriction and enlarge the time step. However, fully implicit methods imply the solution of large nonlinear systems that are computationally very expensive and in which the convergence is numerically very difficult to control. In addition, in many realistic scenarios, both high and low Mach regimes coexist and can arise during the simulation without being predicted in advance, thus needing the design of numerical methods that can deal with all Mach numbers.    

This is the reason behind the research activity carried out in the recent past for investigating an alternative strategy to treat problems with multiple time scales. A successful idea consists in treating implicitly only one part of the system to be solved while keeping the remaining explicit, thus both incompressible and compressible regimes can be handled \cite{Patankar,Turkel,Deg_Del_Sang_Vig,DegTan,Chalons,Bog,DLDV2018}. This approach permits to design space and time discretizations in which the implicit part of the system is relatively easy to be inverted, typically avoiding nonlinear systems, while keeping robustness and shock-capturing properties in the explicit part. There are mainly two classes of schemes that allows to deal with split sub-systems, namely one for the fast and the other for the slow scale phenomena, that are treated implicitly and explicitly, respectively. \revised{The first class is given by implicit-explicit (IMEX) methods \cite{AscRuuSpi,BosRus,PR_IMEX,BPR2017} or, more in general, by the so-called partitioned schemes \cite{Hofer}. IMEX schemes are proven to be very effective for many applications. Their main feature is to achieve high order under a time step stability constraint independent of the values of the fast scale, and to satisfy the Asymptotic Preserving (AP) property, meaning that the limit model is consistently reproduced at the discrete level \cite{JINAP1999, KLARAP1999}.} The other class is represented by semi-implicit methods \cite{Klein,Casulli1990,Casulli1999,BosFil2016,BosRus2019}, which have also gained visibility in the past years. 
Here, the idea is to obtain a linearly implicit scheme for the stiff terms in the governing equations, thus avoiding any need of iterative methods. In \cite{Toro_Vazquez12} a flux splitting for the Euler equations is proposed that aims at obtaining an advection and a pressure sub-system. Both sub-systems are demonstrated to be hyperbolic and in \cite{Dumbser_Casulli16} the advection system is treated explicitly, while the pressure sub-system is discretized implicitly, so that the time step is only limited by the fluid velocity and not by the acoustic speed. Staggered meshes are employed, hence allowing for compact stencils and permitting to recover by construction the divergence free constraint of the velocity field in the low Mach limit, along the lines of \cite{Casulli1990,Casulli1999}. Following these ideas, high order semi-implicit methods coupled with discontinuous Galerkin (DG) space discretizations on unstructured staggered meshes have been forwarded for compressible and incompressible flows \cite{TavelliIncNS,TavelliCNS}, on dynamic adaptive meshes \cite{Fambri2017,Fambri2017b} and for axially symmetric flows \cite{Ioriatti2015,Ioriatti2018}.

\revised{It is worth to notice, that the two above approaches, IMEX and semi-implicit, have been generalized under a unified framework in \cite{BosFil2016}. This generalization, permits the construction of high order linearly implicit schemes by using the standard formalism of IMEX Runge-Kutta methods. In this paper, we will rely on such methodology to achieve efficient high-order accurate semi-implicit discretizations of the Navier-Stokes system for all Mach number flows.} 

Recently, an all Mach solver for the 3D Euler equations has been designed \cite{BDLTV2020}, which is a cell-centered second order accurate finite volume method with IMEX time stepping. An elliptic equation on the pressure is solved at the aid of an iterative Picard algorithm. Once the pressure is computed, it is used for advancing in time the momentum and the total energy. The work has then been extended in \cite{BDT_cns} to the full Navier-Stokes equations with implicit viscosity treatment. In \cite{Avgerinos2019} the two-dimensional Euler equations are considered at all Mach number regimes, presenting second order IMEX schemes based on the solution of an elliptic equation on the energy. A finite volume solver has been forwarded in \cite{Degond2} for inviscid and viscous compressible fluids in two space dimensions at high and low Mach flows. There, an elliptic equation on the enthalpy allows to treat also general equations of state (EOS) that link the internal energy with the density and the enthalpy. A similar approach based on the solution of a pressure wave equation can be found in \cite{Dumbser_Casulli16}, where ideal gas and general cubic EOS are considered. The aforementioned references are at most second order accurate in space and time. High order semi-implicit schemes for the isentropic Euler equations on two-dimensional Cartesian meshes are described in \cite{BosRus2019}, where the seminal work presented in \cite{BosFil2016} is applied to low Mach flows.

In this article we present a third-order semi-implicit scheme on collocated Cartesian grids for the solution of the compressible Navier-Stokes equations at all Mach numbers. The flux splitting proposed in \cite{Toro_Vazquez12} requires an implicit sub-system to be solved for the pressure. The novel semi-implicit discretization proposed in this work splits the kinetic energy contribution as well as the enthalpy fluxes in the energy equation into an explicit and an implicit part. Differently from what presented in \cite{Dumbser_Casulli16,BDLTV2020}, no iterative solvers are needed anymore and the solution of the pressure system is directly computed. General equations of state yield a mildly nonlinear system that can then be handled relying on a nested Newton technique developed in \cite{NestedNewton}. An efficient CWENO reconstruction that is carried out dimension-by-dimension is used for achieving high order of accuracy in space, while IMEX Runga-Kutta time stepping is adopted following \cite{BosFil2016}. Finally, a quadrature-free finite volume scheme is developed for the convective explicit part of the system, which also contains the viscous fluxes of the Navier-Stokes equations. The scheme is also proven to be asymptotically preserving, hence it recovers the limit model in the stiff limit. To overcome the appearance of spurious oscillations due to high order discretization in time, a new limiting strategy is proposed which is based on a convex combination between high and first order numerical solution. Applications to inviscid and viscous compressible flows in low and high Mach regimes are shown, demonstrating the accuracy and the robustness of the novel method.

The rest of the paper is organized as follows. In Section \ref{sec:equations} the compressible Navier-Stokes equations are described. We also introduce the low Mach scaling of the governing partial differential equations (PDE) and the continuous model retrieved in the stiff limit. The novel semi-implicit scheme is detailed in Section \ref{sec:num_scheme}. Firstly, a semi-discretization in time is explained, then the fully discrete first order scheme is derived. Details of the high order extension in time and in space are then given and finally the quadrature-free finite volume scheme for the treatment of explicit terms is presented. The limiting strategy adopted to reduce spurious oscillations when high order time discretizations are adopted is detailed at the end of this section. Numerical convergence studies and applications to a wide set of test problems is shown in Section \ref{sec:num_results}. A concluding section finalizes the article where we draw some conclusions and present
an outlook to future research.


\section{Governing equations} \label{sec:equations}
Let $\Omega \in \R^d$ represent a bounded domain in the space dimension $d \in \{ 1, 2, 3 \}$, which is defined by spatial position vector $\xx \in \Omega$ and time variable $t \in \R_+$. The compressible Navier-Stokes equations write
\begin{equation}   
\frac{\partial }{\partial t} \left( \begin{array}{c}
\rho \\ \rho \uu \\ \rho E
\end{array} \right) + \nabla \cdot \left( \begin{array}{c}
\rho \uu \\ \rho \uu \otimes \uu + p \mathbf{I} \\ \rho k \uu + h \rho \uu
\end{array} \right) = \nabla \cdot \left( \begin{array}{c} 0 \\ \bss \\ \bss \uu + \lambda \nabla T \end{array} \right),
\label{eqn.cns}
\end{equation}
with $\mathbf{I}$ being the identity matrix. $\rho(t,\xx) > 0$ is the density of the fluid, $\uu(t,\xx) \in \R^d$ denotes the velocity vector, $\rho E(t,\xx)$ represents the total energy with the specific kinetic energy $k$, the specific internal energy $e$ and the specific enthalpy $h$. The fluid pressure is denoted by $p(t,\xx) > 0$ and $T(t,\xx) > 0$ refers to the fluid temperature with $\lambda$ representing the thermal conductivity. In $\R^3$ one has $\xx=(x,y,z)$ and $\uu=(u,v,w)$. The right hand side of system \eqref{eqn.cns} is conveniently formulated by introducing the stress tensor $\bss$ which under Stokes hypothesis is 
\begin{equation}
\bss = \mu \left(\nabla \uu + \nabla \uu^\top \right) - \frac{2}{3} \left( \mu \nabla \cdot \uu \right) \, \mathbf{I},
\label{eqn.sigma}
\end{equation} 
where $\mu$ is the viscosity of the fluid. A thermal equation of state $p=p(T,\rho)$ and a caloric equation of state $e=e(T,\rho)$ are required to close system \eqref{eqn.cns}. Typically, the temperature is canceled from these two equations of state, yielding one single relation of the form $e=e(p,\rho)$, which will be adopted in this work. We assume that the internal energy is a non-negative and non-decreasing function of the fluid pressure $p$. Furthermore, the relation between the viscosity coefficient and the fluid temperature is governed by Sutherland's law, that is
\begin{equation}
\mu(T) = \mu_0 \left(\frac{T}{T_0}\right)^{\beta} \frac{T_0+s}{T+s},
\label{eqn.Sutherland}
\end{equation}
with parameters $\mu_0$, $T_0$, $\beta$ and $s$. Notice that constant viscosity is retrieved if $\beta=1$ and $s=0$. The ratio of specific heats of the gas at constant pressure $c_p$ and at constant volume $c_v$ is $\gamma=c_p/c_v$ and the specific heat at constant volume $c_v$ is determined by $c_v=R/(\gamma-1)$ with $R$ being the gas constant which is assumed to be $R=0.4$. Finally, the specific kinetic energy $k$ and the specific enthalpy $h$ are given by the following relations:
\begin{equation}
k=\Frac12 \uu^2, \qquad h=e + \frac{p}{\rho}.
\label{eqn.eDef}
\end{equation}
Let observe that the total energy flux in \eqref{eqn.cns} is written as
\begin{equation}
\uu(\rho E + p) = \rho k\, \uu  + \rho h\,\uu,
\label{eqn.fluxsplit}
\end{equation}
according to \cite{Toro_Vazquez12}, thus introducing a \textit{flux splitting} which will be extremely important for the numerical methods developed in this work.  

\subsection{Ideal gas EOS} \label{ssec.idealgas}	
If an ideal gas is considered, the thermal and caloric equation of state (EOS) are given by
\begin{equation}
\frac{p}{\rho}=RT, \qquad e=c_v T.
\label{eqn.idealEOS_twoeqn}
\end{equation}
The temperature can be eliminated using both expressions in \eqref{eqn.idealEOS_twoeqn}, thus leading to an equation of state of the form $e(p,\rho)$, that is
\begin{equation}
e(p,\rho) = \frac{p}{(\gamma-1)\rho}.
\label{eqn.idealEOS}
\end{equation}
Here the relation between pressure $p$ and internal energy $e$ is linear.

\subsection{Redlich-Kwong EOS} \label{ssec.RKgas}
Let now consider a general cubic equation of state, which according to \cite{Vidal} takes the form 
\begin{equation}
p(T,\rho)=\frac{RT}{v-b} - \frac{a(T)}{(v-b r_1) \, (v-b r_2)}.
\label{eqn.pEOS}
\end{equation}
Here, $v=1/\rho$ is the specific volume, $b$ represents the co-volume and $r_1,r_2$ are two parameters. The function $a(T)$ is related to the attraction term in the EOS. The caloric equation of state which corresponds to \eqref{eqn.pEOS} writes \cite{Vidal} 
\begin{equation}
e(T,\rho) = c_v T + \frac{a(T)-T a^\prime(T)}{b} \, U(v,b,r_1,r_2),
\label{eqn.eEOS}
\end{equation}
with
\begin{equation}
a^\prime(T) = \frac{d a(T)}{dT}, \qquad U(v,b,r_1,r_2) = \frac{1}{r_1-r_2} \ln \left( \frac{v-b r_1}{v-b r_2} \right).
\end{equation}

Different equation of states can be derived by appropriate choices of the parameters in \eqref{eqn.pEOS}-\eqref{eqn.eEOS}. For instance, the Redlich-Kwong EOS is obtained by setting $r_1=0$ and $r_2=-1$, with the attraction term given by $a(T)=1/(2\sqrt{T})$. In this case the EOS yields a nonlinear relation between internal energy and pressure. In order to compute a function of the form $e(p,\rho)$, we first need to determine the temperature from the thermal equation of state \eqref{eqn.pEOS}. In our approach we solve the nonlinear equation $p(T,\rho)$ numerically relying on a simple and efficient Newton method. Once the temperature is known, it can be inserted into the caloric EOS \eqref{eqn.eEOS} to obtain the relation $e(p,\rho)$.

\subsection{Scaling of the Navier-Stokes equations} \label{ssec.scalingNS}
The governing equations \eqref{eqn.cns} can be rescaled relying on the following scaled variables:
\begin{equation}
\tilde \rho=\rho/\rho_0, \qquad \tilde \uu=\uu/u_0, \qquad \tilde p=p/p_0, \qquad \tilde E=\frac{\rho_0 E}{p_0}, \qquad \tilde T=T/T_0, \qquad \tilde \xx=\xx/x_0, \qquad \tilde t=t/t_0,
\label{eqn.rescaling}
\end{equation}
where $\rho_0$, $p_0$, $x_0$, $t_0$, $u_0=x_0/t_0$ and $T_0=p_0/\rho_0$ are typical values referred to the problems under consideration. Furthermore, let $\tilde \mu=\mu/\mu_0$ and $\tilde \lambda=\lambda/\lambda_0$ be the rescaled coefficients for the viscosity and the thermal conductivity, respectively, and let us introduce a stiffness parameter $\varepsilon$ 
\begin{equation}
\varepsilon=\frac{\rho_0 u_0^2}{p_0}, 
\label{eqn.eps}
\end{equation} 
which is related to the global Mach number $M=u_0/c$ and characterizes the flow and the nondimensionalisation. The sound speed $c$ is then given by $c^2 = \left( \Frac{\partial p}{\partial \rho}\right)_s$ with $s$ representing the entropy. Using the definitions \eqref{eqn.rescaling}-\eqref{eqn.eps} and omitting the tildes, the rescaled Navier-Stokes equations read
\begin{equation}   
\frac{\partial }{\partial t} \left( \begin{array}{c}
\rho \\ \rho \uu \\ \rho E
\end{array} \right) + \nabla \cdot \left( \begin{array}{c}
\rho \uu \\ \rho \uu \otimes \uu + \frac{1}{\varepsilon} p \mathbf{I} \\ \varepsilon \, \rho  k \uu + h \rho \uu
\end{array} \right) = \nabla \cdot \left( \begin{array}{c} 0 \\ \bss \\ \varepsilon \, \bss \uu + \lambda \nabla T \end{array} \right).
\label{eqn.cns_rescaled}
\end{equation}
The right hand side of \eqref{eqn.cns_rescaled} can be expanded at the aid of the the Reynolds and Prandtl numbers
\begin{equation}
Re = \frac{u_0 x_0}{\nu}, \qquad Pr = \frac{\mu \gamma c_v }{\lambda},
\end{equation}
with the kinematic viscosity $\nu=\frac{\mu}{\rho}$, hence yielding
\begin{equation}
\nabla \cdot \left( \begin{array}{c} 0 \\ \bss \\ \varepsilon \, \bss \uu + \lambda \nabla T \end{array} \right) = \nabla \cdot \left( \begin{array}{c}
	0 \\ \frac{1}{Re}  \nabla \cdot \left[ \rho \left(\nabla \uu + \nabla \uu^\top \right) - \frac{2}{3} \left( \rho \nabla \cdot \uu \right) \, \mathbf{I} \right] \\ \frac{\varepsilon}{Re}  \nabla \cdot \left[ \left( \rho \left(\nabla \uu + \nabla \uu^\top \right) - \frac{2}{3} \left( \rho \nabla \cdot \uu \right) \, \mathbf{I} \right) \cdot \uu \right]  + \frac{c_p}{Re \cdot Pr} \nabla T 
\end{array} \right).
\end{equation}
Looking at the rescaled equations \eqref{eqn.cns_rescaled} it is evident that the stiffness is originated in the momentum equations by the pressure waves. As a consequence, the energy equation turns out to be stiff as well. In particular, the energy can be decomposed into a pressure and a kinetic part, corresponding to internal and kinetic energy contribution. In the stiff regime the pressure evolves very fast, implying the same for the component of the energy related to the pressure, i.e. the internal energy $e$. 

\subsection{Low Mach limit of the Navier-Stokes equations} \label{ssec.lowMachNS}	
Let us now investigate the limit of the Navier-Stokes equations in the case $\varepsilon \to 0$. In \cite{Dou} the low Mach limit is studied in a bounded domain, while a fully three dimensional space is considered in \cite{Ala2}. Here, we only briefly recall the formal limit obtained with an ideal gas EOS \eqref{eqn.idealEOS}. On the boundary $\partial \Omega$ of the computational domain the following conditions must be imposed: 
\begin{equation}\label{bound_cond}
\uu(t,\xx)\cdot \mathbf{n}=0, \quad \frac{\partial T}{\partial n}(t,\xx)=0, \quad \forall \xx\in \partial \Omega, \ t >0,
\end{equation}
with $\mathbf{n}$ denoting the unit outward normal vector to the boundary and $n$ its direction. The limit for $\varepsilon \to 0$ of system \eqref{eqn.cns} then writes \cite{Ala2}
\begin{eqnarray}   
\label{eq:cns_rho_lim}
&&\partial_t \rho + \nabla \cdot (\rho \uu) =  0, \\
\label{eq:cns_q_lim}
&&\partial_t (\rho \uu)   + \nabla \cdot \lp \rho \uu \otimes \uu \rp + \nabla p_1 = \nabla\cdot \bss, \\ 
\label{eq:cns_divu_lim}
&& \gamma \nabla \cdot \uu = (\gamma-1)\nabla\cdot\left(\frac{\lambda}{R}\nabla\left( \frac{1}{\rho}\right)\right),  
\end{eqnarray}
assuming that the limit pressure $p_1=\lim_{\eps\rightarrow 0}\frac{1}{\eps}\left(p- p_0\right)$ exists. Notice that if we set $\mu=0$ and $\lambda=0$, that is viscous forces and thermal conductivity are neglected, the well known low Mach limit for the compressible Euler equations is retrieved. Specifically, this implies a divergence free condition on the velocity field, i.e. $\nabla\cdot \uu=0$, which derives from the energy equation. In the case of the Navier-Stokes equations, in the low Mach limit the fluid is no more incompressible because of large temperature variations and heat conduction effects. 

Regardless the viscous or inviscid property of the fluid, in the low Mach regime the sound speed is much bigger than the fluid velocity, thus it corresponds to small values of $\varepsilon$ in \eqref{eqn.cns_rescaled}. From the numerical viewpoint, the maximum admissible time step $\Delta t=t^{n+1}-t^n$ for fully explicit schemes is given by a CFL-type stability condition that writes
\begin{equation}
\Delta t\leq  \textnormal{CFL} \min_{\Omega} \, \left(\Frac{\Max(|u\pm c/\sqrt{\varepsilon}|)}{\Delta x}+\Frac{\Max(|v\pm c/\sqrt{\varepsilon}|)}{\Delta y}+\Frac{\Max(|w\pm c/\sqrt{\varepsilon}|)}{\Delta z}+\max\left(\Frac{\lambda_v}{\Delta x^2}+\Frac{\lambda_v}{\Delta y^2}+\Frac{\lambda_v}{\Delta z^2}\right)\right)^{-1},
\label{eqn.dtex}
\end{equation}
where $\Delta x$, $\Delta y$, $\Delta z$ are the the characteristic mesh spacing along each spatial direction in 3D. The eigenvalues of the viscous sub-system for an ideal gas are given according to \cite{ADERNSE} by 
\begin{equation}
\lambda_v=\max \left( \frac{4}{3} \frac{\mu}{\rho}, \, \frac{\gamma \mu}{Pr \, \rho}\right).
\label{eqn.lambdav}
\end{equation}
$\Delta t$ is of order $\sqrt{\varepsilon}$ and tends to $0$ with $\varepsilon$, thus dictating severe limits in the maximum size of the time step. Furthermore, even if this constraint is satisfied and the scheme runs with very small time steps, explicit schemes are not capable to capture the correct asymptotic regime as discussed in \cite{Guillard,GH,Dellacherie1}.


\section{Numerical scheme} \label{sec:num_scheme}
For the sake of simplicity we present the discretization for the compressible Euler equations, that are retrieved by neglecting the terms on the right hand side of system \eqref{eqn.cns}. The compressible Navier-Stokes model will then be included with fully explicit discretization of the viscous forces in the momentum equations and the work of the viscous stress tensor in the energy equation, while keeping untouched the semi-implicit scheme developed for the inviscid hydrodynamics model. For an ideal gas \eqref{eqn.idealEOS} the heat flux in the energy equation can be treated implicitly because temperature can be written as $T=p/(R\rho)$, therefore it can be easily embedded in the semi-implicit solver for the pressure.

\subsection{First order semi-discrete scheme in time} \label{ssec.semi-disc}
The time discretization is based on a \textit{semi-implicit} approach which leads to the following scheme:
\begin{eqnarray}
&& \frac{\rho^{n+1}-\rho^{n}}{\dt} + \nabla \cdot \left(\rho \uu \right)^n = 0, \label{eqn.SI-rho}\\
&& \frac{(\rho \uu)^{n+1}-(\rho \uu)^{n}}{\dt} + \nabla \cdot \left( (\rho \uu)^n \otimes (\rho \uu)^n \right) + \frac{1}{\varepsilon} \nabla p^{n+1} = 0, \label{eqn.SI-mom}\\
&& \frac{(\rho e)^{n+1} + \frac{\varepsilon}{2} \frac{(\rho \uu)^n}{2\rho^{n}} \, (\rho \uu)^{n+1} -(\rho E)^{n}}{\dt} + \nabla \cdot \left(\varepsilon \, \rho  k \uu \right)^n + \nabla \cdot \left(h^{n}  (\rho \uu)^{n+1}  \right)= 0, \label{eqn.SI-e}\\
\end{eqnarray}
where the kinetic energy in the total energy definition splits into an explicit and an implicit contribution, namely
\begin{eqnarray}
(\rho E)^{n+1} &:=& (\rho e)^{n+1} + \varepsilon \frac{(\rho \uu)^n}{2\rho^{n}} \, (\rho \uu)^{n+1}.
\label{eqn.ek_si}
\end{eqnarray}

The scheme \eqref{eqn.SI-rho}-\eqref{eqn.SI-e} is written in flux form for all variables, hence it is locally and globally conservative. The terms involving pressure are treated implicitly, while the convective part of the system is discretized explicitly. This allows the time step to be free from any restriction based on acoustic waves, that is a desirable property in the low Mach regime when $\varepsilon \to 0$. As a consequence, the time step must satisfy a milder CFL stability condition which is based only on the material speed of the flow $|\uu|$, that is
\begin{equation}
\dt \leq \textnormal{CFL} \frac{\min_{\Omega} (\dx,\dy,\dz)}{\max_{\Omega} (|\uu|)},
\end{equation}  
with $\textnormal{CFL}<1$. The algorithm for the solution of the semi-implicit numerical scheme \eqref{eqn.SI-rho}-\eqref{eqn.SI-e} is made of the following steps. 
\begin{enumerate}
	\item The density equation can be solved explicitly, thus $\rho^{n+1}$ is readily obtained from \eqref{eqn.SI-rho}.
	\item The momentum equation \eqref{eqn.SI-mom} is then inserted into the energy equation \eqref{eqn.SI-e} yielding an \textit{elliptic} equation on the pressure:
	\begin{eqnarray}
	&&(\rho e)^{n+1} + \varepsilon \frac{(\rho \uu)^n}{2\rho^{n}} \, \left( (\rho \uu)^{n} - \dt \nabla \cdot \left( (\rho \uu)^n \otimes (\rho \uu)^n \right) - \frac{\dt}{\varepsilon} \nabla p^{n+1} \right) = \nonumber \\
	&&(\rho E)^{n} - \dt \nabla \cdot \left(\varepsilon \, \rho k \uu \right)^n - \dt  \nabla \cdot \left(h^{n} \left( (\rho \uu)^{n} - \dt \nabla \cdot \left( (\rho \uu)^n \otimes (\rho \uu)^n \right) - \frac{\dt}{\varepsilon} \nabla p^{n+1} \right) \right).
	\label{eqn.p1}
	\end{eqnarray}
    Notice that a semi-implicit discretization of the enthalpy flux in the energy equation \eqref{eqn.SI-e} leads to an explicit evaluation of the enthalpy, i.e. $h^n$, and an implicit treatment of the momentum, that is $(\rho \uu)^{n+1}$. Shifting the unknowns on the left hand side and multiplying by the stiffness factor $\varepsilon$, the pressure wave equation \eqref{eqn.p1} writes
    \begin{equation}
    \varepsilon \,  (\rho e)^{n+1} + \varepsilon \,  \frac{\dt}{2} \frac{(\rho \uu)^n}{\rho^{n}} \, \nabla p^{n+1} - \dt^2 \left( \nabla \cdot h^n \nabla p^{n+1} \right) = \varepsilon \, \left[ (\rho E)^{*} - \varepsilon \frac{\dt}{2} \frac{(\rho \uu)^n}{\rho^{n}} \,  (\rho \uu)^{*} - \dt \nabla \cdot \left( h^n \, (\rho \uu)^{*} \right) \right],
    \label{eqn.p-SI}
    \end{equation} 
    with the explicit quantities
    \begin{eqnarray}
    (\rho E)^{*} &=& (\rho E)^{n}  - \dt \nabla \cdot \left(\varepsilon \, \rho k \uu \right)^n, \\    
    (\rho \uu)^{*} &=& (\rho \uu)^{n} - \dt \nabla \cdot \left( (\rho \uu)^n \otimes (\rho \uu)^n \right).
    \end{eqnarray}
    The internal energy $(\rho e)^{n+1}$ must now be written in terms of the new pressure $p^{n+1}$ using the equation of state.
    
    \paragraph{Ideal gas EOS} According to \eqref{eqn.idealEOS} for a perfect gas the internal energy is given by
    \begin{equation}
    (\rho e)^{n+1} = \frac{p^{n+1}}{\gamma-1},
    \label{eqn.pLin}
    \end{equation}
    thus the elliptic equation \eqref{eqn.p-SI} constitutes a linear system that can be directly solved.
    	
    \paragraph{General EOS} For a general equation of state the relation between internal energy and pressure might be nonlinear, hence requiring the solution of the following nonlinear equation for the pressure:
    \begin{equation}
    g(p^{n+1}) = \mathcal{P}(p^{n+1}) + \mathcal{R} \, p^{n+1} - b^n = 0,
    \label{eqn.gp}
    \end{equation}
    with the definitions
    \begin{eqnarray}
    \mathcal{R} \, p^{n+1} := \varepsilon \,  \frac{\dt}{2} \frac{(\rho \uu)^n}{\rho^{n}} \, \nabla p^{n+1} - \dt^2 \left( \nabla \cdot h^n \nabla p^{n+1} \right), \\
    b^n := \varepsilon \, \left[ (\rho E)^{*} - \varepsilon \frac{\dt}{2} \frac{(\rho \uu)^n}{\rho^{n}} \,  (\rho \uu)^{*} - \dt \nabla \cdot \left( h^n \, (\rho \uu)^{*} \right) \right].
    \end{eqnarray}
    The term $\mathcal{P}(p^{n+1})$ contains the nonlinearity of \eqref{eqn.p-SI} due to the EOS \eqref{eqn.eEOS} for the internal energy $(\rho e)^{n+1}$. Recall that the new density $\rho^{n+1}$ is already known thanks to \eqref{eqn.SI-rho}. A Newton method is then used for solving the \textit{piecewise linear} equation \eqref{eqn.gp}, along the lines of the algorithm presented in \cite{NestedNewton}. The solution for the new pressure is iteratively obtained as
    \begin{eqnarray}
    g(p^{n+1,k+1}) = g(p^{n+1,k}) + \Delta p^k \frac{d g(p^{n+1,k})}{p^{n+1,k}} = 0,
    \label{eqn.pNewton}
    \end{eqnarray}
    with $k$ denoting the iteration index and $\Delta p^k=(p^{n+1,k+1}-p^{n+1,k})$. In practice, equation \eqref{eqn.pNewton} is directly solved for $\Delta p^k$, then the new pressure at the next Newton iteration is given by $p^{n+1,k+1}=p^{n+1,k} - \Delta p^k$. The Newton method stops when the prescribed tolerance $\delta=10^{-10}$ has been reached, e.g. $\Delta p^k<\delta$.
    \item The new pressure $p^{n+1}$ is used in \eqref{eqn.SI-mom} to compute the momentum $(\rho \uu)^{n+1}$ at the next time level.
    \item Finally, the total energy is simply updated relying on \eqref{eqn.ek_si}, which ensures thermodynamic compatibility between the new pressure $p^{n+1}$ and momentum $(\rho \uu)^{n+1}$.
    
\end{enumerate}

For the full Navier-Stokes system, the viscous contribution $\bss$ in the momentum equation as well as the work of the viscous stress tensor $\bss \cdot \uu$ in the energy equation are discretized explicitly and are formally embedded in the 
explicit quantities $(\rho \uu)^{*}$ and $(\rho E)^{*}$, respectively. If an ideal gas is considered, an implicit discretization is likely to be assumed for the temperature gradient in the energy equation. Since temperature can be easily written in terms of pressure, i.e. $T^{n+1}=p^{n+1}/(R \rho^{n+1})$, this contribution is added to the pressure wave equation \eqref{eqn.p-SI} and implicitly solved. 

\paragraph{Asymptotic preserving property} The limit of the governing equations \eqref{eqn.cns_rescaled} is given when $\varepsilon \to 0$. The expansion of a generic variable $m$ in powers of the stiffness parameter $\varepsilon$, that is $m = m_{(0)} + \varepsilon m_{(1)} + \varepsilon^2 m_{(2)} + \ldots$, is applied to all variables involved in the governing rescaled Navier-stokes model, which is here assumed with $\mu=\lambda=0$. These expressions are then inserted into the semi-discrete scheme \eqref{eqn.SI-rho}-\eqref{eqn.SI-e} and only leading order terms are considered, thus obtaining
\begin{eqnarray}
&& \frac{\rho_{(0)}^{n+1}-\rho_{(0)}^{n}}{\dt} + \nabla \cdot \left(\rho \uu \right)_{(0)}^n = 0, \label{eqn.SI-rho-limit}\\
&& \frac{(\rho \uu)_{(0)}^{n+1} -(\rho \uu)_{(0)}^{n}}{\dt} + \nabla \cdot \left( (\rho \uu)_{(0)}^{n} \otimes (\rho \uu)_{(0)}^{n} \right) + \nabla p_{(1)}^{n+1} = 0, \label{eqn.SI-mom-limit}\\
&& \frac{(\rho e)_{(0)}^{n+1} -(\rho e)_{(0)}^{n}}{\dt} + \nabla \cdot \left(h_{(0)}^n \rho_{(0)}^{n+1} \uu_{(0)}^{n+1} \right) = 0, \label{eqn.SI-e-limit} \\
&& \nabla p_{(0)}^{n+1} = 0. \label{eqn.p-limit}
\end{eqnarray}
The incompressibility constraint 
\begin{equation}
\nabla \cdot \uu_{(0)}^{n+1}=\mathcal{O}(\dt),
\label{eqn.divV}
\end{equation}
with $\mathcal{O}(\dt)$ independent of $\varepsilon$, must now be retrieved in order to demonstrate the asymptotic preserving property of the scheme. Indeed, this recovers the limit of the energy equation at the continuous level \eqref{eq:cns_divu_lim}. To ease the notation the subscript $(0)$ is removed and all terms estimated by $C \dt$ with the constant $C \neq C(\varepsilon)$ will simply be addressed with $\mathcal{O}(\dt)$. Notice that from \eqref{eqn.SI-rho-limit}-\eqref{eqn.SI-mom-limit} one has
\begin{equation}
\rho^{n+1} = \rho^n + \mathcal{O}(\dt), \qquad \uu^{n+1} = \uu^n + \mathcal{O}(\dt).
\label{eqn.rho_u_1}
\end{equation}
Recalling that $e=e(\rho,p)$ according to \eqref{eqn.eDef}, the first term in the limit energy equation \eqref{eqn.SI-e-limit} can be written as
\begin{eqnarray}
\frac{(\rho e)^{n+1} -(\rho e)^{n}}{\dt} &=& \frac{1}{\dt} \left( \frac{\partial (\rho e)}{\partial \rho} (\rho^{n+1}-\rho^n) + \frac{\partial (\rho e)}{\partial p} (p^{n+1}-p^n) + \mathcal{O}(\rho^{n+1}-\rho^n)^2 + \mathcal{O}(p^{n+1}-p^n)^2 \right) \nonumber \\
&=& \frac{\partial (\rho e)}{\partial \rho} \frac{\rho^{n+1}-\rho^n}{\dt} + \mathcal{O}(\dt),
\label{eqn.incLim1}
\end{eqnarray}
where $p^{n+1}-p^{n}=0$ because of \eqref{eqn.p-limit} which implies that pressure is constant when $\varepsilon \to 0$. Now, using the limit continuity equation \eqref{eqn.SI-rho-limit} and the observation \eqref{eqn.rho_u_1}, from expression \eqref{eqn.incLim1} we get
\begin{eqnarray}
\frac{(\rho e)^{n+1} -(\rho e)^{n}}{\dt} &=& \frac{\partial (\rho e)}{\partial \rho} \frac{\rho^{n} - \dt \nabla(\rho^n \uu^n)-\rho^n}{\dt} + \mathcal{O}(\dt) \nonumber \\
&=& - \frac{\partial (\rho e)}{\partial \rho} \left( \uu^n \cdot \nabla \rho^n + \rho^n  \nabla \cdot \uu^n \right) + \mathcal{O}(\dt) \nonumber \\
&=& - \frac{\partial (\rho e)}{\partial \rho} \left( \uu^{n+1} \cdot \nabla \rho^n + \rho^n  \nabla \cdot \uu^{n+1} \right) + \mathcal{O}(\dt).
\label{eqn.SI-e-limit-term1}
\end{eqnarray}
Since $(\rho h)=(\rho e) + p$, the divergence flux in \eqref{eqn.SI-e-limit} is rewritten at the aid of \eqref{eqn.rho_u_1} as
\begin{eqnarray}
\nabla \cdot \left(h^n \rho^{n+1} \uu^{n+1} \right) &=& \uu^{n+1} \cdot \nabla ( h^n \rho^{n+1}) + h^n \rho^{n+1} \nabla \cdot \uu^{n+1} \nonumber \\
&=& \uu^{n+1} \cdot \nabla ( h^n \rho^{n}) + h^n \rho^{n} \nabla \cdot \uu^{n+1} + \mathcal{O}(\dt) \nonumber \\
&=& \uu^{n+1} \cdot \nabla \left( (\rho e)^n + p^n \right) + h^n \rho^{n} \nabla \cdot \uu^{n+1} + \mathcal{O}(\dt) \nonumber \\
&=& \uu^{n+1} \, \left( \frac{\partial (\rho e)}{\partial \rho} \nabla\rho^n + \frac{\partial (\rho e)}{\partial p} \nabla p^n \right) + h^n \rho^{n} \nabla \cdot \uu^{n+1} + \mathcal{O}(\dt) \nonumber \\
&=& \frac{\partial (\rho e)}{\partial \rho} \, \uu^{n+1} \cdot \nabla\rho^n + h^n \rho^{n} \nabla \cdot \uu^{n+1} + \mathcal{O}(\dt).
\label{eqn.SI-e-limit-term2}
\end{eqnarray}
The energy equation \eqref{eqn.SI-e-limit} can therefore be formulated by summing up the terms \eqref{eqn.SI-e-limit-term1}-\eqref{eqn.SI-e-limit-term2}, which yields
\begin{eqnarray}
- \frac{\partial (\rho e)}{\partial \rho} \left( \uu^{n+1} \cdot \nabla \rho^n + \rho^n  \nabla \cdot \uu^{n+1} \right) + \frac{\partial (\rho e)}{\partial \rho} \, \uu^{n+1} \cdot \nabla\rho^n + h^n \rho^{n} \nabla \cdot \uu^{n+1} &=& \mathcal{O}(\dt) \nonumber \\
- \frac{\partial (\rho e)}{\partial \rho} \rho^n  \nabla \cdot \uu^{n+1} + ( (\rho e)^n + p^n) \, \nabla \cdot \uu^{n+1} &=& \mathcal{O}(\dt) \nonumber \\
\left( -(\rho^n)^2 \frac{\partial e}{\partial \rho} + p^n \right) \, \nabla \cdot \uu^{n+1} &=& \mathcal{O}(\dt).
\label{eqn.divFree}
\end{eqnarray}
The incompressibility constraint \eqref{eqn.divV} is therefore satisfied because $\left( -(\rho^n)^2 \frac{\partial e}{\partial \rho} + p^n \right)\neq 0$. For an ideal gas EOS it holds $\frac{\partial e}{\partial \rho}=-\frac{p}{\rho^2 (\gamma-1)}$, thus equation \eqref{eqn.divFree} becomes
\begin{equation}
h^n \, \nabla \cdot \uu^{n+1} = \mathcal{O}(\dt).
\end{equation}

\subsection{First order fully discrete scheme in space and time} \label{ssec.fully-disc-1}
Let us consider a three-dimensional computational domain $\Omega(\xx)=[x_{\min};x_{\max}] \times [y_{\min};y_{\max}] \times [z_{\min};z_{\max}]$ which is discretized by a Cartesian grid composed of a total number $N_e=N_x \times N_y \times N_z$ of cells $C_{i,j,k}$ with volume $|C_{i,j,k}|=\dx \, \dy \, \dz$. Specifically, the characteristic mesh sizes are given by
\begin{equation}
\dx=\frac{x_{\max}-x_{\min}}{N_x}, \qquad \dy=\frac{y_{\max}-y_{\min}}{N_y}, \qquad  \dz=\frac{z_{\max}-z_{\min}}{N_z}.
\end{equation} 
A triple index $(i,j,k)$ referred to each space direction allows a cell to be uniquely identified. The faces in $x$, $y$ and $z$ direction are referred to as $(i+1/2,j,k)$, $(i,j+1/2,k)$ and $(i,j,k+1/2)$, respectively. The associated normal vectors are the canonical unit vectors, that is $\U{n}_{x}=(1,0,0)$, $\U{n}_{y}=(0,1,0)$, $\U{n}_{z}=(0,0,1)$. The cell center is located at point $\U{x}_{i,j,k}=(x_i,y_j,z_k)$ and a face center is at point $\U{x}_{i+1/2,j,k}=\left( \frac{x_i+x_{i+1}}{2},y_j,z_k \right)$. The spatial discretization is based on \textit{collocated grids}, in which all variables of the governing equations are defined at the cell centers of the control volumes. Implicit fluxes are discretized with finite differences with no numerical dissipation, while we rely on finite volume schemes based on numerical fluxes for the explicit terms. 

For the sake of clarity and to improve the readability, the fully discrete scheme will be presented for the one-dimensional case, since extension to multiple space dimensions on Cartesian meshes follows straightforward. Let us introduce the explicit operator $m_{i}^*$ which applies to a generic cell quantity $m_{i}^n$:
\begin{equation}
m_{i}^* = m_{i}^n 
- \Frac{\Delta t}{\Delta x} \left( f_{i+1/2}^{m} - f_{i-1/2}^{m} \right).
\label{eqn.Fop}
\end{equation}
Here, $f_{i \pm 1/2}^{m}$ denote the numerical fluxes, that are explicitly given by a Rusanov--type approximate Riemann solver, thus leading to
\begin{equation}
\begin{array}{l}
f_{i+1/2}^{m} = \Frac12 \left( f(m_{i+1}^n) + f(m_{i}^n) \right) -\Frac12 a_{i+1/2}^n \left( m_{i+1}^n - m_{i}^n \right), \qquad a_{i+1/2}^n = \max \left( |u_{i+1}^n|,|u_{i}^n| \right), \\
\\
f_{i-1/2}^{m} = \Frac12 \left( f(m_{i}^n) + f(m_{i-1}^n) \right) -\Frac12 a_{i-1/2}^n \left( m_{i}^n - m_{i-1}^n \right), \qquad a_{i-1/2}^n = \max \left( |u_{i}^n|,|u_{i-1}^n| \right),
\end{array}
\label{eqn.Fop_flux}
\end{equation}
where $f(\cdot)$ represents the physical flux related to variable $m$. Notice that the numerical viscosity is chosen to be proportional to the material speed, so that if the Mach number is high the speed of sound is bounded by the fluid velocity, whereas for very low Mach number this choice should be sufficient to guarantee stability with the dissipation $a_{i+1/2}^n\approx |u|$.

Now, a spatial discretization of the time semi-discrete scheme \eqref{eqn.SI-rho}-\eqref{eqn.SI-e} is presented. Let us start following the steps of the algorithm detailed in \ref{ssec.semi-disc}.
\begin{enumerate}
	\item The new density $\rho^{n+1}$ is immediately available solving the continuity equation \eqref{eqn.SI-rho} with $m_i^n\equiv \rho_i^n$ in \eqref{eqn.Fop}:
	\begin{equation}\label{eqn.rho}
	\rho_i^{n+1} = \left(\rho_i^n\right)^*. 
	\end{equation}
	\item The next step requires the solution of the elliptic equation for the pressure \eqref{eqn.p-SI}. Two different flux derivative operators need to be discretized, namely $\frac{\partial p}{\partial x}$ and $\frac{\partial}{\partial x} \left(h \frac{\partial p}{\partial x} \right)$. A central finite difference discretization is then given by
	\begin{eqnarray}
	\left. \frac{\partial p}{\partial x} \right|_{x_i}^{n+1} &=& \frac{p_{i+1}^{n+1}-p_{i-1}^{n+1}}{2\dx} + \mathcal{O}(\dx^2), \label{eqn.px2} \\
	\left. \frac{\partial}{\partial x} \left( h \frac{\partial p}{\partial x} \right) \right|_{x_i}^{n,n+1} &=& \frac{1}{\dx^2}\left[ \begin{array}{ccc}
	h_{i-1}^n & h_i^n & h_{i+1}^n
	\end{array} \right] \, \left[ \begin{array}{ccc}
	{3/4} & -1 & {1/4} \\ 0 & 0 & 0 \\ {1/4} & -1 & {3/4}
	\end{array} \right] \, \left[ \begin{array}{c}
	p_{i-1}^{n+1} \\ p_i^{n+1} \\ p_{i+1}^{n+1}
	\end{array} \right] + \mathcal{O}(\dx^2), \label{eqn.hpx2}
	\end{eqnarray}
    which provides up to second order of accuracy. The approximate derivative in \eqref{eqn.hpx2} is based on a finite difference approach with Lagrange interpolation polynomials of degree $2$ on the stencil composed by cells $[i-1,i,i+1]$. This makes the scheme very compact even on collocated grids, hence achieving the same properties typically related to the usage of staggered meshes. Similar discretization is applied to all derivative operators that are discretized implicitly, namely the pressure gradient in the momentum equation, the implicit part of the kinetic energy and the term $\nabla \cdot (h^n (\rho u)^*)$ in the energy equation. The pressure wave equation \eqref{eqn.p-SI} can now be expressed after multiplication by the cell volume $\dx$ as
    \begin{eqnarray}
    &&\varepsilon \dx \,  \rho_i^{n+1} e_i^{n+1} + \varepsilon \,  \frac{\dt}{4} \frac{(\rho u)_i^n}{\rho_i^{n}} \, \left(p_{i+1}^{n+1}-p_{i-1}^{n+1}\right) - \nonumber \\
    &&\frac{\dt^2}{\dx} \left( p_{i-1}^{n+1} \left( \frac{3}{4} h_{i-1}^n + \frac{1}{4} h_{i+1}^n \right) - p_{i}^{n+1} \left( h_{i-1}^n + h_{i+1}^n \right) + p_{i+1}^{n+1} \left( \frac{1}{4} h_{i-1}^n + \frac{3}{4} h_{i+1}^n \right) \right) = \varepsilon \dx \, b_i^n, 
    \label{eqn.p_SI-discr}
    \end{eqnarray}
    with the known right hand side
    \begin{equation}
    b_i^n = (\rho E)_i^{*} - \varepsilon \frac{\dt}{2} \frac{(\rho u)_i^n}{\rho_i^{n}} \,  (\rho u)_i^{*} - \frac{\dt}{2\dx} \left( h_{i+1}^n \, (\rho u)^{*}_{i+1} - h_{i-1}^n \, (\rho u)^{*}_{i-1} \right).
    \label{eqn.b-SI-discr}
    \end{equation}
    Depending on the equation of state, the pressure wave equation \eqref{eqn.p_SI-discr}-\eqref{eqn.b-SI-discr} involves either a linear or a nonlinear system that is solved following the ansatz given by \eqref{eqn.pLin} or \eqref{eqn.pNewton}, respectively. This allows to compute the new pressure $p_i^{n+1}$. 
    \item The momentum is updated at the next time level as follows:
    \begin{equation}
    (\rho u)_i^{n+1} = (\rho u)_i^{*} - \frac{1}{\varepsilon} \, \frac{\dt}{2\dx} \left( p_{i+1}^{n+1} - p_{i-1}^{n+1} \right).
    \label{eqn.mom-1}
    \end{equation}
    \item The total energy is then given by
    \begin{equation}
    (\rho E)_i^{n+1} = \rho_i^{n+1} e^{n+1} + \varepsilon \frac{(\rho u)_i^n}{2\rho_i^{n}} \, (\rho u)_i^{n+1},
    \label{eqn.en-1}
    \end{equation}
     where the internal energy is computed at the aid of the equation of state $e=e(\rho,p)$.
\end{enumerate}

\paragraph{Remark on the discretization of the enthalpy} The enthalpy in the energy flux is discretized explicitly in each control volume, i.e. $h_i^n$. Particular care must be taken in order to achieve preservation of constant velocity and pressure flows at the discrete level. Therefore, the enthalpy is not simply evaluated according to the definition \eqref{eqn.eDef}, but it is discretized as
\begin{equation}
h_i^n = \frac{\rho_i^n \, h_i^n}{\rho_i^{n+1}},
\label{eqn.hc}
\end{equation}  
which guarantees structure preserving properties that will be explained in Section \ref{ssec.AC}.


\paragraph{Remark on the compactness of the stencil}
A direct discretization of the total energy equation in \eqref{eqn.cns_rescaled} would lead to
\begin{eqnarray}
(\rho E)_i^{n+1} &=& (\rho E)_i^{*} - \frac{\dt}{2 \dx} \left( h_{i+1}^n (\rho u)_{i+1}^{n+1} - h_{i-1}^n (\rho u)_{i-1}^{n+1} \right) \nonumber \\
&=& (\rho E)_i^{*} - \frac{\dt}{2 \dx} \left[ h_{i+1}^n \left( (\rho u)_{i+1}^* - \frac{\dt}{2 \dx} \left(p_{i+2}^{n+1} - p_{i}^{n+1}\right) \right) - h_{i-1}^n \left( (\rho u)_{i-1}^* - \frac{\dt}{2 \dx} \left(p_{i}^{n+1} - p_{i-2}^{n+1}\right) \right)   \right] \nonumber \\
&:=& E_1, 
\label{eqn.E1}
\end{eqnarray}
where the viscous source terms have been neglected and the discretization of the momentum equation \eqref{eqn.mom-1} has been directly inserted. This is a consequence of the fully discrete approach, in which both momentum and energy equations are first discretized in space and time, and then formal substitution of the discrete momentum into the energy equation yields the scheme \eqref{eqn.E1}, as done in \cite{Dumbser_Casulli16,BDLTV2020}. It is evident that in this case the stencil is larger, therefore the pressure wave equation \eqref{eqn.p1} would involve those cells spanning the interval $[i-2,i-1,i,i+1,i+2]$. On the other hand, the total energy which results from the wave equation for the pressure \eqref{eqn.p_SI-discr} with the discrete operator \eqref{eqn.hpx2} is indeed given by
\begin{eqnarray}
(\rho E)_i^{n+1} &=&(\rho E)_i^{*} - \frac{\dt}{2\dx} \left( h_{i+1}^n \, (\rho u)^{*}_{i+1} - h_{i-1}^n \, (\rho u)^{*}_{i-1} \right) \nonumber \\
&+& \frac{\dt^2}{\dx^2} \left( p_{i-1}^{n+1} \left( \frac{3}{4} h_{i-1}^n + \frac{1}{4} h_{i+1}^n \right) - p_{i}^{n+1} \left( h_{i-1}^n + h_{i+1}^n \right) + p_{i+1}^{n+1} \left( \frac{1}{4} h_{i-1}^n + \frac{3}{4} h_{i+1}^n \right) \right) \nonumber \\
&:=&E_2,
\label{eqn.E2}
\end{eqnarray}
where we have first performed only a time discretization, then plugged the momentum into the energy equation and finally discretized in space. It is also possible to quantify a kind of correction $\mathcal{C}_E$ to retrieve equation \eqref{eqn.E1} from \eqref{eqn.E2}, that is
\begin{eqnarray}
C_E &=& E_2 - E_1 \nonumber \\
&=& \frac{\dt^2}{\dx^2} \left( h_{i-1}^n p_{i}^{n+1}-h_{i-1}^n p_{i-2}^{n+1} + h^n_{i+1} p_{i}^{n+1} - h_{i+1}^n p_{i+2}^{n+1} \right) \nonumber \\
&+&\frac{\dt^2}{\dx^2} \left[ 4 p_{i-1}^{n+1}\left(\frac{3}{4} h_{i-1}^n + \frac{1}{4} h_{i+1}^n \right) - 4 p_{i}^{n+1} (h_{i-1}^n + h_{i+1}^n) + 4 p_{i+1}^{n+1} \left( \frac{1}{4} h_{i-1}^n + \frac{3}{4} h_{i+1}^n \right) \right].
\end{eqnarray}
We point out that this discretization does not compromise the incompressibility constraint that must be preserved in the low Mach limit, as shown in the previous Section \ref{ssec.semi-disc}. Indeed, it allows to maintain the stencil more compact, thus improving the computational efficiency for parallel simulations.

\subsection{Exact preservation of pressure and velocity across a contact discontinuity} \label{ssec.AC}
As pointed out in \cite{BilletAbgrall2003}, a consistent numerical scheme should be able to preserve a constant pressure and a constant velocity field through a discontinuity in the fluid density during the time evolution of the solution. Let us consider a one-dimensional computational domain $\Omega=[-x_L;x_R]$ which is filled with an ideal fluid assigned with the following initial condition ($t=t_0$) for all control volumes $i = 1, \ldots N_x$:
\begin{equation}
\rho_i(t_0)=\rho_i^0 = \left\{\begin{array}{ll}
\rho_L & x \leq x_D \\ \rho_R & x > x_D
\end{array} \right. , \quad u_i(t_0)=u_i^0 = u_0, \quad p_i(t_0)=p_i^0 = p_0, \qquad i = 1, \ldots N_x,
\label{eqn.IniAC}
\end{equation}
with $x_D$ representing the location of the discontinuity inside the domain and $\rho_L \not= \rho_R$ being non-negative real numbers. The explicit operator \eqref{eqn.Fop} for density, momentum and energy explicitly writes
\begin{eqnarray}
(\rho_i)^* &=& \rho_i^n - \dtdx \left( f^{\rho,n}_{i+1/2} - f^{\rho,n}_{i-1/2} \right) 
= \rho_i^{n+1},  \label{eqn.rhoAC}\\
(\rho u)_i^* &=& (\rho u)_i^n - \dtdx \left( f^{q,n}_{i+1/2} - f^{q,n}_{i-1/2} \right) 
=  \rho_i^n \, u_i^n - \dtdx \, u_0 \left( f^{\rho,n}_{i+1/2} - f^{\rho,n}_{i-1/2} \right) 
= \rho_i^{n+1} \, u_0, \label{eqn.rhouAC} \\
(\rho E)_i^* &=&
\underbrace{\rho_i^n e_i^n + \rho_i^n k_i^n}_{=E_i^n} - \dtdx \, \Frac{u_0^2}{2} \left( f^{\rho,n}_{i+1/2} - f^{\rho,n}_{i-1/2} \right) \nonumber \\
&=& \rho_i^n e_i^n + \rho_i^n \, \Frac{(u_i^n)^2}{2} - \dtdx \, \Frac{(u_i^n)^2}{2} \left( f^{\rho,n}_{i+1/2} - f^{\rho,n}_{i-1/2} \right) = \rho_i^n e_i^n + \rho_i^{n+1} \, \Frac{u_0^2}{2} . \label{eqn.rhokAC}
\end{eqnarray}	

Assuming an ideal gas EOS \eqref{eqn.idealEOS} and using the enthalpy definition \eqref{eqn.hc} together with the above explicit operators, the right hand side \eqref{eqn.b-SI-discr} of the pressure system reads
\begin{eqnarray}
b_i &=&  \left(  \Frac{p_0}{\gamma-1} + \rho_i^{n+1} \, \Frac{u_0^2}{2} - \varepsilon \, \frac{\rho_i^{n} \, u_0}{2\rho_i^{n}} \, \rho_i^{n+1} \, u_0 \right) -\Frac{\Delta t}{2\Delta x} \left( \frac{\gamma \, p_0^n}{(\gamma-1)\rho_{i+1}^{n+1}}  \rho_{i+1}^{n+1} \, u_0 - \frac{\gamma \, p_0^n}{(\gamma-1)\rho_{i-1}^{n+1}}  \rho_{i-1}^{n+1} \, u_0 \right) \nonumber \\
&=&\Frac{p_0}{\gamma-1} + (1-\varepsilon) \, \rho_i^{n+1} \, \Frac{u_0^2}{2},
\label{eqn.rhspAC}
\end{eqnarray}  
while the left hand side \eqref{eqn.p_SI-discr} simply reduces to 
\begin{eqnarray}
\Frac{\varepsilon\Delta x}{\gamma-1} p_0 + \varepsilon \,  \frac{\dt}{4} \frac{(\rho u)_i^n}{\rho_i^{n}} \, \left(p_0-p_0\right) - \left[ \frac{3}{4} h_{i-1}^n + \frac{1}{4} h_{i+1}^n - h_{i-1}^n - h_{i+1}^n  + \frac{1}{4} h_{i-1}^n + \frac{3}{4} h_{i+1}^n \right] \, p_0 &=& \varepsilon \Delta x b_i, \nonumber \\
\Frac{\varepsilon\Delta x}{\gamma-1} p_0 + 0 + [0] p_0 &=& \varepsilon \Delta x b_i,
\label{linear_syst_pAC}
\end{eqnarray}
thus $p_0$ is the solution of the linear system \eqref{eqn.p_SI-discr}-\eqref{eqn.b-SI-discr} independently of $\varepsilon$ and the constant pressure field is preserved. As a consequence, the update of the momentum with \eqref{eqn.mom-1} writes
\begin{equation} \label{eqn.rhouAC2}
(\rho u)_i^{n+1} = \rho_i^{n+1} \, u_0 -\Frac{\Delta t}{2 \, \varepsilon \, \Delta x} \cdot 0, 
\end{equation}
hence the constant velocity is maintained as well. The test case RP0 in Section \ref{ssec:RP} gives numerical evidences of this property achieved by the novel semi-implicit scheme \eqref{eqn.p_SI-discr}-\eqref{eqn.mom-1}.

\subsection{Extension to high order of accuracy} \label{ssec.highOrder}
To reduce the effects of numerical dissipation the first order semi-implicit scheme detailed in Section \ref{ssec.fully-disc-1} is extended to high order of accuracy in space and time. A semi-implicit IMEX discretization \cite{BosFil2016} is adopted for achieving high order in time, while we rely on a CWENO reconstruction \cite{LPR:99} for gaining high accuracy in space. Finally, the numerical scheme for the explicit convective terms is implemented in a new \textit{quadrature-free} formulation. 

The governing equations \eqref{eqn.cns} can be cast into a compact and general form that writes
\begin{equation}
\frac{\partial \Q}{\partial t} + \nabla \cdot \mathbf{F} + \nabla p = \mathbf{0},
\label{eqn.pde}
\end{equation}
where $\Q=(\rho,\rho\uu,\rho E)$ is the vector of conserved variables and $\mathbf{F}=\mathbf{F}(\Q,\nabla \Q)$ represents the nonlinear flux tensor which includes both convective and viscous fluxes of the Navier-Stokes equations, that is
\begin{equation}
\mathbf{F} = \left( \begin{array}{c}
\rho \uu \\ \rho \uu \otimes \uu - \bss \\ \rho k \uu + h \rho \uu
- \bss \uu - \lambda \nabla T \end{array} \right).
\label{eqn.Ftot}
\end{equation}

\subsubsection{High order in time} \label{ssec.timeHO}
Following \cite{BosFil2016}, the governing PDE are written under the form of an autonomous system, that is
\begin{equation}
\frac{\partial \Q}{\partial t} = \mathcal{H}\left(\Q(t), \Q(t) \right), \qquad \forall t > t_0, \qquad \textnormal{with} \qquad \Q(t_0)= \Q_0,
\end{equation}
with the initial condition $\Q_0$ defined at time $t_0$. The function $\mathcal{H}$ represents the spatial approximation of the terms $\nabla \cdot \mathbf{F} + \nabla p$ in \eqref{eqn.pde}. An explicit treatment is assumed for the first argument of $\mathcal{H}$ denoted with $\Q_E$, whereas an implicit discretization is adopted for the second argument referred to as $\Q_I$, thus obtaining a partitioned system with $\Q=(\Q_E,\Q_I)$, hence
\begin{equation}
\left\{\begin{aligned}
\frac{\partial \Q_E}{\partial t} &=  \mathcal{H}\left(\Q_E, \Q_I \right) \\[0.7pt]
\frac{\partial \Q_I}{\partial t} &=  \mathcal{H}\left(\Q_E, \Q_I \right) 
\end{aligned} \right. ,
\label{eqn.Hauto}
\end{equation}
where the number of unknowns has been doubled. However, for a specific choice of time discretizations and for autonomous systems this duplication is indeed only apparent \cite{BosFil2016}. The Navier-Stokes equations with the flux splitting \eqref{eqn.fluxsplit} fulfill the formalism \eqref{eqn.Hauto}, i.e.
\begin{equation}
\mathcal{H}\left(\Q_E, \Q_I \right) = \left\{  \begin{array}{c}
(\rho \uu)_E \\ (\rho \uu \otimes \uu)_E + p_I - \bss_E \\ \rho_I \, (k \uu)_E + \rho_I \, (h \uu)_E 
- (\bss \uu)_E - \lambda \nabla T_I \end{array} \right. ,
\label{eqn.partSys}
\end{equation}
where $\Q_E=(\rho_E,(\rho \uu)_E,(\rho E)_E)$ and $\Q_I=(\rho_I,(\rho \uu)_I,(\rho E)_I)$. High order in time is achieved making use of implicit-explicit (IMEX) Runge-Kutta schemes \cite{PR_IMEX}, that are multi-step methods based on $s$ stages and typically represented with the double Butcher tableau:
\begin{equation}
\begin{array}{c|c}
\tilde{c} & \tilde{A} \\ \hline & \tilde{b}^\top
\end{array} \qquad
\begin{array}{c|c}
c & A \\ \hline & b^\top
\end{array},
\end{equation}
with the matrices $(\tilde{A},A) \in \R^{s \times s}$ and the vectors $(\tilde{c},c,\tilde{b},b) \in \R^s$. The tilde symbol refers to the explicit scheme and matrix $\tilde{A}=(\tilde{a}_{ij})$ is a lower triangular matrix with zero elements on the diagonal, while $A=({a}_{ij})$ is a triangular matrix which accounts for the implicit scheme, thus having non-zero elements on the diagonal. Here, we adopt IMEX schemes with $\tilde{b}=b$ and the stiffly accurate property in the implicit part, that is crucial for assuring asymptotic consistency and
accuracy of the scheme \cite{Boscarino2019}. Applying the partitioned Runge-Kutta method to \eqref{eqn.partSys} under the assumption that the system is autonomous, only one set of stage fluxes needs to be computed and the fluxes at each stage $i = 1, \ldots, s$ can be evaluated as
\begin{equation}
k_i = \mathcal{H} \left( \, \, \Q_E^n + \dt \sum \limits_{i=1}^s \tilde{a}_{ij} \, k_j, \quad \Q_I^n + \dt \sum \limits_{i=1}^s a_{ij} \, k_j \, \, \right), \qquad 1 \leq i \leq s.
\label{eqn.ki}
\end{equation}
A semi-implicit IMEX Runge-Kutta method is obtained as follows. Let us first set $\Q_E^n=\Q_I^n=\Q^n$, then the stage fluxes for $i = 1, \ldots, s$ are calculated as
\begin{subequations}
\begin{align}
\Q_E^i &= \Q_E^n + \dt \sum \limits_{j=1}^{i-1} \tilde{a}_{ij} k_j, \qquad 2 \leq i \leq s, \label{eq.QE} \\[0.5pt]
\tilde{\Q}_I^i &= \Q_E^n + \dt \sum \limits_{j=1}^{i-1} a_{ij} k_j, \qquad 2 \leq i \leq s, \label{eq.QI}  \\[0.5pt]
k_i &= \mathcal{H} \left( \Q_E^i, \tilde{\Q}_I^i + \dt \, a_{ii} \, k_i \right), \qquad 1 \leq i \leq s. \label{eq.k} 
\end{align}
\end{subequations}
Finally, the numerical solution is updated with
\begin{equation}
\Q^{n+1} = \Q^n + \dt \sum \limits_{i=1}^s b_i k_i.
\label{eqn.QRKfinal}
\end{equation}

Notice that equation \eqref{eq.k} implies an implicit step with the solution of a system for $k_i$, that corresponds to the pressure wave equation \eqref{eqn.p_SI-discr}-\eqref{eqn.b-SI-discr}. The final update of the solution \eqref{eqn.QRKfinal} is done using the implicit weights $b^\top$ that are assumed to be equal to the explicit ones $\tilde{b}^\top$. Furthermore, the stage fluxes $k_i$ in \eqref{eqn.ki} are the same for both explicit and implicit conserved vectors $\Q_E$ and $\Q_I$,  therefore the system is actually not doubled, since there is indeed only one set of numerical solution. 

The IMEX schemes used in this work have been developed in \cite{PR_IMEX,PR_IMEXHO} and are listened hereafter. Stiffly accurate schemes are addressed with SA, while SSP stands for Strong Stability Preserving methods, which perform better if shock waves or strong discontinuities appear in the flow. Each scheme is described with a triplet $(s,\tilde{s},p)$ which characterizes the number $s$ of stages of the implicit scheme, the number $\tilde{s}$ of stages of the explicit scheme and the order $p$ of the resulting scheme.
\begin{itemize}
	\item SP(1,1,1)
	
	\begin{equation}
	\begin{array}{c|c}
	0 & 0 \\ \hline & 1
	\end{array} \qquad
	\begin{array}{c|c}
	1 & 1 \\ \hline & 1
	\end{array}
	\label{eqn.IMEX1}
	\end{equation}
	
	\item SA SSP(3,3,2)
	
	\begin{equation}
	\begin{array}{c|ccc}
	0 & 0 & 0 & 0 \\ 1/2 & 1/2 & 0 & 0 \\ 1 & 1/2 & 1/2 & 0 \\ \hline & 1/3 & 1/3 & 1/3
	\end{array} \qquad
	\begin{array}{c|ccc}
	1/4 & 1/4 & 0 & 0 \\ 1/4 & 0 & 1/4 & 0 \\ 1 & 1/3 & 1/3 & 1/3 \\ \hline & 1/3 & 1/3 & 1/3
	\end{array}
	\label{eqn.IMEX2}
	\end{equation}
	
	\item SA DIRK(3,4,3)
	
	\begin{equation}
	\begin{array}{c|cccc}
	0 & 0 & 0 & 0 & 0 \\ \delta & \delta & 0 & 0 & 0 \\ 0.717933 & 1.437745 & −0.719812 & 0 & 0 \\ 1 & 0.916993 & 1/2 & − 0.416993 & 0 \\ \hline  & 0 & 1.208496 & -0.644363 & \delta
	\end{array} \qquad
	\begin{array}{c|cccc}
	\delta & \delta & 0 & 0 & 0 \\ \delta & 0 & \delta & 0 & 0 \\ 0.717933 & 0 & 0.282066 & \delta & 0  \\ 1 & 0 & 1.208496 & -0.644363 & 0.416993 \\ \hline  & 0 & 1.208496 & -0.644363 & \delta
	\end{array}
	\label{eqn.IMEX3a}
	\end{equation}
	$\delta=0.435866$
	
	\item SSP3(4,3,3)
	
	\begin{equation}
	\begin{array}{c|cccc}
	0 & 0 & 0 & 0 & 0 \\ 0 & 0 & 0 & 0 & 0 \\ 1 & 0 & 1 & 0 & 0 \\ 1/2 & 0 & 1/4 & 1/4 & 0 \\ \hline  & 0 & 1/6 & 1/6 & 2/3
	\end{array} \qquad
	\begin{array}{c|cccc}
	\alpha & \alpha & 0 & 0 & 0 \\ 0 & -\alpha & \alpha & 0 & 0 \\ 1 & 0 & 1-\alpha & \alpha & 0 \\ 1/2 & \delta & \eta & 1/2-\delta-\eta-\alpha & \alpha \\ \hline  & 0 & 1/6 & 1/6 & 2/3
	\end{array}
	\label{eqn.IMEX3b}
	\end{equation}
	$\alpha=0.241694$, $\delta=0.060424$, $\eta=0.129153$
	
\end{itemize} 

The first order scheme \eqref{eqn.IMEX1} corresponds to the implicit Euler method and is stiffly accurate and stability preserving. Both properties are also exhibited by the second order scheme \eqref{eqn.IMEX2}, while for third order accurate IMEX RK methods we use either \eqref{eqn.IMEX3a} or \eqref{eqn.IMEX3b} for low or high Mach number flows, respectively. Indeed, in the stiff limit only \eqref{eqn.IMEX3a} can be used in order to be consistent with the limit model at the discrete level \cite{Dellacherie1,Boscarino2019}.

\subsubsection{High order in space} \label{ssec.spaceHO}
To achieve high order of accuracy in space centered finite difference schemes are adopted for the treatment of the implicit terms, whereas a novel CWENO reconstruction is employed for the explicit convective and viscous terms in the governing equations \eqref{eqn.cns}.

\paragraph{Implicit terms} The spatial discretization of the flux derivative operators \eqref{eqn.px2}-\eqref{eqn.hpx2} presented in Section \ref{ssec.fully-disc-1} accounts for up to second order spatial accuracy. Standard finite differences are employed for higher order discretizations, hence yielding
\begin{eqnarray}
\left. \frac{\partial p}{\partial x} \right|_{x_i}^{n+1} &=& \frac{-p_{i+2}^{n+1} + 8 p_{i+1}^{n+1} - 8 p_{i-1}^{n+1} + p_{i-2}^{n+1}}{12\dx} + \mathcal{O}(\dx^4), \label{eqn.px4} \\
\left. \frac{\partial}{\partial x} \left( h \frac{\partial p}{\partial x} \right) \right|_{x_i}^{n,n+1} &=& \frac{1}{\dx^2}\left[ \begin{array}{ccccc}
h_{i-2}^n & h_{i-1}^n & h_i^n & h_{i+1}^n & h_{i+2}^n
\end{array} \right] \cdot \nonumber \\
& &\left[ \begin{array}{ccccc}
-25/144 & {1/3} & -{1/4} & {1/9} & -{1/48} \\ {1/6} & {5/9} & -1 & {1/3} & -{1/18} \\
0 & 0 & 0 & 0 & 0 \\
-{1/18} & {1/3} & -1 & {5/9} & {1/6} \\
-{1/48} & {1/9} & -{1/4} & {1/3} & -{25/144} \\
\end{array} \right] \, \left[ \begin{array}{c}
p_{i-2}^{n+1} \\ p_{i-1}^{n+1} \\ p_i^{n+1} \\ p_{i+1}^{n+1} \\ p_{i+2}^{n+1}
\end{array} \right] + \mathcal{O}(\dx^4). \label{eqn.hpx4}
\end{eqnarray}
Let observe that the operator \eqref{eqn.hpx4} is an extension of the second order operator \eqref{eqn.hpx2}, hence it is derived in the same way by using Lagrange interpolation polynomials of higher order for the derivative of the pressure. We recover the same discretization proposed in \cite{Boscarino2019} and the scheme maintains the compactness of the stencil, that is now bounded in the interval $[i-2,i-1,i,i+1,i+2]$. The fourth order finite difference approximation \eqref{eqn.px4} is applied to all first derivatives that appear in the elliptic equation on the pressure \eqref{eqn.p_SI-discr}-\eqref{eqn.b-SI-discr} as well as for the pressure flux in the momentum equation \eqref{eqn.mom-1}.

\paragraph{Explicit terms} High order shock-capturing finite volume methods are usually built upon a nonlinear reconstruction procedure that allows to stabilize the numerical scheme and avoid spurious oscillations in the vicinity of strong discontinuities. Here, we propose to develop a CWENO-type algorithm \cite{LPR:99,CPSV:2018} because it permits to keep a compact stencil, that is determined by the polynomial degree $M$ of the reconstruction. The spatial discretization makes use of Cartesian control volumes, so that the entire reconstruction algorithm can be performed in a reference element with dimensional splitting, that is we first obtain a high order polynomial of degree $M$ in $x$ direction, then in $y$ and finally in $z$ direction. This results in a computationally more efficient method compared to fully multidimensional reconstruction algorithms. A similar approach has been forwarded in \cite{AMR2013} for WENO reconstructions, requiring larger stencils and thus more memory consumption compared to the algorithm proposed in the following. 

The reconstruction procedure aims at generating high order polynomials $\ww(t,\xx)$, which are written using a \textit{nodal basis} of polynomials of degree $M$ defined in a reference unit interval $\mathcal{I}=[0;1]$. 
A cell is rescaled on the reference interval at the aid of the following change of coordinates:
\begin{equation}
\xi   = \xi(x,i) = \frac{1}{\Delta x} \left( x-x_{i-1/2} \right), \quad 
\eta  = \eta(y,j) = \frac{1}{\Delta y} \left( y-y_{i-1/2} \right), \quad 
\zeta = \zeta(z,k) = \frac{1}{\Delta z} \left( z-z_{i-1/2} \right).
\end{equation}
In particular, since the reconstruction will then be employed for the evaluation of explicit numerical fluxes across cell boundaries, the basis consists of $M+1$ linearly independent Lagrange interpolating polynomials of maximum degree $M$, i.e. $\left\{\psi_l\right\}_{l=1}^{M+1}$, passing through a set of $M+1$ nodal points $\left\{\xi_k\right\}_{k=1}^{M+1}$, which are assumed to be the Gauss-Lobatto nodes. The interpolation property holds by construction, that is
\begin{equation}
\psi_l(\xi_k) = \delta_{lk}, \qquad l,k=1, \ldots, M+1.
\end{equation}
In this way the reconstruction degrees of freedom automatically provide the values of the high order numerical solution for each conserved variable at the nodes. Furthermore, with our choice the degrees of freedom coincide with the Gauss-Lobatto nodes, thus very efficient quadrature-free computations can be designed for general integration over the reference interval, e.g. the numerical flux integration. The final reconstruction polynomial will then take the form
\begin{equation}
\ww(t^n,\xx) = \psi_l(\xi) \, \psi_q(\eta) \, \psi_r(\zeta) \, \hat{\ww}^n_{ijk,pqr},
\label{eqn.WENOpoly}
\end{equation}
with the unknown degrees of freedom $\hat{\ww}^n_{ijk,pqr}$ that must be determined. Einstein summation convention, implying summation over indices appearing twice, is adopted. The CWENO reconstruction is carried out in a dimension-by-dimension manner for each cell $C_{ijk}$ and the starting point is the definition of the one-dimensional reconstruction stencils. Contrarily to WENO schemes \cite{AMR2013}, here we always consider a total number of $N_s=3$ reconstruction stencils, namely one central stencil ($s=0$) for reconstructing a polynomial of degree $M$ and two fully one-sided stencils, one to the left ($s=1$) and the other one to the right ($s=2$), for obtaining second order polynomials that are only used for nonlinear stabilization in the presence of discontinuous profiles of the numerical solution. The central stencil $s=0$ is assembled for each Cartesian direction as 
\begin{equation}
\mathcal{S}^{0,x}_{ijk} = \bigcup_{e=i-L}^{i+L} C_{ejk}, \quad \mathcal{S}^{0,y}_{ijk} = \bigcup_{e=j-L}^{j+L} C_{iek}, \quad \mathcal{S}^{0,z}_{ijk} = \bigcup_{e=k-L}^{k+L} C_{ije},
\label{eqn.centralStencil}
\end{equation}
where the spatial extension of the stencil to the left $L$ and to the right $R$ is given by \\
\begin{center}
\begin{tabular}{c|cc}
$M$ & $L$ & $R$ \\
\hline 
even & 	$-M/2$ & $M/2$ \\
odd & 	$-(M+1)/2$ & $(M+1)/2$
\end{tabular}	
.
\end{center}

The low order one-sided stencil $s=1$ (left-sided) and $s=2$ (right-sided) are simply assembled with the element under consideration $C_{ijk}$ and the direct neighbor to the left and to the right, that is
\begin{equation}
\mathcal{S}^{1,x}_{ijk} = \bigcup_{e=i-1}^{i} C_{ejk}, \quad \mathcal{S}^{2,x}_{ijk} = \bigcup_{e=i}^{i+1} C_{ejk},
\label{eqn.sidedStencil}
\end{equation}
the same holding for $y$ and $z$ direction. The reconstruction is based on integral conservation of all conserved quantities stored in the state vector $\Q^n$ and is firstly performed along the $x$ direction. Therefore, we look for a reconstruction polynomial defined on each reconstruction stencil $s=0,1,2$ of the form
\begin{equation}
\ww^{s,x}(t^n,x) = \psi_p(\xi) \, \hat{\ww}_{ijk,p}^{n,s},
\label{eqn.wx}
\end{equation}
for which integral conservation holds, that is
\begin{equation}
\frac{1}{\Delta x} \int_{x_{e-1/2}}^{x_{e+1/2}} \psi_p(\xi(x)) \, \hat{\ww}_{ijk,p}^{n,s} = \Q_{ejk}^n, \qquad \forall C_{ejk} \in \mathcal{S}^{s,x}_{ijk},
\label{eqn.intCons}
\end{equation}
that must be prescribed for each stencil $s \in [0,1,2]$. Recall that for the one-sided stencils ($s=1,2$) the reconstruction polynomial is of degree one, therefore the nodal basis are defined accordingly by locally setting $M=1$. Equations \eqref{eqn.intCons} lead to a linear system which might become overdetermined in the case of even order schemes. We rely on a constrained least squares (CLSQ) technique \cite{DumbserKaeser2007,ADERFSE} for determining the unknowns  $\hat{\ww}_{ijk,p}^{n,s}$, where the linear constraint is given by requiring that integral conservation \eqref{eqn.intCons} exactly holds true for the cell $C_{ijk}$ under consideration. In the CWENO framework the polynomial $\ww^{0,x}(t^n,x)$ defined on the central stencil is often referred to as \textit{optimal polynomial}, because among all the possible polynomials of degree $M$, it is the only one that shares the same cell average $\Q_{ijk}^n$ in the element, while being close in the least-square sense to the other cell averages in the stencil. According to \cite{ADER_CWENO}, the central polynomial $\tilde{\ww}^{0,x}(t^n,x)$ is then obtained by difference between the polynomial $\ww^{0,x}(t^n,x)$ and the linear combination of the one-sided polynomials $\ww^{1,x}(t^n,x)$ and $\ww^{2,x}(t^n,x)$ of lower degree \cite{CPSV:2018}, that is
\begin{equation}
\label{CWENO:P0}
\tilde{\ww}^{0,x}(t^n,x)= \frac{1}{\lambda_0}\left(\ww^{0,x}(t^n,x) - \sum_{s=1}^{2} \lambda_{s} \ww^{s,x}(t^n,x) \right), \qquad \tilde{\ww}^{1,x}(t^n,x) = {\ww}^{1,x}(t^n,x), \quad \tilde{\ww}^{2,x}(t^n,x) = {\ww}^{2,x}(t^n,x),
\end{equation}
where $\lambda_{s}$ with $s=0,1,2$ are positive coefficients such that
\begin{equation}
\sum_{s=0}^{2}\lambda_{s}=1:=\lambda_{sum}.
\label{eqn.sumCWENO}
\end{equation}
Here, we do not use the pointwise WENO formulation originally introduced in \cite{shu_efficient_weno}, but the polynomial WENO schemes forwarded in \cite{DumbserKaeser2007}. As a consequence, the linear weights are a normalization which sums up to unity and we set $\lambda_{0}=200/\lambda_{sum}$ for $\mathcal{S}^{0,x}$ and $\lambda_{1}=\lambda_{2}=1/\lambda_{sum}$ for all one-sided polynomials \cite{ADER_CWENO}. Once the polynomials $\ww^{s,x}(t^n,x)$ in \eqref{eqn.wx} are available, we proceed by constructing a nonlinear data-dependent hybridization among the three polynomials obtained for each stencil, that is
\begin{equation}
\label{CWENOx}
\tilde{\ww}^{x}(t^n,x) = \psi_p(\xi) \, \hat{\ww}^n_{ijk,p}, \qquad \textnormal{with} \qquad \hat{\ww}^n_{ijk,p} = \sum_{s=0}^{2} \omega_s \tilde{\ww}^{s,x}(t^n,x),
\end{equation}
where the nonlinear weights $\omega_s$ are given by 
\begin{equation}
\label{eqn.weights}
\omega_s = \frac{\tilde{\omega}_s}{\sum \limits_{s=0}^{2} \tilde{\omega}_s}, 
\qquad \textnormal{ with } \qquad 
\tilde{\omega}_s = \frac{\lambda_s}{\left(\sigma_s + \epsilon \right)^r}. 
\end{equation} 
The parameter $\epsilon=10^{-14}$ avoids division by zero and the exponent $r=4$ is chosen according to \cite{ADER_CWENO}. The oscillation indicators $\sigma_s$ are given by
\begin{equation}
\sigma_s = \Sigma_{lm} \hat{\ww}_{l}^{n,s} \hat{\ww}_{m}^{n,s},
\label{eqn.OI}
\end{equation}
where the oscillation matrix $\Sigma_{lm}$ can be computed as done in \cite{DumbserEnauxToro} once and for all on the reference interval $\mathcal{I}$, hence
\begin{equation}
\Sigma_{lm} = \sum_{\alpha=1}^{M} \int_{0}^{1} \frac{\partial^\alpha \psi_l(\xi)}{\partial \xi^\alpha} \cdot \frac{\partial^\alpha \psi_m(\xi)}{\partial \xi^\alpha} \, d\xi. 
\end{equation}
Notice that the integrals appearing in \eqref{eqn.intCons}, which then constitute the so-called \textit{reconstruction matrix}, only depend on the geometry, i.e. on the interval over which integration is carried out. Since this corresponds to the reference element $\mathcal{I}$, the reconstruction matrix can be evaluated, inverted and stored during the pre-processing stage and it remains the same throughout the entire computation. Furthermore, only one reconstruction matrix is needed because all control volumes are rescaled to the reference element and our CWENO reconstruction is carried out one by one for each spatial dimension.

The polynomials $\ww^{x}(t^n,x)$ obtained so far are high order accurate in $x$ direction, but they still remain a cell average along the $y$ and $z$ direction. Therefore the CWENO reconstruction procedure illustrated above needs to be performed again along $y$ and finally along $z$ direction (see \cite{AMR2013} for further details). The final element-wise reconstruction polynomials $\ww(t^n,\xx)$ in \eqref{eqn.WENOpoly} represent entire polynomials defined by a nodal basis, which makes use of the Gauss-Lobatto interpolation points. Consequently, the degrees of freedom associated to the high order reconstruction are nothing but the high order extrapolated values of the conserved quantities at a set of quadrature nodes, thus they are ready for performing integration as a direct result of the reconstruction procedure. In other words, no further reconstruction evaluations will be needed while integrating over cell boundaries for the numerical flux computation.  

\paragraph{Quadrature-free finite volume scheme for the explicit fluxes} The interpolation property of the CWENO reconstruction polynomials, that are expressed in terms of a nodal basis defined through a set of Gauss-Lobatto points in the reference interval $\mathcal{I}$, can be fully exploited for designing a quadrature-free finite volume solver for the computation of the explicit fluxes in \eqref{eqn.Fop}. For the sake of clarity let us consider a one-dimensional setting with the generic cell quantity $m_i^n$. The integral of the basis functions $\psi_l(\xi)$ over the reference interval $\mathcal{I}$ is simply given by
\begin{equation}
\mathcal{F}_l:=\int_{0}^{1} \psi_l(\xi) \, d\xi, \qquad l=1,\ldots,M+1,
\end{equation}
which will then be used as \textit{universal flux matrix}. The high order version of the finite volume scheme \eqref{eqn.Fop} writes
\begin{equation}
m_{i}^* = m_{i}^n 
- \Frac{\Delta t}{\Delta x} \left( \mathcal{F}_l \hat{f}_{l,i+1/2}^{m} - \mathcal{F}_q \hat{f}_{q,i-1/2}^{m} \right),
\label{eqn.QFFV}
\end{equation}
where the expansion coefficients of the fluxes, i.e. $\hat{f}_{l,i+1/2}^{m}$ and $\hat{f}_{q,i-1/2}^{m}$, are obtained by computing the corresponding fluxes defined in \eqref{eqn.Fop_flux} at the Gauss-Lobatto nodes, that is $l$ and $q$ for the face $i+1/2$ and $i-1/2$, respectively. Because the reconstruction values are directly available at quadrature points, no computation is needed but it is sufficient to pick the correct degree of freedom (either $l$ or $q$ in \eqref{eqn.QFFV}) out of the CWENO polynomial $\ww(t^n,\xx)$ and calculate the Rusanov flux \eqref{eqn.Fop_flux}. From \eqref{eqn.QFFV} it is evident that the high order finite volume scheme is quadrature-free, thus requiring only a matrix-vector multiplication for obtaining high order numerical fluxes for all conserved variables. The same applies for the fluxes in $y$ and $z$ direction.

\subsection{Viscous fluxes of the Navier-Stokes equations} \label{ssec.viscFlux}
The extension of the algorithm to the Navier-Stokes model is simply performed relying on an explicit discretization of the viscous fluxes which are then added to the explicit operators $m_i^*$. For first order schemes the discrete velocity gradients on
the control volume boundaries are computed as
\begin{equation}
\nabla u_{i+1/2}^n = \frac{1}{\Delta x} \left( \nabla u_{i+1}^n - \nabla u_{i}^n \right), \qquad \nabla u_{i-1/2}^n = \frac{1}{\Delta x} \left( \nabla u_{i}^n - \nabla u_{i-1}^n \right).
\label{eqn.gradU1}
\end{equation}
For higher order ($M>0$) the CWENO reconstruction polynomials are exploited, since they automatically provide gradients of the conserved variables. Specifically, the Rusanov flux \eqref{eqn.Fop_flux} is slightly modified along the lines of \cite{GassnerDiffusion,ADERNSE} in order to include both the convective and the viscous terms, hence obtaining the following numerical flux $f^{\ww}$ for the reconstructed conserved variables $\ww(t^n,x)$ at the interface $x_{i+1/2}$:
\begin{eqnarray}
f^{\ww}_{i+1/2} &=& \Frac12 \left[ f\left( \ww_{i+1}(t^n,x_{i+1/2}),\nabla \ww_{i+1}(t^n,x_{i+1/2}) \right) + f \left(\ww_{i}(t^n,x_{i+1/2}),\nabla \ww_{i}(t^n,x_{i+1/2}) \right) \right] \nonumber \\
&-&\Frac12 \left( a_{i+1/2}^n + 2 \eta \lambda_{v,i+1/2}^n \right) \cdot \left[ \ww_{i+1}(t^n,x_{i+1/2}) - \ww_{i}(t^n,x_{i+1/2}) \right],
\label{eqn.Fop_fluxHO}
\end{eqnarray} 
with $\lambda_{v,i+1/2}^n $ representing the maximum eigenvalue of the viscous operator defined in \eqref{eqn.lambdav}. Here, the physical fluxes $f(\cdot)$ contain all terms of the nonlinear flux tensor $\mathbf{F}$ \eqref{eqn.Ftot} and the numerical viscosity is supplemented with a dissipative coefficient $\eta$ which arises from the solution of the generalized diffusive Riemann problem \cite{GassnerDiffusion} and is evaluated as
\begin{equation}
\eta = \frac{2 M + 1}{\Delta x \sqrt{\pi/2}}.
\end{equation} 

\subsection{\textit{A posteriori} stabilization at high Mach number} \label{ssec.limiter}
In the case of strong discontinuities in the flow, implicit high order time discretizations are not able to remove overshooting and undershooting of
the numerical solution. This aspect has been investigated in \cite{Gottlieb2001} and numerically observed in \cite{BDLTV2020} for a second order IMEX scheme applied to the Euler equations. Spurious oscillations are generated by the violation of the explicit CFL stability condition and they do not vanish but remain limited in time. In this sense high order implicit schemes are not $L_\infty$ stable but only $L_2$ stable.

To overcome this problem, we propose to use a stabilization technique that is based on a \textit{convex combination} of high order and first order schemes, which are proven to ensure monotonicity. Differently from the \textit{a posteriori} approach that has been developed in \cite{DLDV2018}, here we employ a limiting procedure which makes use of an \textit{a priori} strategy. The first objective is to detect troubled cells, i.e. those regions of the computational domain $\Omega$ which are characterized by strong shocks. We rely on the flattener variable described in \cite{BalsaraFlattener} as \textit{shock indicator}. A shock can be identified by comparing the divergence of the velocity field $\nabla \cdot \uu^n$ with the minimum of the sound speed $c_{\min}^n$ obtained by considering the element $C_{ijk}$ itself as well as its Neumann neighborhood $\mathcal{D}_{ijk}$, i.e. all elements which share one face with $C_{ijk}$:
\begin{equation}
\mathcal{D}_{ijk} = \bigcup_{e=i-1}^{e=i+1} C_{ejk} \cup \bigcup_{e=j-1}^{e=j+1} C_{iek} \cup \bigcup_{e=k-1}^{e=k+1} C_{ije} \qquad \textnormal{with} \quad e\neq \{i,j,k\}.
\end{equation}
The divergence of the velocity field is then evaluated as follows: 
\begin{eqnarray}
(\nabla \cdot \uu^n)_{ijk} &=& \frac{1}{\Delta x}  \left[ (u^n_{i+1jk} - u^n_{ijk}) - (u^n_{i-1jk} - u^n_{ijk})\right] \nonumber \\
&+& \frac{1}{\Delta y}  \left[ (v^n_{ij+1k} - v^n_{ijk}) - (v^n_{ij-1k} - v^n_{ijk})\right]\nonumber \\
&+& \frac{1}{\Delta z}  \left[ (w^n_{ijk+1} - w^n_{ijk}) - (w^n_{ijk-1} - w^n_{ijk}) \right]. 
\label{eqn.divVflattener}
\end{eqnarray}
Among all neighbors we compute the minimum sound speed $c_{\min}^n \in \mathcal{D}_{ijk}$, which is a function of the pressure and the density. The divergence of the velocity field \eqref{eqn.divVflattener} is estimated from the cell-averaged states $\Q^n$ which are known at the current time. The flattener variable $\chi_{ijk}^n$ can now be computed:
\begin{equation}
\chi^n = \min {\left[ 1, \max {\left(0, -\frac{\nabla \cdot \uu^n + k_1 c_{\min}^n}{k_1 c_{\min}^n}\right)}\right]},
\label{eqn.flattener}
\end{equation}
with the coefficient $k_1=10^{-3}$ set for all our computations. To ensure further stabilization, the flattener is extended also to those elements which are about to be crossed by a shock, but have still to enter the wave, as done in \cite{BoscheriWAO}. The flattener variable is interpreted as a \textit{detector}, therefore the cell is flagged as troubled if $\chi^n>0$. Let observe that in the case of rarefaction waves, where the divergence of the velocity field is positive in \eqref{eqn.flattener}, and when shocks of modest strength occur, that is $-k_1 c_{\min}^n \leq \nabla \cdot \uu^n \leq 0$, the flattener variable is zero. Moreover, the flattener is bounded in the interval $[0;1]$.

Once the flattener indicator has been computed for all cells, the semi-implicit scheme presented in Section \ref{ssec.fully-disc-1} is run with high order time and space discretizations following the algorithm detailed in Section \ref{ssec.highOrder}. As a result one obtains a so-called \textit{candidate solution} $\Q^{n+1,\mathcal{O}(M+1)}$ that is of order $M+1$. Then, if at least one cell is marked as troubled by the flattener, a first order numerical solution is computed, i.e. $\Q^{n+1,\mathcal{O}(1)}$. Finally, the new solution at the next time level is given by the convex combination 
\begin{equation}
\Q^{n+1} = \chi^n \, \Q^{n+1,\mathcal{O}(1)} + (1-\chi^n) \, \Q^{n+1,\mathcal{O}(M+1)}.
\label{eqn.Q1flattener}
\end{equation}
If no cells are marked as troubled, then the new solution corresponds to the fully high order candidate solution. We underline that very few cells are typically flagged by the flattener and if no shocks occur, like in the low Mach regime, the flattener is never activated and the semi-implicit scheme is always run without any further stabilization. 


\section{Numerical results} \label{sec:num_results}
The new high order semi-implicit pressure solver (SI-P) is applied to a large set of different test cases in order to asses the accuracy and the robustness of the numerical scheme. Firstly, 
the accuracy of the method is validated at different Mach number regimes. Secondly, shock tube problems with ideal gas law \eqref{eqn.idealEOS} and Redlich-Kwong EOS \eqref{eqn.pEOS}-\eqref{eqn.eEOS} are considered, thus showing the capability of the high order semi-implicit method to deal with both linear and nonlinear equation of state. Finally, multidimensional test cases for inviscid and viscous flows involving shocks and other discontinuities are presented. All simulations are run in a fully three-dimensional setting and the time step is always computed according to a CFL-type stability condition that is only based on the maximum absolute value of the flow velocity and eventually the viscous eigenvalues, i.e.
\begin{equation}
\Delta t\leq  \textnormal{CFL} \min_\Omega \left(\Frac{|u|}{\Delta x}+\Frac{|v|}{\Delta y}+\Frac{|w|}{\Delta z}+ \left(\Frac{\lambda_v}{\Delta x^2}+\Frac{\lambda_v}{\Delta y^2}+\Frac{\lambda_v}{\Delta z^2}\right)\right)^{-1},
\label{eqn.dtsi}
\end{equation}
which does no longer involve any dependency on the rescaled sound speed $c/\sqrt{\varepsilon}$ compared to the time step \eqref{eqn.dtex} of fully explicit discretizations. The computational domain is addressed with $\Omega$ and is discretized with a total number of $N_x \times N_y \times N_z$ Cartesian control volumes. If not specified, the ideal gas EOS is assumed with $\gamma=1.4$, the flattener variable presented in Section \ref{ssec.limiter} is not activated and the third order version of the method in space and time is adopted. For viscous flows we assume constant viscosity, hence we set $\beta=1$ and $s=0$ in \eqref{eqn.Sutherland}. The vector of conserved variables is $\Q=(\rho, \rho u, \rho v, \rho w, \rho E)$, while the vector of primitive variables is addressed with $\mathbf{U}=(\rho,u,v,w,p)$.

\subsection{Numerical convergence studies} \label{ssec:conv}
The convergence of the novel semi-implicit pressure solver presented in this article is studied by considering a modified version of the smooth isentropic vortex \cite{HuShuTri} governed by the compressible Euler equations, thus we set $\mu=\lambda=0$ in the Navier-Stokes system \eqref{eqn.cns}. The computational domain is given by $\Omega=[0;10] \times [0;10] \times [0;1]$ with periodic boundaries. The fluid is characterized by a homogeneous background field on the top of which some perturbations are added, thus
\begin{equation}
\mathbf{U}(t_0,\xx) = (1+\delta \rho, 1, 1, 0, 1+\delta p),
\label{eq.ConvEul-IC}
\end{equation}
with the perturbations for temperature $\delta T$, density $\delta \rho$ and pressure $\delta p$ that read
\begin{equation}
\label{ShuVortDelta}
\delta T = -\frac{(\gamma-1)\epsilon^2}{8\gamma\pi^2}e^{1-r^2}, \quad
\delta \rho = (1+\delta T)^{\frac{1}{\gamma-1}}-1,  \quad 
\delta p = (1+\delta T)^{\frac{\gamma}{\gamma-1}}-1.
\end{equation}
The vortex maintains perfect equilibrium and the flow is stationary, thus the exact solution $\mathbf{U}_{ex}$ is simply given by the initial condition at any time $t>0$, i.e. $\mathbf{U}_{ex} = \mathbf{U}(t_0,\xx)$. The final time of the simulation is $t_f=1$ and the test is run on a sequence of successively refined computational meshes. The grids are refined in the $x-y$ plane while keeping constant the number of cells $N_z=4$ along the $z$ direction. The error $L_m$ is normalized with respect to the exact solution, hence it is computed at the final time as
\begin{equation}
L_m(\Q) = \frac{\sqrt[m]{\int_\Omega \left|\left|  \ww(t_f,\xx) - \Q(t_0,\xx)\right|\right|^m \, d\xx}}{\sqrt[m]{\int_\Omega \left|\left|  \Q(t_0,\xx) \right|\right|^m \, d\xx}},
\end{equation} 
where the integrals are evaluated with Gaussian quadrature formulae of suitable order of accuracy (see \cite{stroud}) and the exponent $m$ determines the type of error norm that is computed. The numerical solution $\Q(t_f,\xx)$ is reconstructed with the high order accurate CWENO procedure detailed in Section \ref{ssec.highOrder}, that is $\ww(t_f,\xx)$. The time step is computed according to \eqref{eqn.dtsi} with $\textnormal{CFL}=0.9$. Numerical convergence studies are firstly carried out in the normal regime $\varepsilon=1$ for second and third order accurate semi-implicit schemes. The results are reported in Table \ref{tab.conv} where errors are measured in $L_1$, $L_2$ and $L_\infty$ norm for the conserved variables $(\rho,\rho u, \rho E)$. Both IMEX schemes \eqref{eqn.IMEX3a} and \eqref{eqn.IMEX3b} are proven to achieve the formal order of accuracy as well as the second order SA-SSP2 scheme \eqref{eqn.IMEX2}. 

\begin{table}[!htbp]  
	\caption{Numerical convergence results for the compressible Euler equations using second and third order SI-P schemes with $\varepsilon=1$ and different IMEX time stepping. The errors are measured in $L_1$ norm and refer to the variables $\rho$ (density), $\rho u$ (horizontal momentum) and $\rho E$ (energy) at time $t=1$.}  
	\begin{center} 
		\begin{small}
			\renewcommand{\arraystretch}{1.0}
			\begin{tabular}{c|cccccc} 
				\multicolumn{7}{c}{SI-P $\mathcal{O}2$} \\
				\hline
				$N_x$ ($N_y$) & $L_1(\rho)$ & $\mathcal{O}(\rho)$ & $L_1(\rho u)$ & $\mathcal{O}(\rho u)$ &  $L_1(\rho E)$ & $\mathcal{O}(\rho E)$ \\ 
				\hline
				16  & 8.274E-04 & -    & 1.861E-02 & -    & 8.897E-04 & -     \\ 
				32  & 3.374E-04 & 1.29 & 4.307E-03 & 2.11 & 3.041E-04 & 1.55  \\ 
				64  & 8.838E-05 & 1.93 & 1.045E-03 & 2.04 & 8.468E-05 & 1.84  \\ 
				128 & 2.199E-05 & 2.01 & 2.583E-04 & 2.02 & 2.195E-05 & 1.95  \\ 
				\hline
				\multicolumn{7}{c}{} \\
				\multicolumn{7}{c}{SI-P $\mathcal{O}3$ with SA-DIRK(3,4,3)} \\
				\hline
				$N_x$ ($N_y$) & $L_1(\rho)$ & $\mathcal{O}(\rho)$ & $L_1(\rho u)$ & $\mathcal{O}(\rho u)$ &  $L_1(\rho E)$ & $\mathcal{O}(\rho E)$ \\ 
				\hline
				16  & 7.677E-05 & -    & 2.134E-03 & -    & 9.448E-05 & -     \\ 
				32  & 1.164E-05 & 2.72 & 2.884E-04 & 2.89 & 1.296E-05 & 2.87  \\ 
				64  & 1.508E-06 & 2.95 & 3.898E-05 & 2.89 & 1.806E-06 & 2.84  \\ 
				128 & 1.741E-07 & 3.11 & 6.320E-06 & 2.62 & 8.737E-08 & 4.37  \\ 
				\hline
				\multicolumn{7}{c}{} \\
				\multicolumn{7}{c}{SI-P $\mathcal{O}3$ with SSP(4,3,3)} \\
				\hline
				$N_x$ ($N_y$) & $L_1(\rho)$ & $\mathcal{O}(\rho)$ & $L_1(\rho u)$ & $\mathcal{O}(\rho u)$ &  $L_1(\rho E)$ & $\mathcal{O}(\rho E)$ \\ 
				\hline
				16  & 7.428E-05 & -    & 2.137E-03 & -    & 9.495E-05 & -     \\ 
				32  & 1.128E-05 & 2.72 & 2.880E-04 & 2.89 & 1.292E-05 & 2.88  \\ 
				64  & 1.498E-06 & 2.91 & 3.900E-05 & 2.88 & 1.806E-06 & 2.84  \\ 
				128 & 2.346E-07 & 2.67 & 5.887E-06 & 2.73 & 2.856E-07 & 2.66  \\ 
				\hline
			\end{tabular}
		\end{small}
	\end{center}
	\label{tab.conv}
\end{table}

\begin{figure}[!htbp]
	\begin{center}
		\begin{tabular}{ccc} 
			\includegraphics[width=0.33\textwidth]{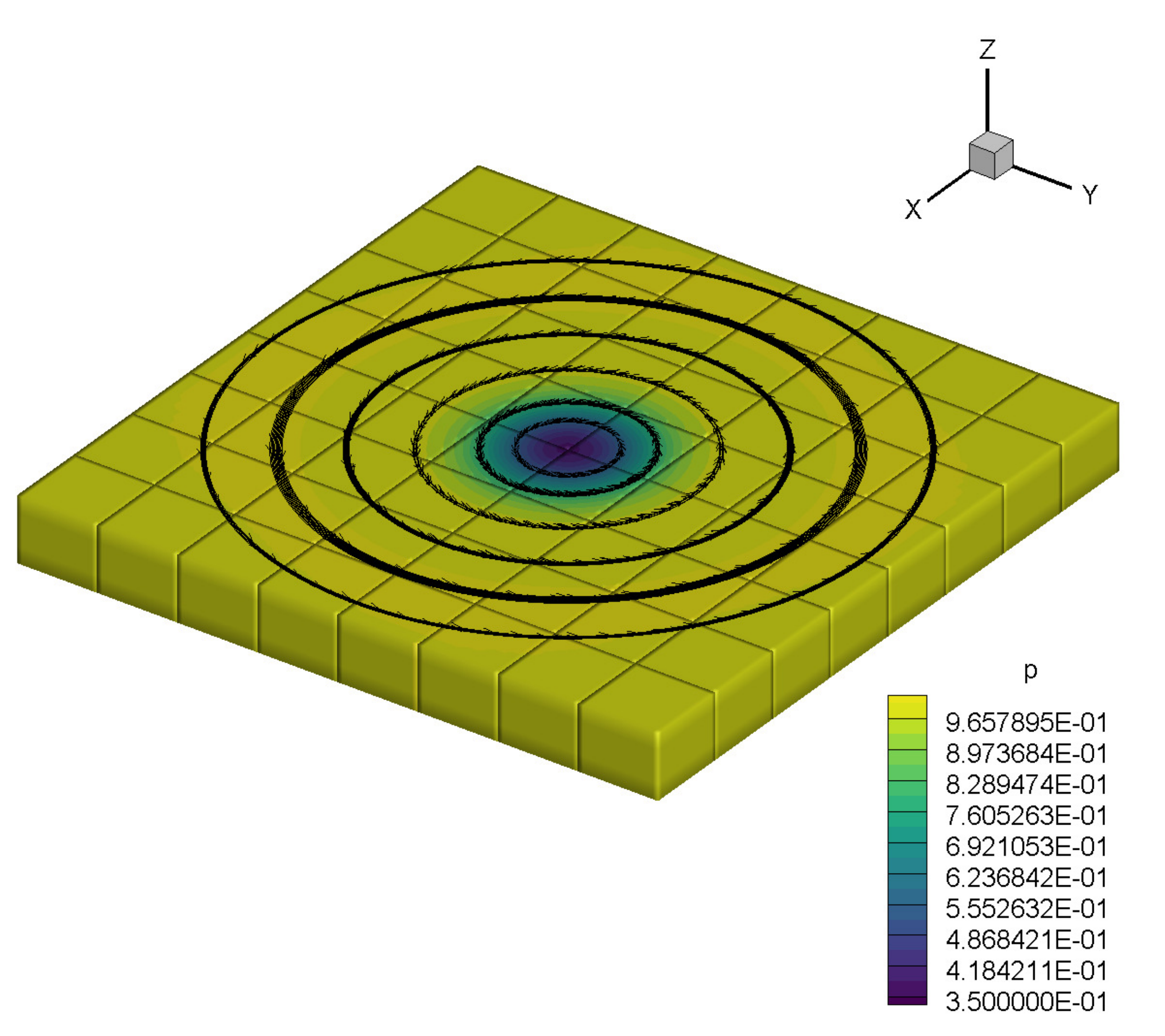} & 
			\includegraphics[width=0.33\textwidth]{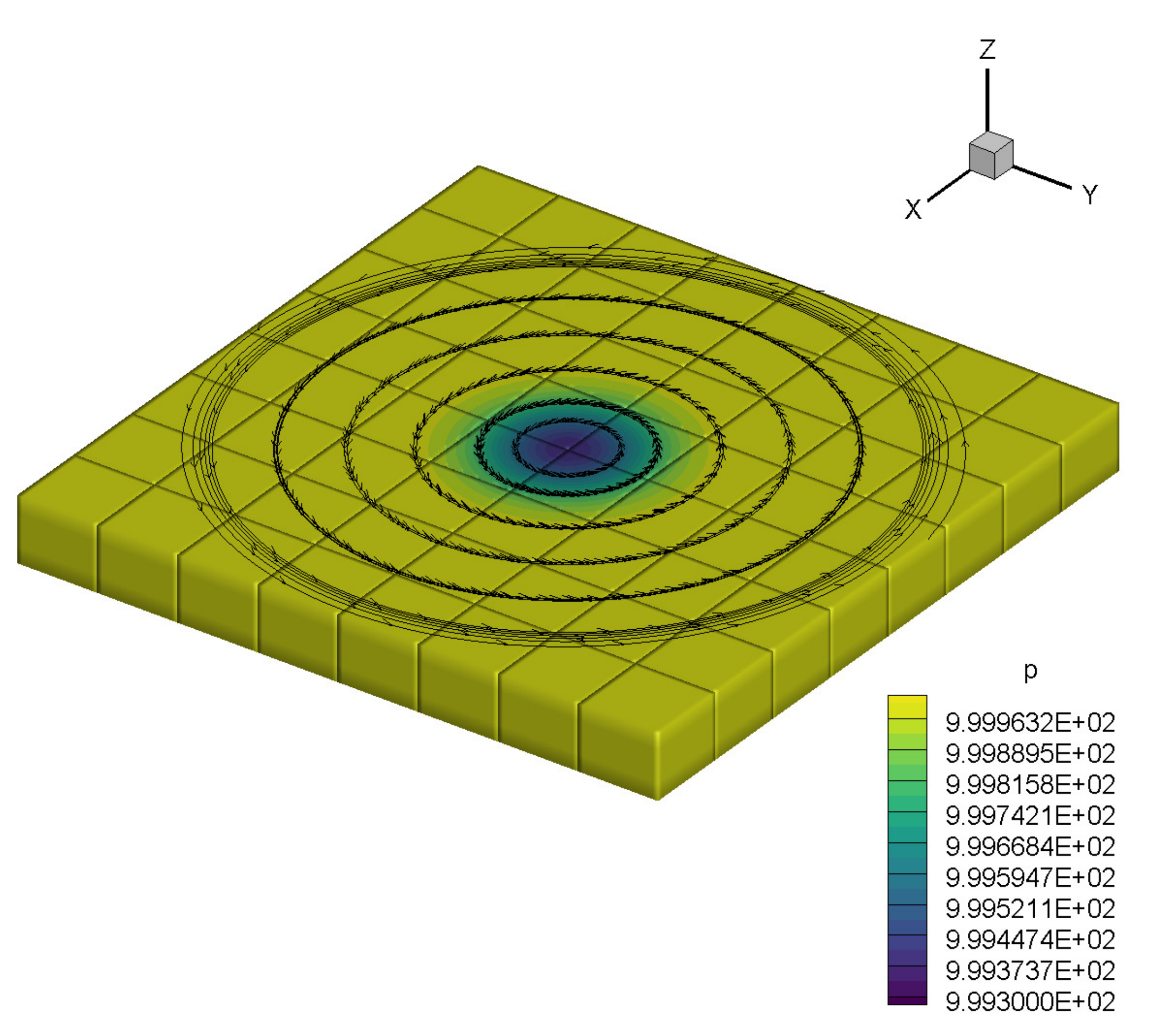} &
			\includegraphics[width=0.33\textwidth]{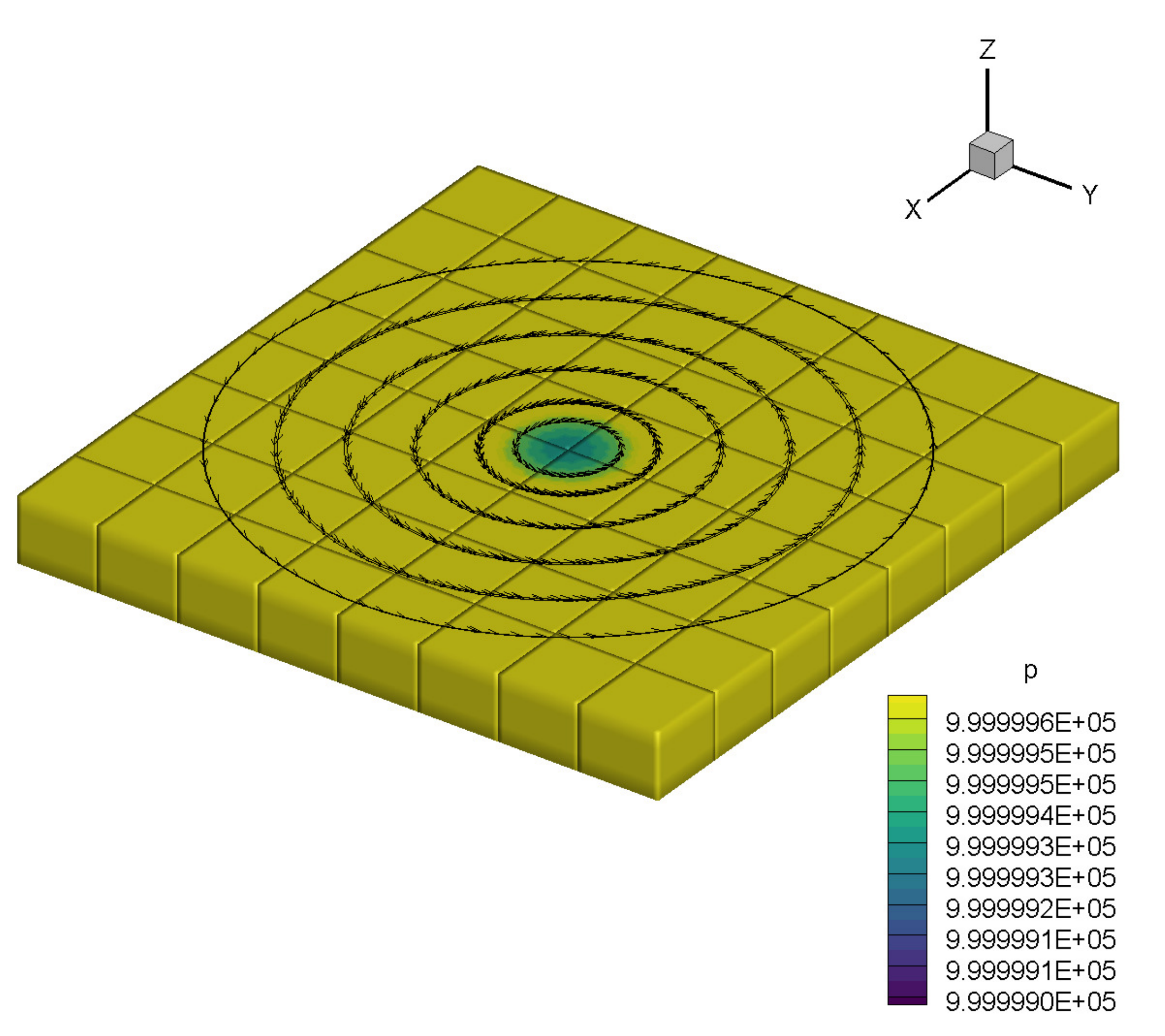}\\   
		\end{tabular} 
		\caption{Third order pressure contours ($20$ levels have been used in the range  bounded by the minimum and maximum value of the pressure) for the isentropic vortex test with $N_x=N_y=128$ at time $t=1$. Stream-traces of the velocity field for $\varepsilon=10^{0}$ (left), $\varepsilon=10^{-3}$ (middle) and $\varepsilon=10^{-6}$ (right).}
		\label{fig.ShuVortex3D}
	\end{center}
\end{figure}

Secondly, the behavior of the scheme at low Mach regimes is investigated by considering different values of the stiffness parameter $\varepsilon$, namely we consider $\varepsilon \in [10^{-6};10^{-1}]$ and the convergence rates are shown for the horizontal momentum $\rho u$ in Table \eqref{tab.conv_lowMach}. Second and third order of accuracy are well preserved even in the limit $\varepsilon=10^{-6}$ in which pressure is almost constant and the total energy is entirely constituted by its kinetic part. Figure \ref{fig.ShuVortex3D} depicts the pressure contours and the velocity stream-traces for the smooth isentropic vortex in the low Mach regime, highlighting that the numerical solution is independent of the Mach number, as expected. The stiffly accurate IMEX scheme \eqref{eqn.IMEX3a} has been used at third order for retrieving the correct asymptotic behavior in the stiff limit.

\begin{table}[!htbp]  
	\caption{Numerical convergence results for the compressible Euler equations using second and third order SI-P schemes at different low Mach regimes with stiffness parameters ranging in the interval $[\varepsilon=10^{-1};\varepsilon=10^{-6}]$. The errors are measured in $L_1$ norm and refer to the variables $\rho u$ (horizontal momentum) at time $t=1$.}  
	\begin{center} 
		\begin{small}
			\renewcommand{\arraystretch}{1.0}
			\begin{tabular}{c|cccccc} 
				\multicolumn{7}{c}{SI-P $\mathcal{O}2$} \\
				\hline
				\hline
				& \multicolumn{2}{c}{$\varepsilon=10^{-1}$} & \multicolumn{2}{c}{$\varepsilon=10^{-2}$} & \multicolumn{2}{c}{$\varepsilon=10^{-3}$} \\
				$N_x$ ($N_y$) & ${L_1}$ & Order & ${L_1}$ & Order& ${L_1}$ & Order \\ 
				\hline
				16  & 1.791E-02 & -    & 1.731E-02 & -    & 1.744E-02 & -    \\ 
				32  & 3.558E-03 & 2.33 & 3.875E-03 & 2.16 & 3.442E-03 & 2.34 \\ 
				64  & 7.835E-04 & 2.18 & 9.854E-04 & 1.98 & 7.485E-04 & 2.20 \\ 
				128 & 1.856E-04 & 2.08 & 2.116E-04 & 2.22 & 1.811E-04 & 2.05 \\ 
				\multicolumn{7}{c}{} \\
				& \multicolumn{2}{c}{$\varepsilon=10^{-4}$} & \multicolumn{2}{c}{$\varepsilon=10^{-5}$} & \multicolumn{2}{c}{$\varepsilon=10^{-6}$} \\
				$N_x$ ($N_y$) & ${L_1}$ & Order & ${L_1}$ & Order& ${L_1}$ & Order \\ 
				\hline
				16  & 1.743E-02 & -    & 1.738E-02 & -    & 1.738E-02 & -    \\ 
				32  & 3.473E-03 & 2.33 & 3.474E-03 & 2.32 & 3.473E-03 & 2.32 \\ 
				64  & 7.303E-04 & 2.25 & 7.327E-04 & 2.25 & 7.333E-04 & 2.24 \\ 
				128 & 1.663E-04 & 2.13 & 1.668E-04 & 2.14 & 1.669E-04 & 2.14 \\
				\hline
				\multicolumn{7}{c}{} \\
				\multicolumn{7}{c}{SI-P $\mathcal{O}3$ with SA-DIRK(3,4,3)} \\
				\hline
				\hline
				& \multicolumn{2}{c}{$\varepsilon=10^{-1}$} & \multicolumn{2}{c}{$\varepsilon=10^{-2}$} & \multicolumn{2}{c}{$\varepsilon=10^{-3}$} \\
				$N_x$ ($N_y$) & ${L_1}$ & Order & ${L_1}$ & Order& ${L_1}$ & Order \\ 
				\hline
				16  & 2.222E-03 & -    & 2.227E-03 & -    & 2.188E-03 & -    \\ 
				32  & 3.381E-04 & 2.72 & 3.319E-04 & 2.75 & 3.292E-04 & 2.73 \\ 
				64  & 4.608E-05 & 2.88 & 4.440E-05 & 2.90 & 4.350E-05 & 2.92 \\ 
				128 & 6.320E-06 & 2.87 & 5.746E-06 & 2.95 & 5.511E-06 & 2.98 \\ 
				\multicolumn{7}{c}{} \\
				& \multicolumn{2}{c}{$\varepsilon=10^{-4}$} & \multicolumn{2}{c}{$\varepsilon=10^{-5}$} & \multicolumn{2}{c}{$\varepsilon=10^{-6}$} \\
				$N_x$ ($N_y$) & ${L_1}$ & Order & ${L_1}$ & Order& ${L_1}$ & Order \\ 
				\hline
				16  & 2.188E-03 & -    & 2.188E-03 & -    & 2.325E-03 & -    \\ 
				32  & 3.280E-04 & 2.74 & 3.328E-04 & 2.74 & 3.279E-04 & 2.83 \\ 
				64  & 4.342E-05 & 2.92 & 4.342E-05 & 2.92 & 4.341E-05 & 2.92 \\ 
				128 & 5.448E-06 & 2.99 & 5.448E-06 & 2.99 & 5.448E-06 & 2.99 \\ 
				\hline
			\end{tabular}
		\end{small}
	\end{center}
	\label{tab.conv_lowMach}
\end{table}


\subsection{Shock tube problems} \label{ssec:RP}
The novel numerical method is here validated against a set of well-known Riemann problems for the compressible Euler equations taken from \cite{ToroBook}. The initial condition of the gas consists in a left (L) and a right (R) state that are separated by a discontinuity located at $x=x_d$. The computational domain is the box $\Omega=[xL;xR] \times [0;0.1]\times [0;0.1]$ with Dirichlet boundary conditions imposed along the $x$ direction and periodic boundaries set elsewhere. Table \ref{tab:init} summarizes the extension of the computational domain as well as the initial condition for density, horizontal velocity and pressure for all shock tube problems considered in the following. Riemann problems RK1 and RK2 are concerned with the nonlinear Redlich-Kwong EOS and the computational domain is discretized with $N_x \times N_y \times N_z = 400 \times 4 \times 4$ control volumes, while the other tests involve an ideal gas and the computational grid is composed of $N_x \times N_y \times N_z = 200 \times 4 \times 4$ cells. The computation stops at the final time indicated in Table \ref{tab:init} and we set $\textnormal{CFL}=0.9$ for the first four test cases, whereas $\textnormal{CFL}=0.5$ is adopted for the simulations involving the Redlich-Kwong EOS. The reference solution for all test problems is computed with an explicit second order MUSCL-TVD scheme run on a very fine mesh composed of 10'000 cells. The numerical solution is plot considering a 1D cut through the $x$ direction of the computational domain with $200$ equidistant sample points.

\begin{table}[!htbp]  
	\caption{Initialization of shock tube problems. Initial states left (L) and right (R) are reported as well as the final time of the simulation $t_f$, the computational domain $[x_L;x_R]$ and the position of the initial discontinuity $x_d$. The equation of state (EOS) is also specified. }  
	\begin{center} 
		\begin{small}
			\renewcommand{\arraystretch}{1.0}
			\numerikNine
			\begin{tabular}{l|c|ccc|ccc|ccc|c} 
				\hline
				Name & $t_{f}$ & $x_L$ & $x_R$ & $x_d$ & $\rho_L$ & $u_L$ & $p_L$ & $\rho_R$ & $u_R$ & $p_R$ & EOS \\
				\hline
				RP0 (Contact)          & 0.50 & 0.0 & 1.0 & 0.25 & 1000  & 0     & 10$^5$ & 0.01 & 0.0 & 10$^5$ & ideal gas \\
				RP1 (Lax)              & 0.14 & 0.0 & 1.0 & 0.50 & 0.445 & 1.698 & 3.528  & 0.5 &  0.0 & 0.571  & ideal gas \\
				RP2 (Two shocks)       & 0.80 & 0.0 & 1.0 & 0.50 & 1.0   &  2.0  & 0.1    & 1.0 & -2.0 & 0.1    & ideal gas \\
				RP3 (Two rarefactions) & 0.15 & 0.0 & 1.0 & 0.50 & 1.0   & -1.0  & 0.4    & 1.0 &  1.0 & 0.4    & ideal gas \\
				RK1                    & 0.10  & -0.5 & 0.5 & 0.0 & 1.0 & 1.0 & 2.0 & 1.0   & -1.0 & 1.0 & Redlich-Kwong \\
				RK2                    & 0.20  & -0.5 & 0.5 & 0.0 & 1.0 & 0.0 & 1.0 & 0.125 & 0.0  & 0.1 & Redlich-Kwong \\
				\hline
			\end{tabular}
		\end{small}
	\end{center}
	\label{tab:init}
\end{table}

\begin{figure}[!htbp]
	\begin{center}
		\begin{tabular}{cc} 
			\includegraphics[width=0.47\textwidth]{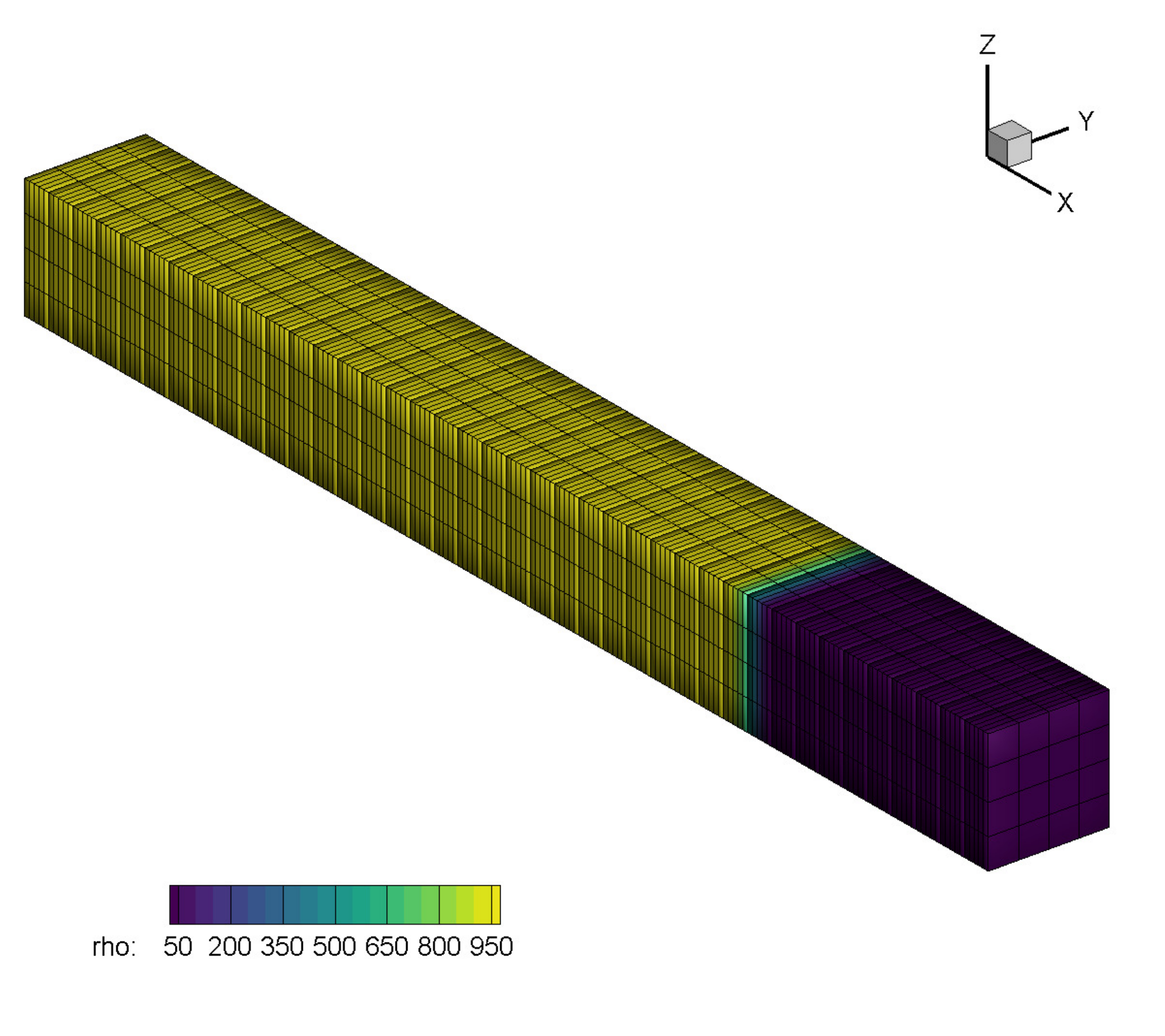}  & 
			\includegraphics[width=0.47\textwidth]{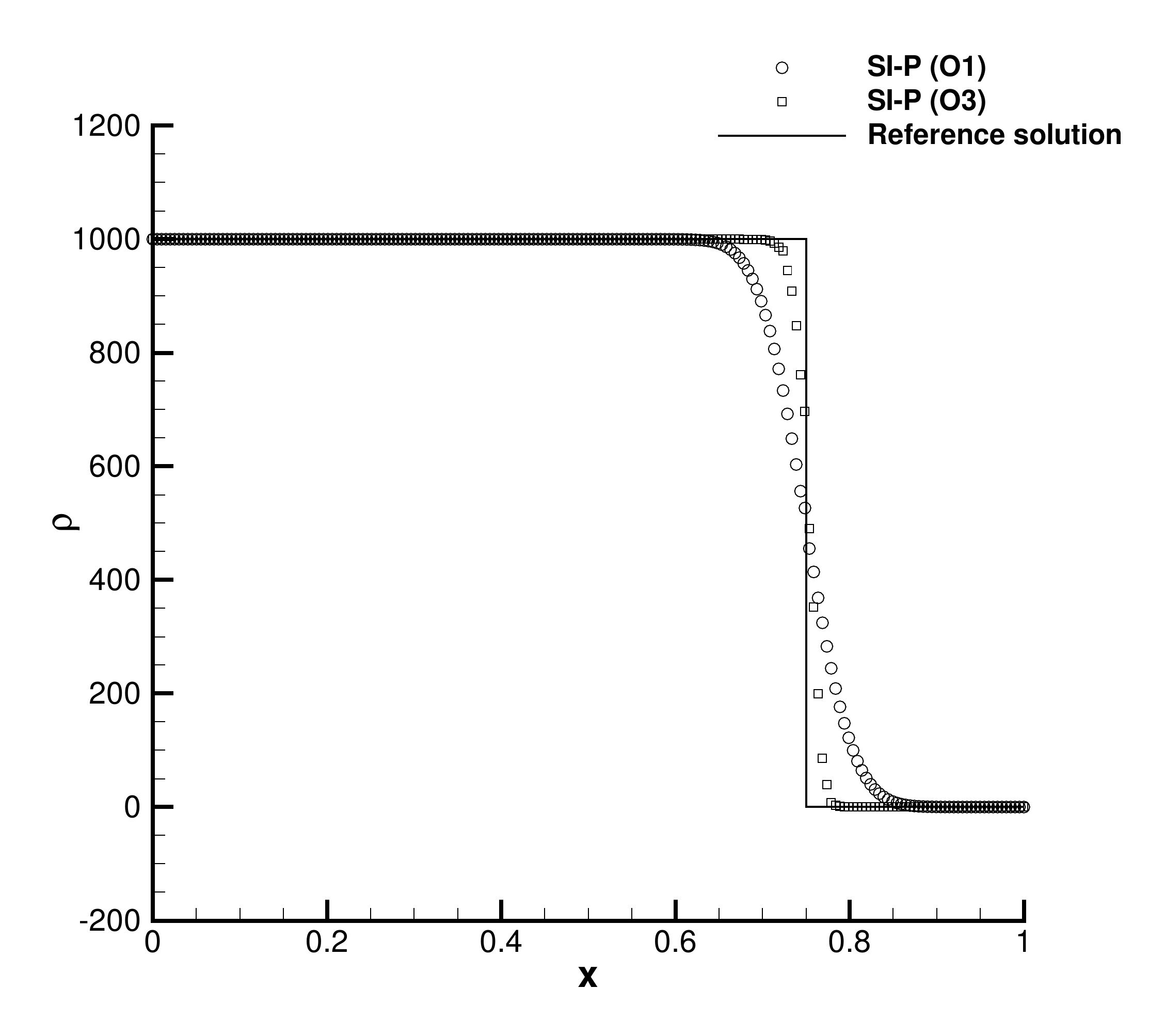} \\  
			\includegraphics[width=0.47\textwidth]{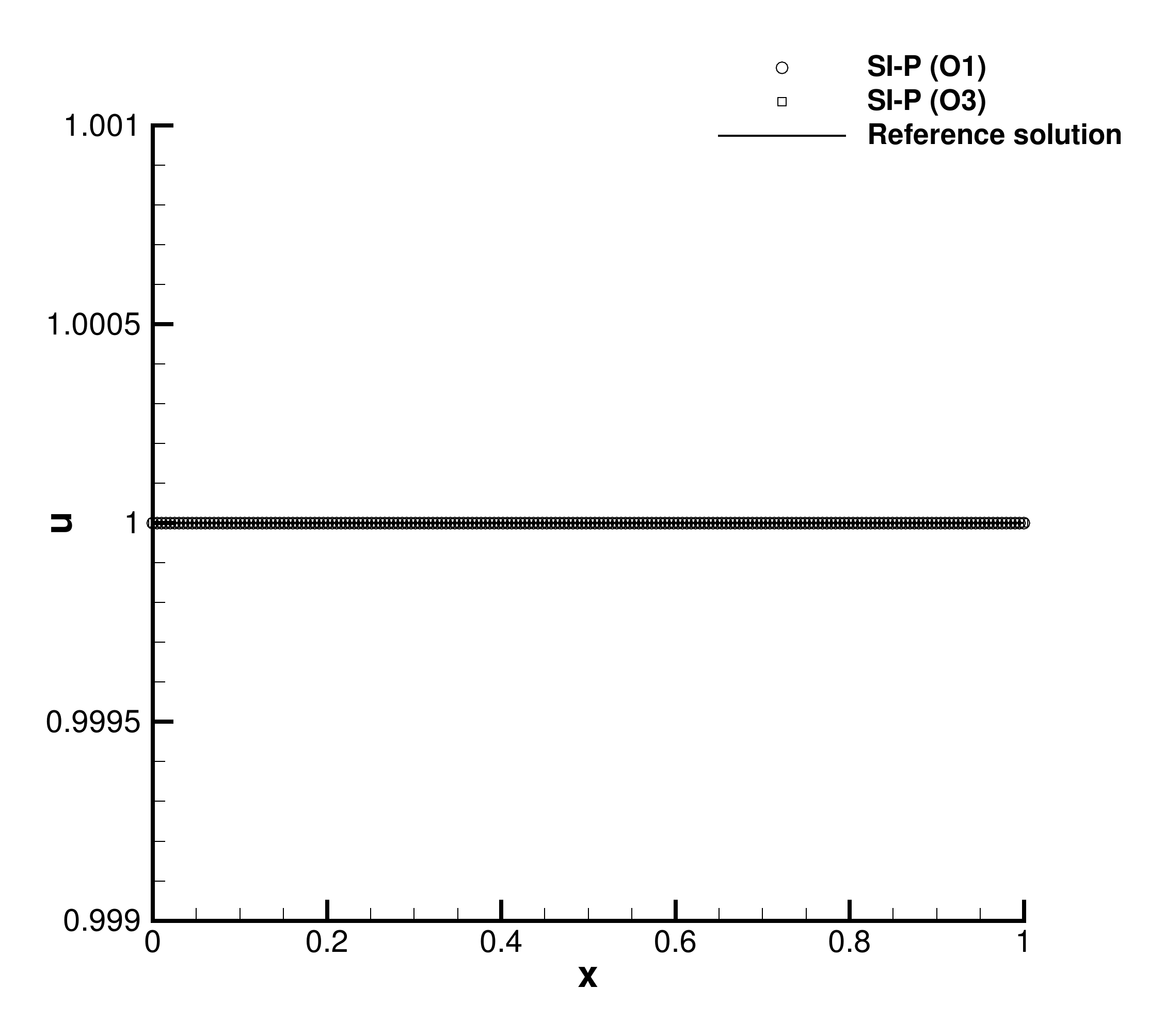} &  
			\includegraphics[width=0.47\textwidth]{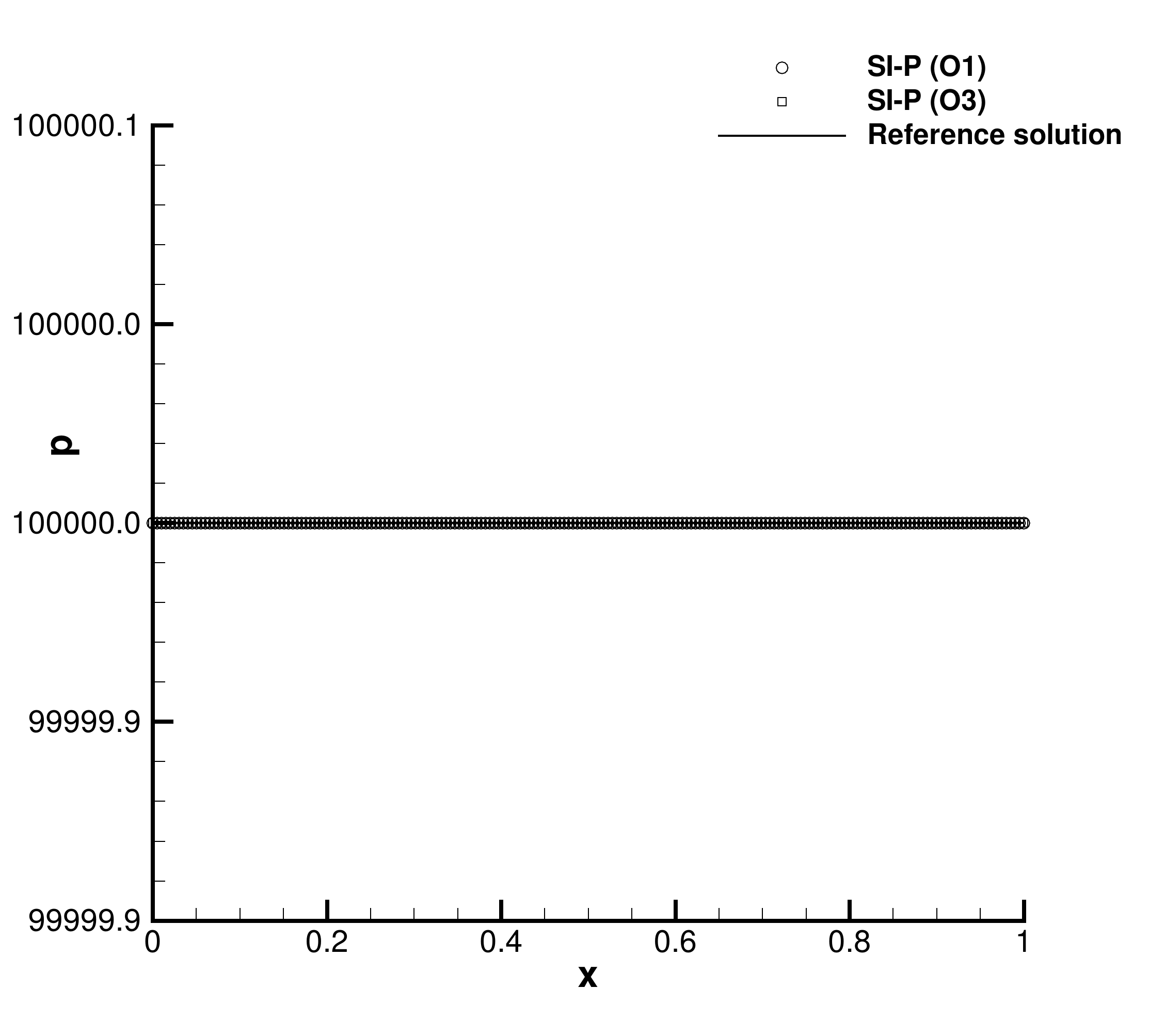}   
		\end{tabular} 
		\caption{RP0 involving a moving contact discontinuity at final time $t_f=0.5$. From top left to bottom right: 3d computational grid with density contours and comparison of density, velocity and pressure (symbols) versus the reference solution (straight line).}
		\label{fig.RP0}
	\end{center}
\end{figure}

Test RP0 provides numerical evidences about the property of the SI-P scheme of maintaining an exact preservation of constant pressure and velocity across a contact discontinuity, see Section \ref{ssec.AC}. Indeed, a contact discontinuity involving a density step of five orders of magnitude is moving at constant velocity and pressure to the right of the domain at Mach number ranging in the interval $[2.7\cdot 10^{-4};8.5\cdot 10^{-2}]$. Figure \ref{fig.RP0} shows a comparison between first and third order numerical solution and it highlights that pressure and horizontal velocity are kept constant up to machine precision. 

\begin{figure}[!htbp]
	\begin{center}
		\begin{tabular}{cc} 
			\includegraphics[width=0.47\textwidth]{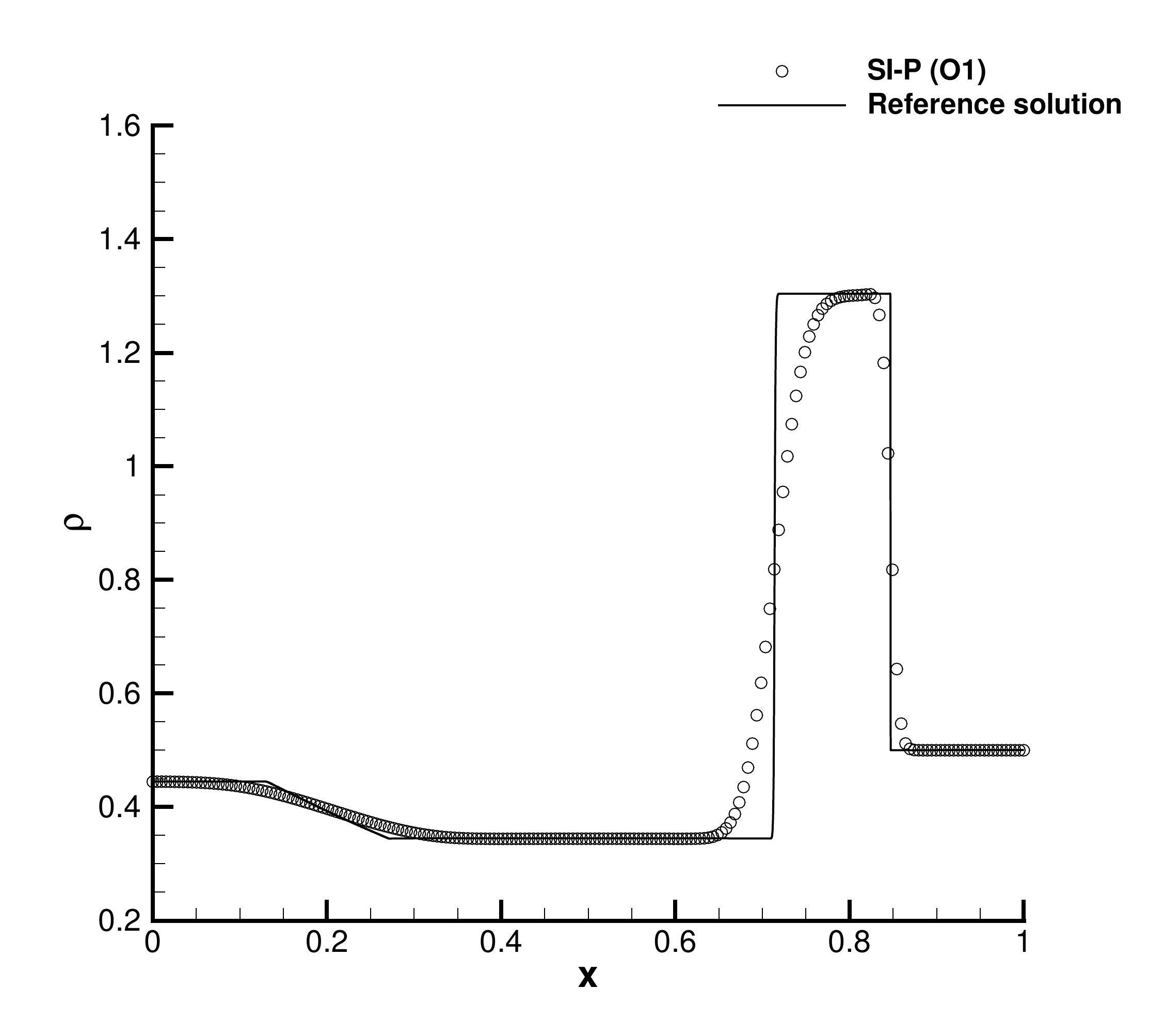}  & 
			\includegraphics[width=0.47\textwidth]{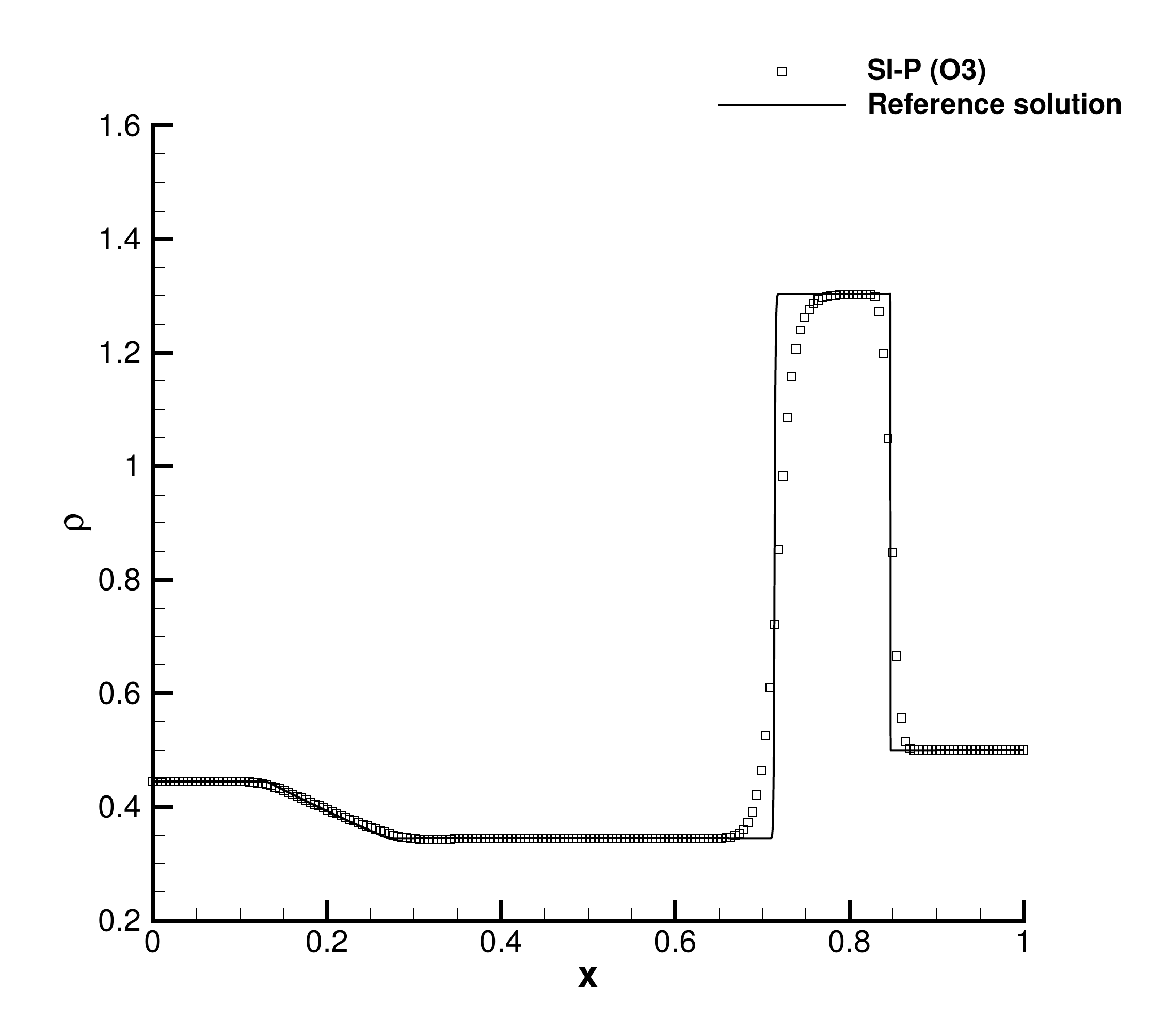} \\  
			\includegraphics[width=0.47\textwidth]{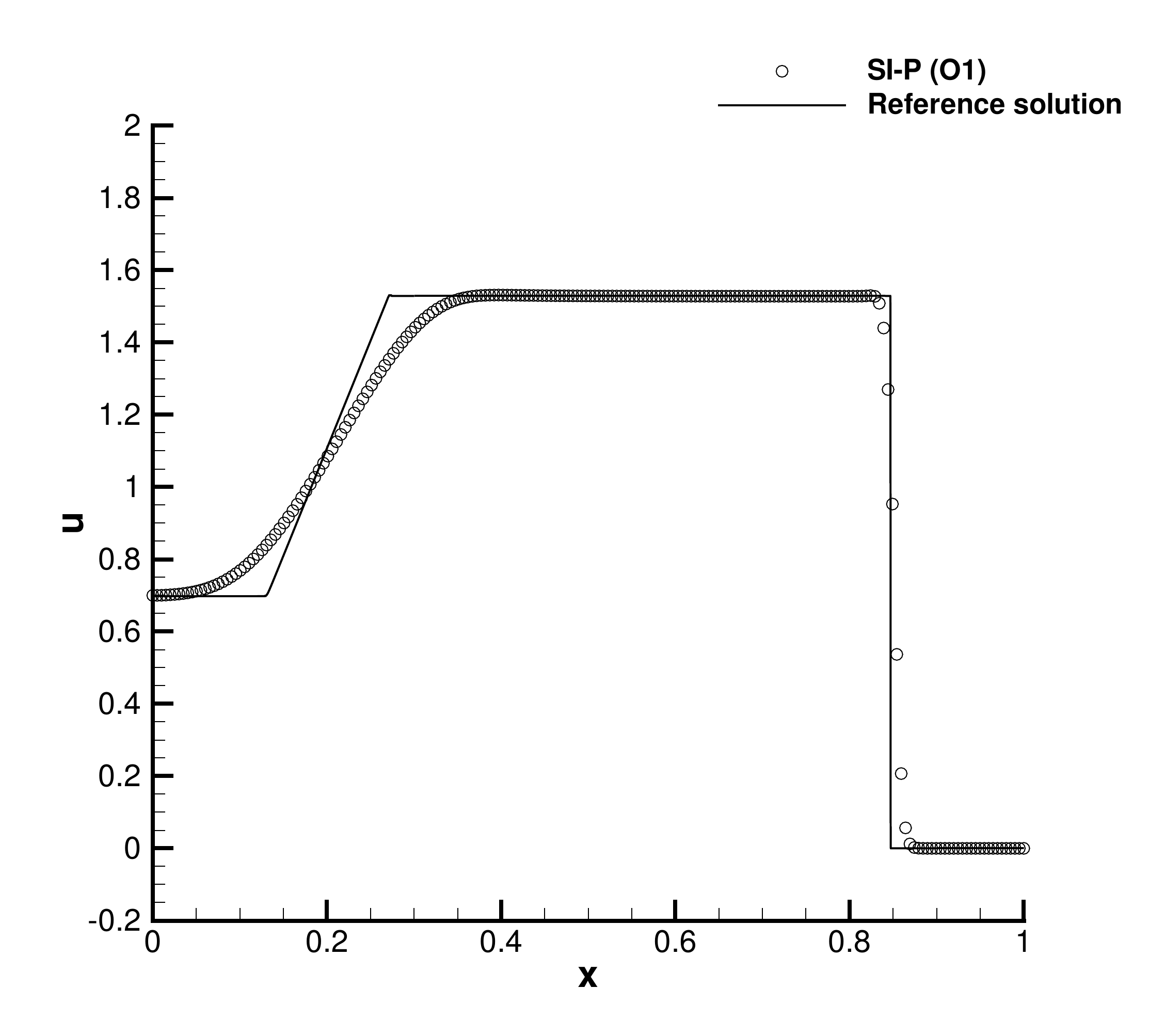}  & 
			\includegraphics[width=0.47\textwidth]{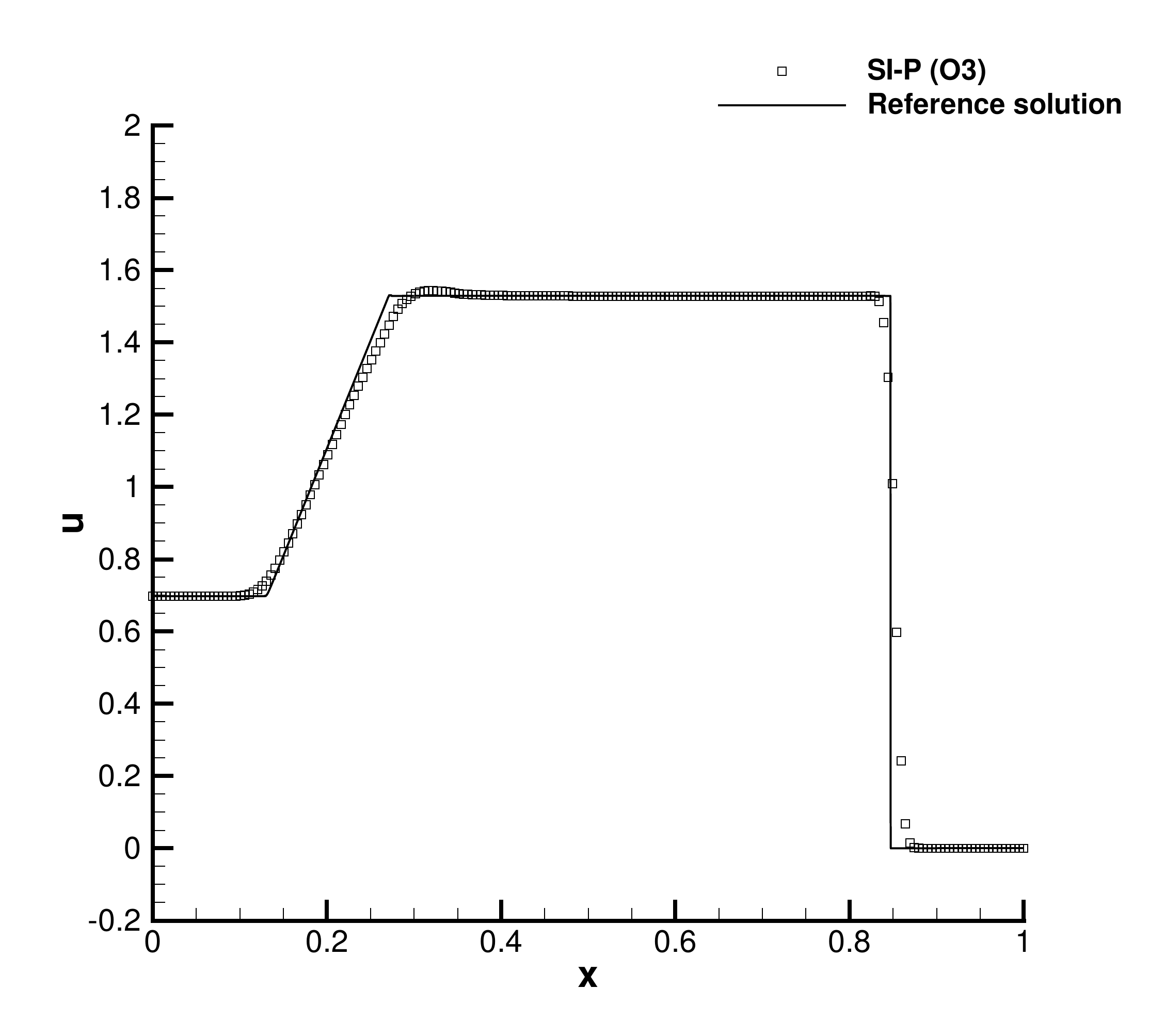} \\  
			\includegraphics[width=0.47\textwidth]{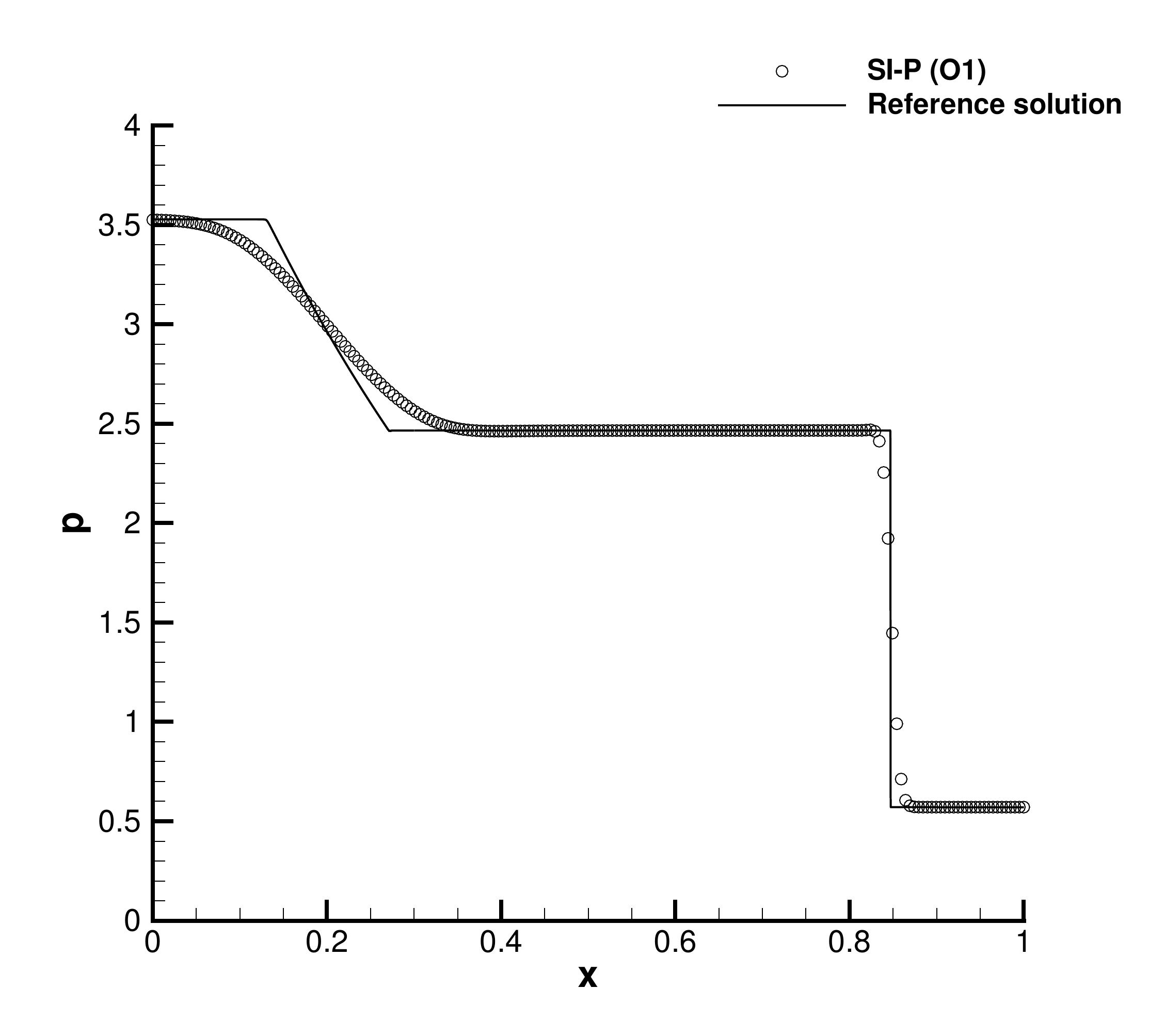}  & 
			\includegraphics[width=0.47\textwidth]{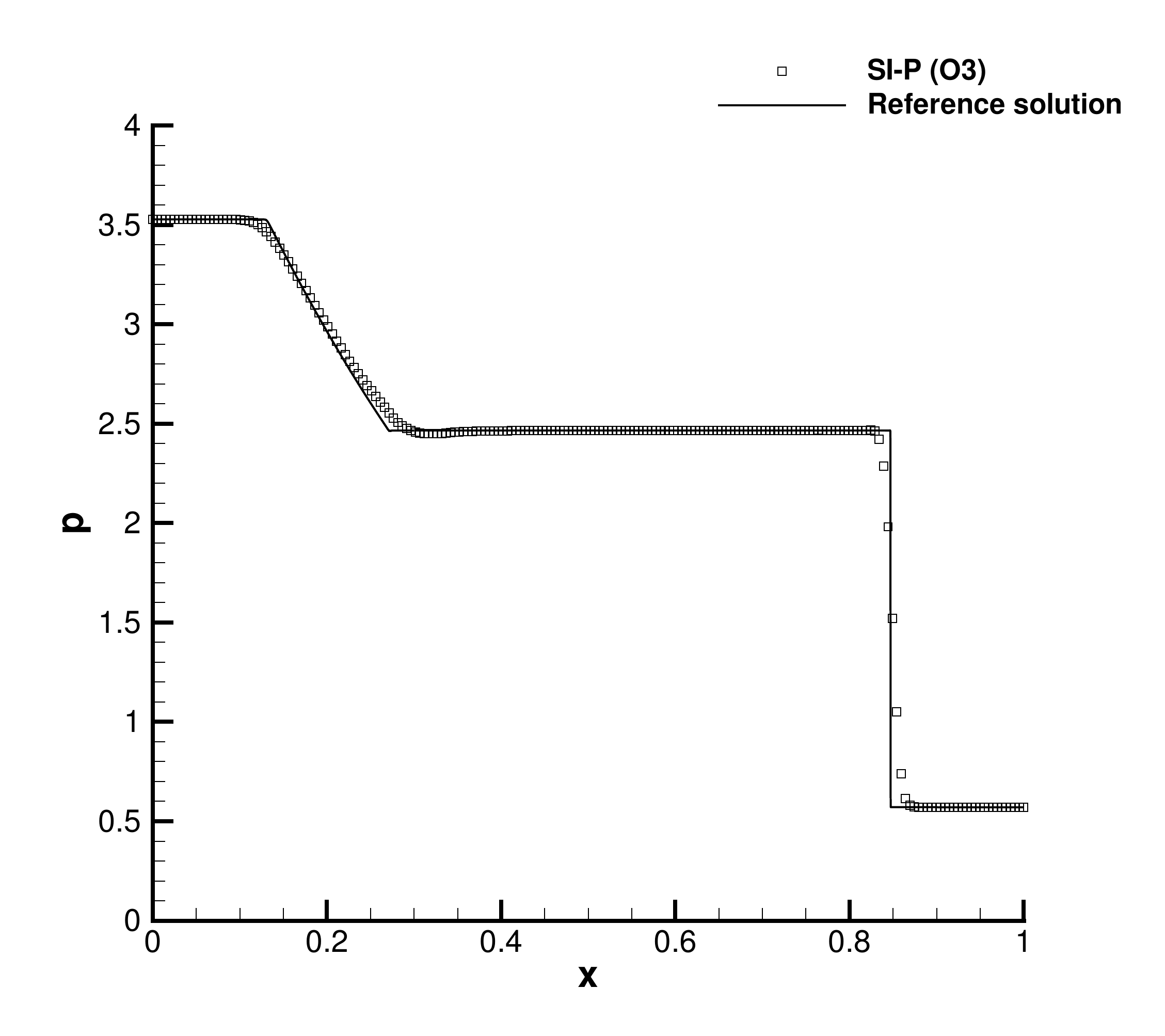} \\     
		\end{tabular} 
		\caption{Lax shock tube problem (RP1) at final time $t_f=0.14$. Comparison of density, velocity and pressure (symbols) versus the reference solution (straight line) for first and third order SI-P schemes.}
		\label{fig.Lax}
	\end{center}
\end{figure}

The second Riemann problem represents a benchmark test for Godunov-type finite volume methods, namely the Lax shock tube problem. The results are gathered in Figure \ref{fig.Lax} where a very good agreement with the reference solution can be appreciated. Moreover, the high order solution is much less dissipative compared to the first order results, hence enhancing the benefits of high order discretizations in terms of accuracy especially across rarefaction waves.

\begin{figure}[!htbp]
	\begin{center}
		\begin{tabular}{ccc} 
			\includegraphics[width=0.33\textwidth]{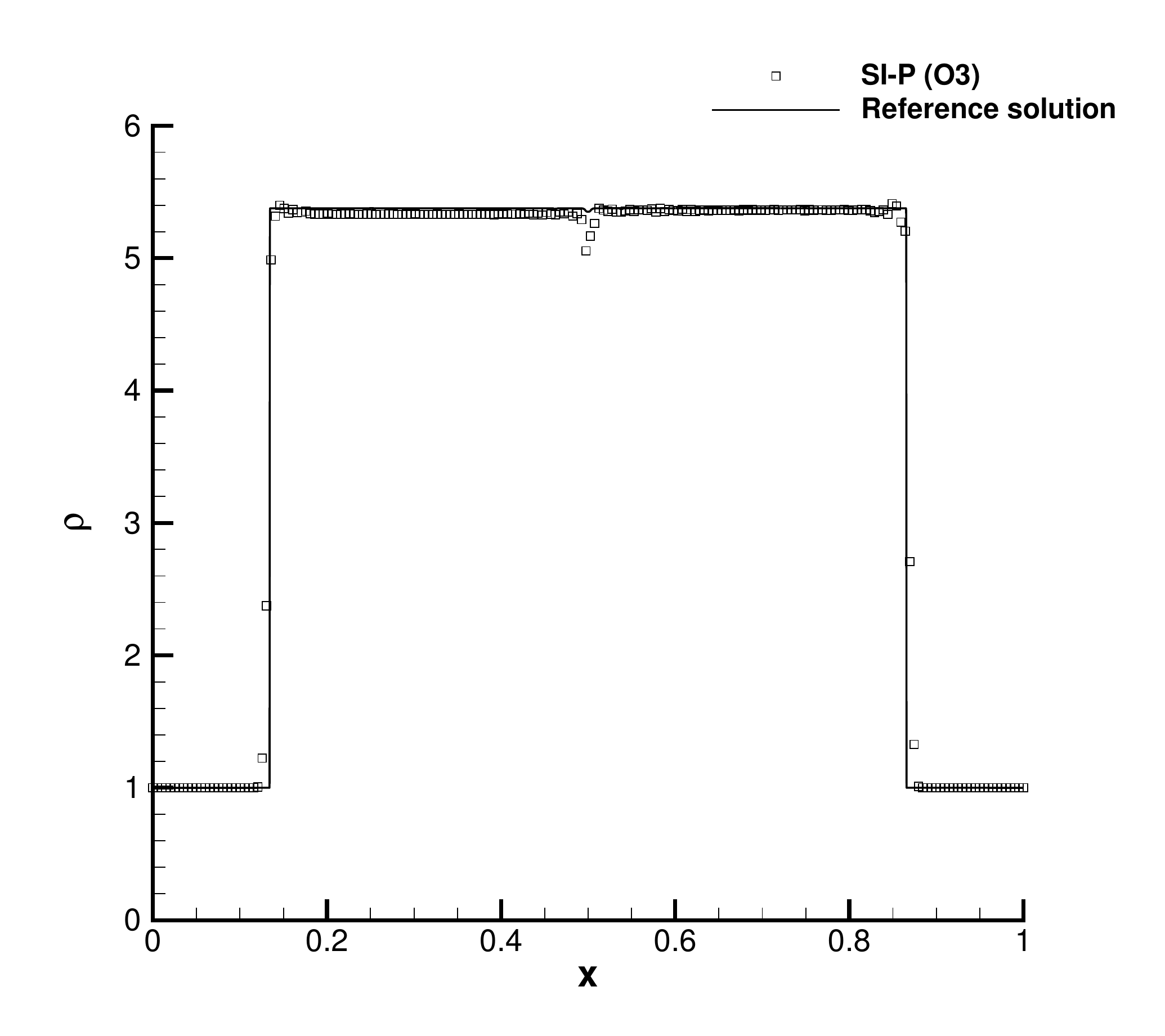}  & 
			\includegraphics[width=0.33\textwidth]{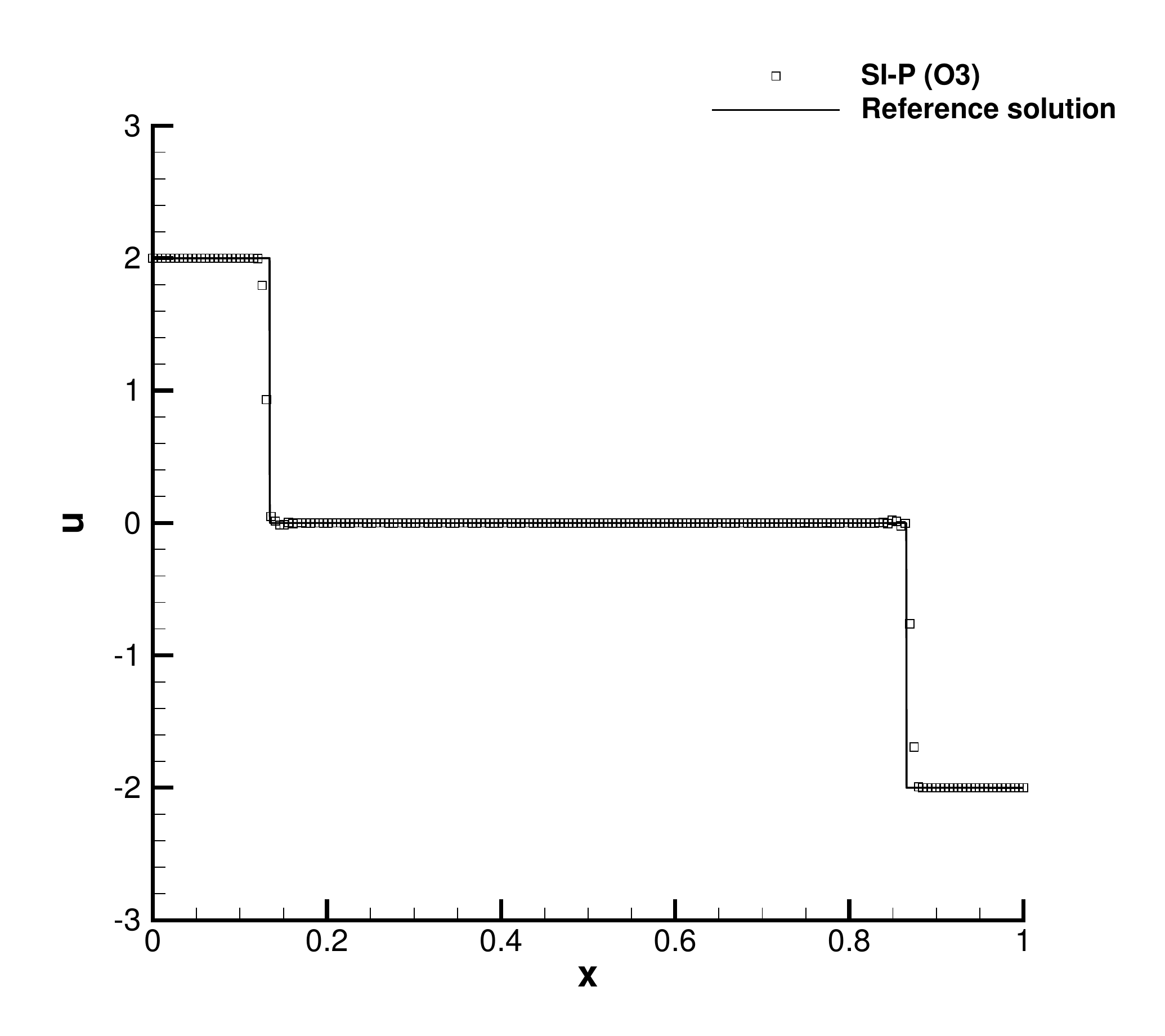}  &
			\includegraphics[width=0.33\textwidth]{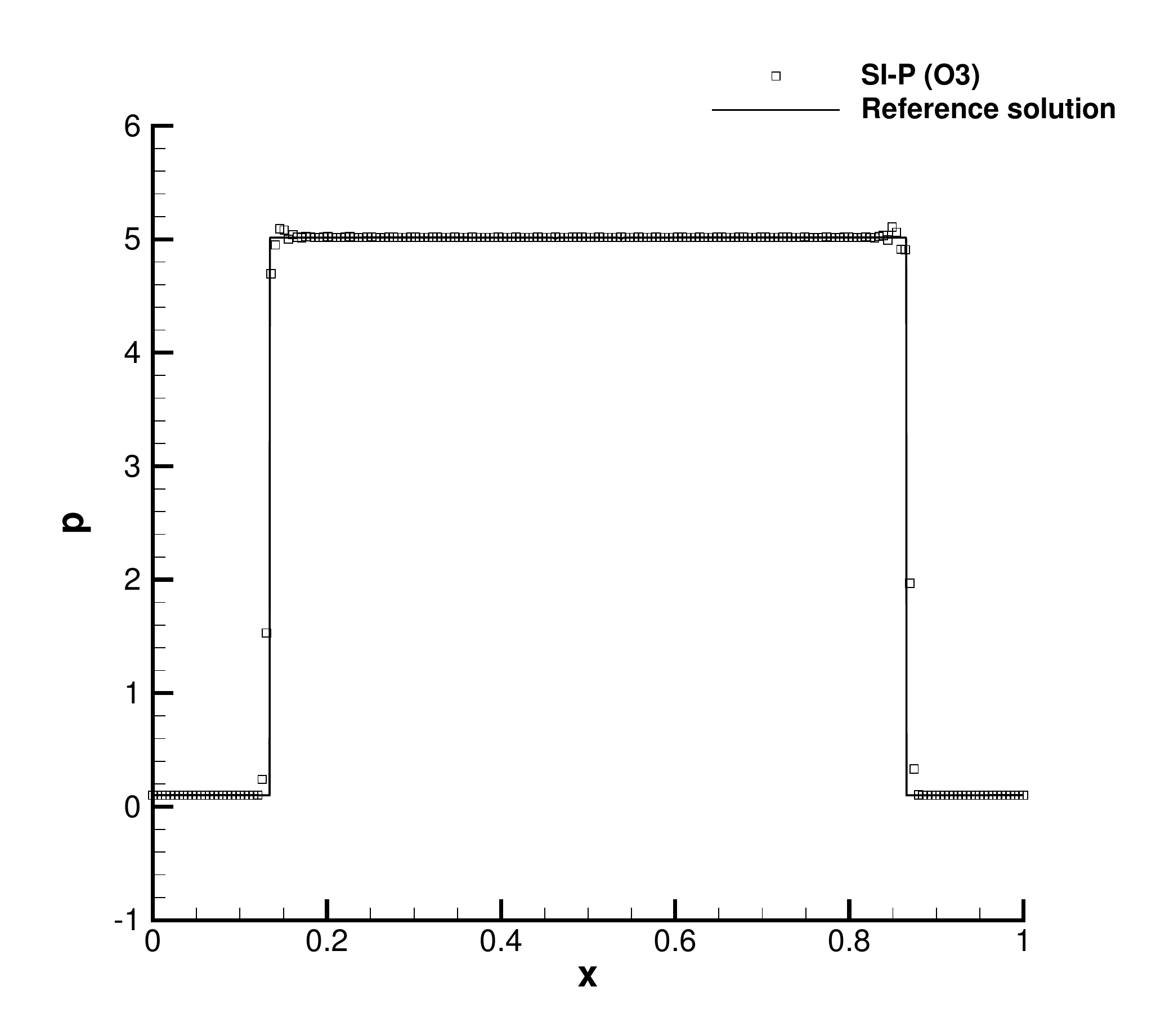} \\   
		\end{tabular} 
		\caption{Colliding shock test (RP2) at final time $t_f=0.8$. Comparison of density, velocity and pressure (symbols) versus the reference solution (straight line) for third order SI-P schemes.}
		\label{fig.2shock}
	\end{center}
\end{figure}

\begin{figure}[!htbp]
	\begin{center}
		\begin{tabular}{ccc} 
			\includegraphics[width=0.33\textwidth]{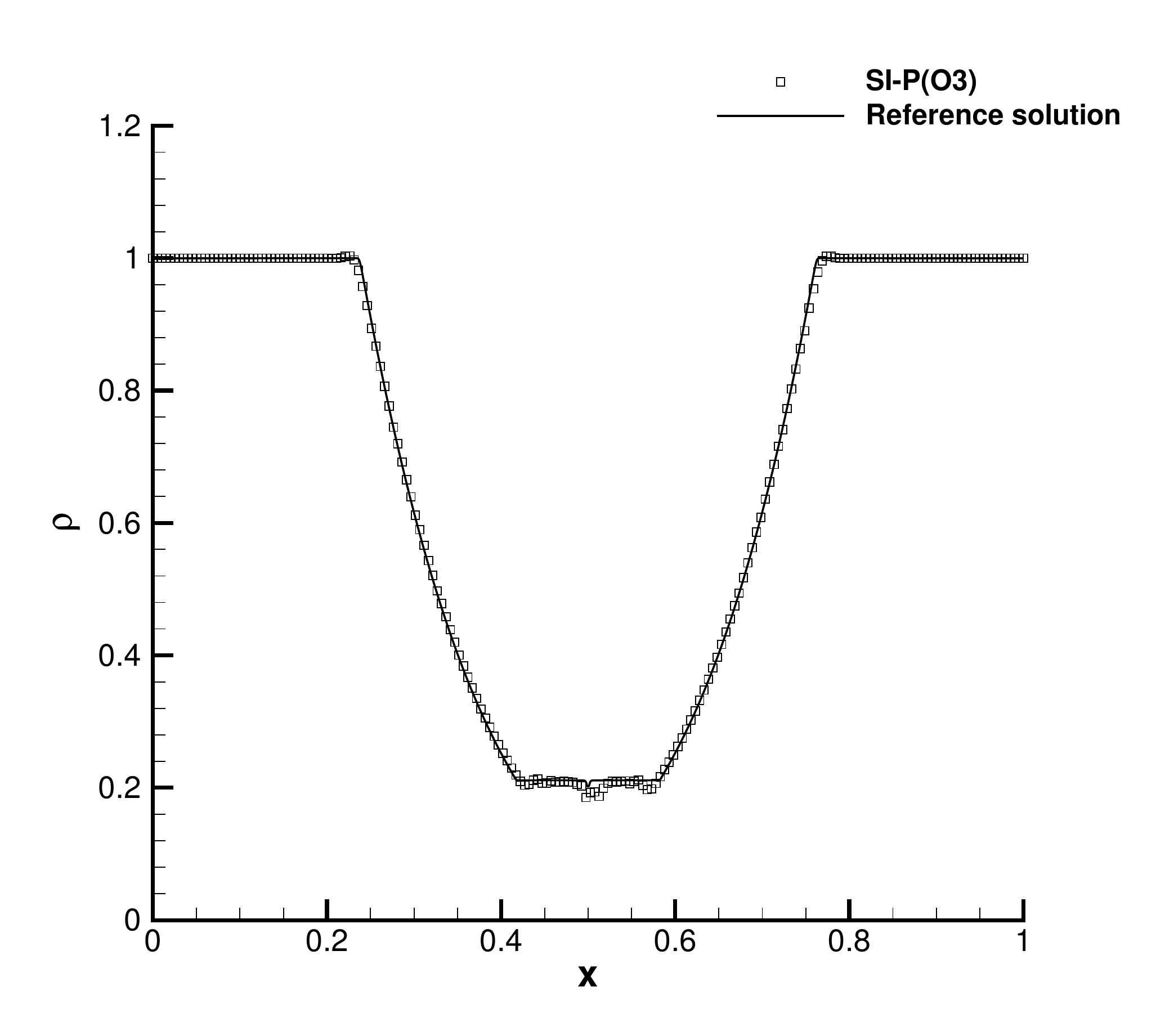}  & 
			\includegraphics[width=0.33\textwidth]{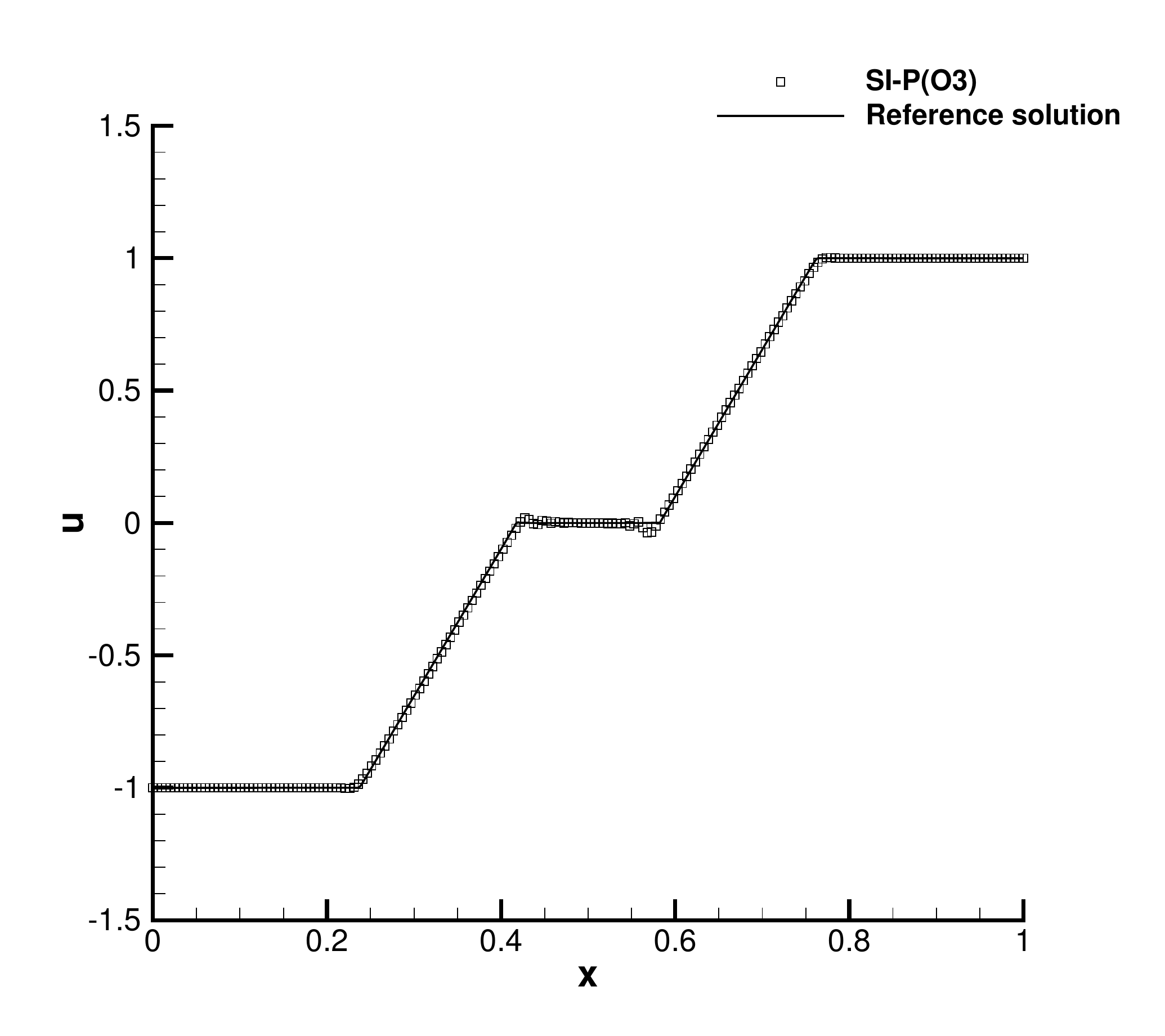}  &
			\includegraphics[width=0.33\textwidth]{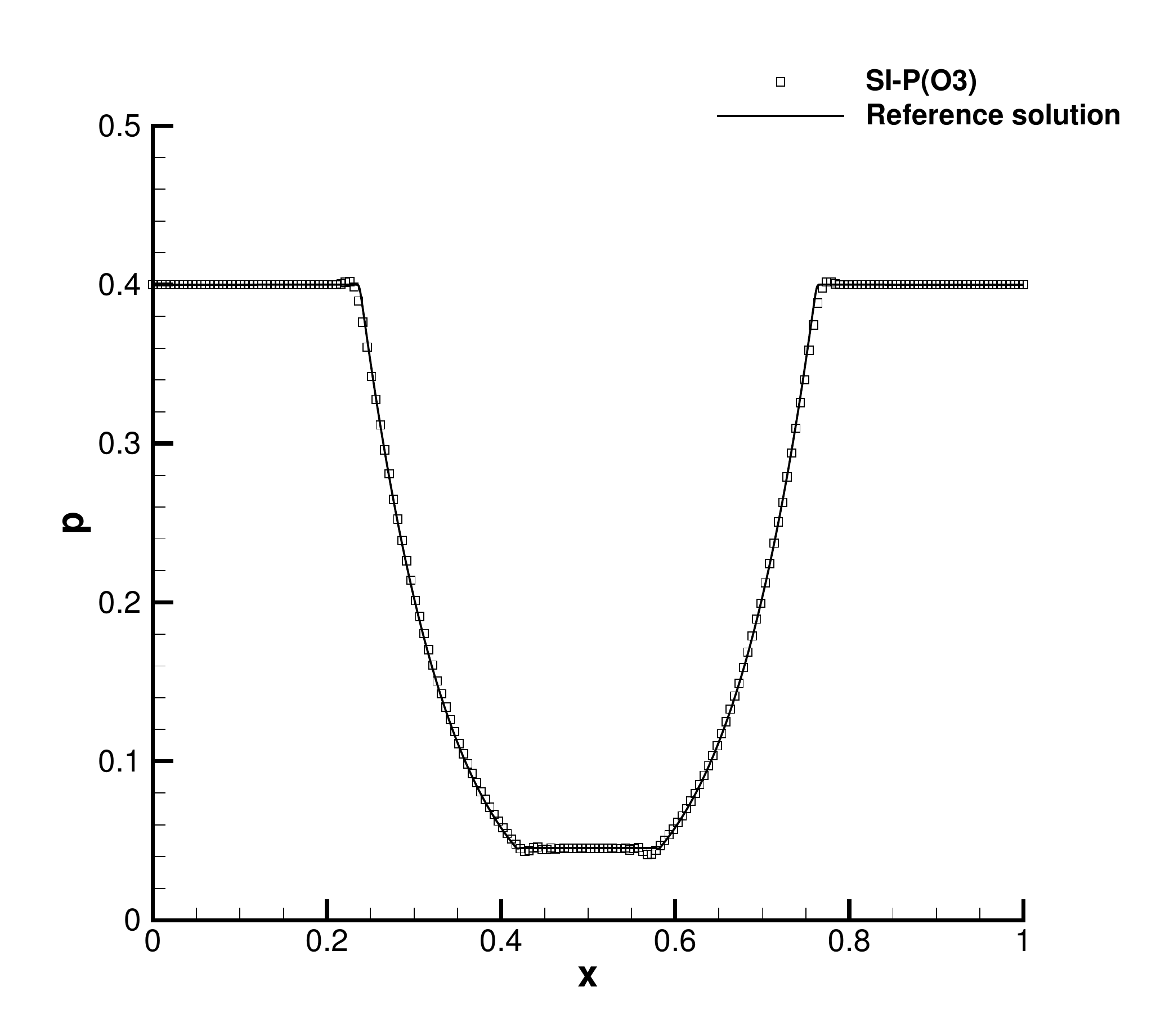} \\   
		\end{tabular} 
		\caption{Double rarefaction test (RP3) at final time $t_f=0.15$. Comparison of density, velocity and pressure (symbols) versus the reference solution (straight line) for third order SI-P schemes.}
		\label{fig.2rarefaction}
	\end{center}
\end{figure}

RP2 and RP3 are concerned with two strong colliding shocks and a symmetric double rarefaction, respectively. Figures \ref{fig.2shock} and \ref{fig.2rarefaction} plot the numerical solution at the final time of the simulation for RP2 and RP3, respectively. For both tests some nonphysical oscillations can be noticed in the density profile which are also present in \cite{Dumbser_Casulli16}, while the numerical solution for velocity and pressure is overall in good agreement with the reference solution. Very small perturbations occur at the shock waves, but this is well known also for explicit Godunov-type finite volume schemes as pointed out in \cite{ToroBook}.

\begin{figure}[!htbp]
	\begin{center}
		\begin{tabular}{cc} 
			\includegraphics[width=0.47\textwidth]{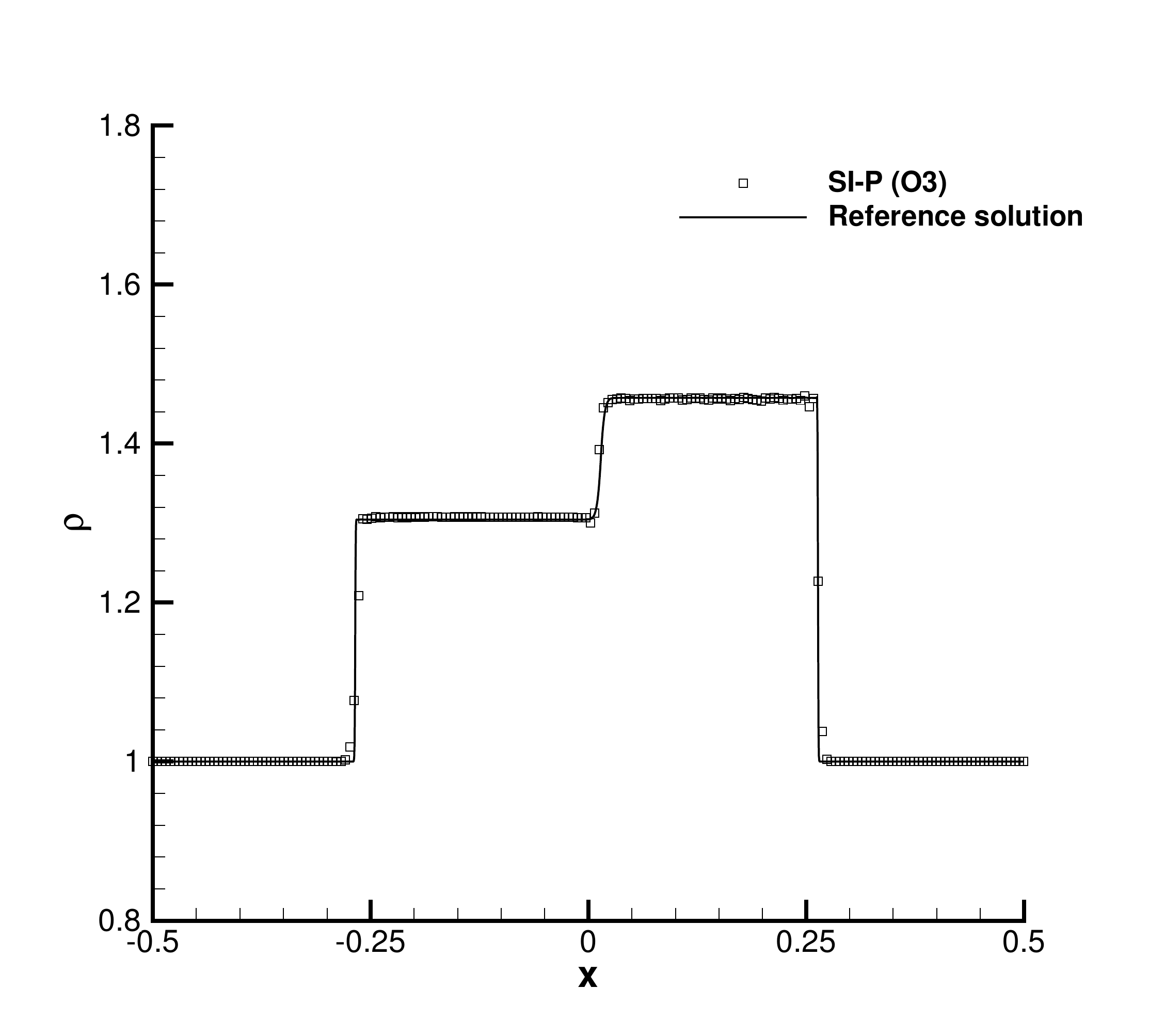} &           
			\includegraphics[width=0.47\textwidth]{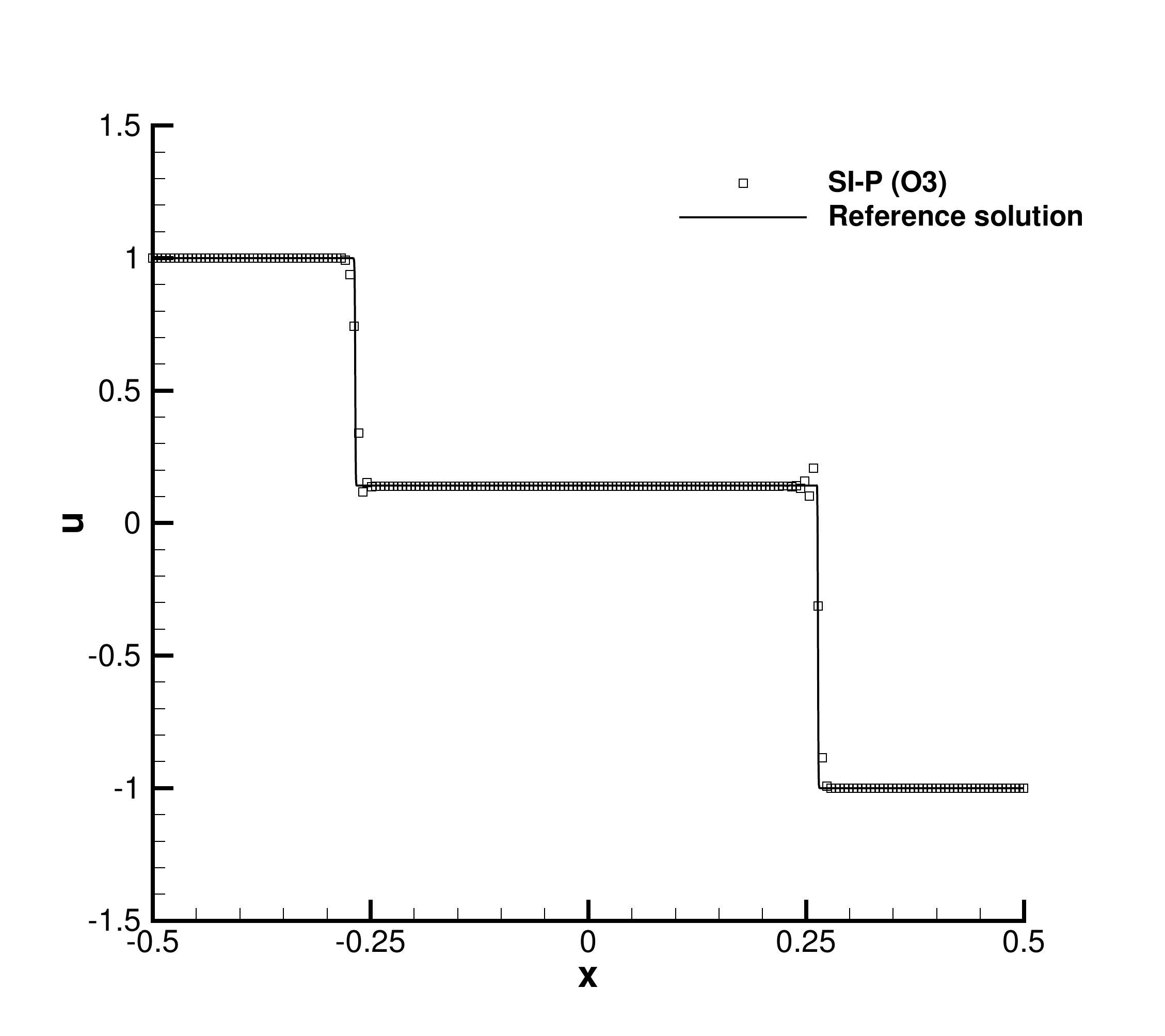} \\
			\includegraphics[width=0.47\textwidth]{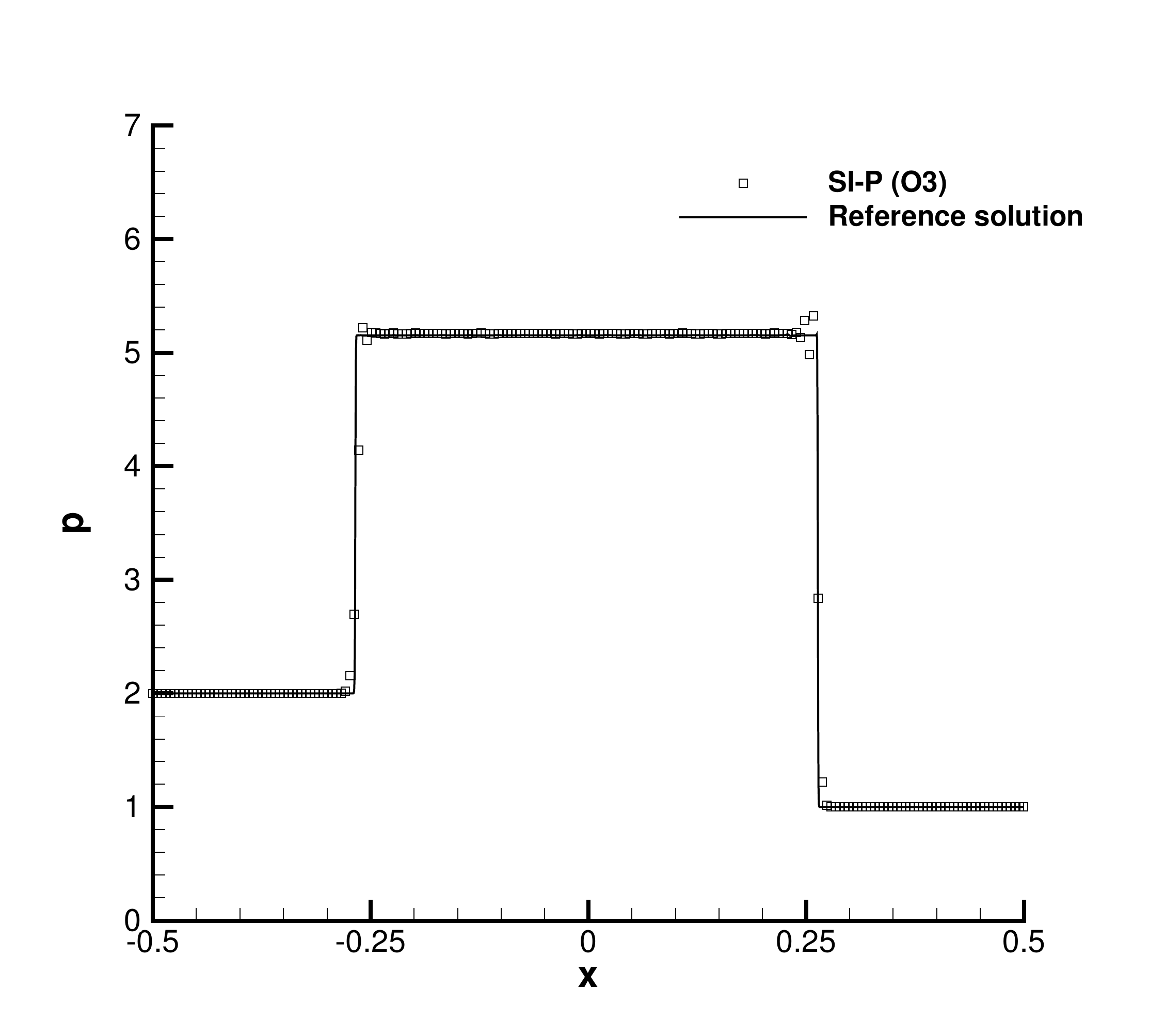} &           
			\includegraphics[width=0.47\textwidth]{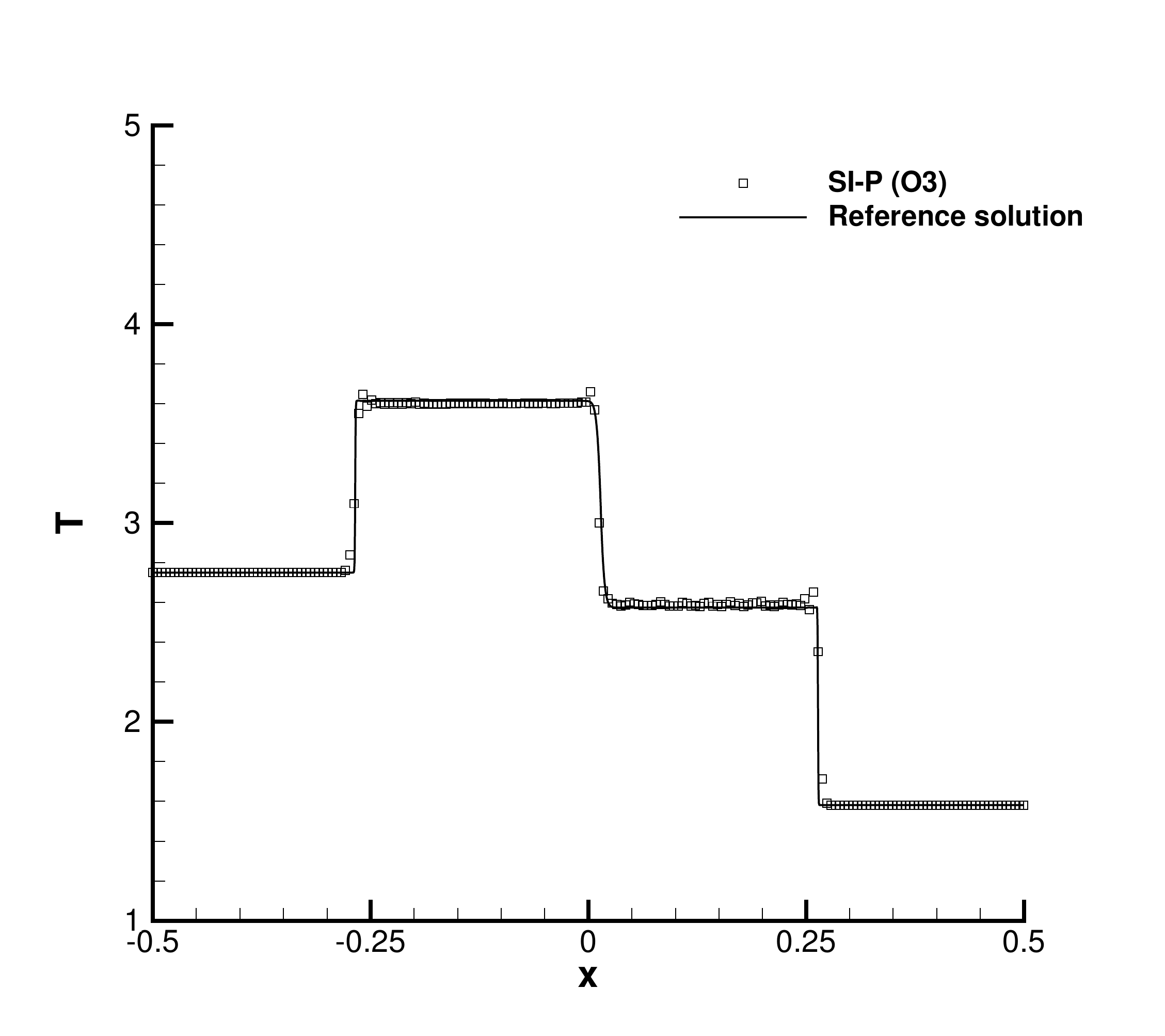} \\
		\end{tabular} 
		\caption{Shock tube problem RK1 at final time $t_f=0.1$ for the Redlich-Kwong EOS. Comparison of density, velocity, pressure and temperature (symbols) versus the reference solution (straight line) for third order SI-P schemes.}
		\label{fig.RK1}
	\end{center}
\end{figure}

\begin{figure}[!htbp]
	\begin{center}
		\begin{tabular}{cc} 
			\includegraphics[width=0.47\textwidth]{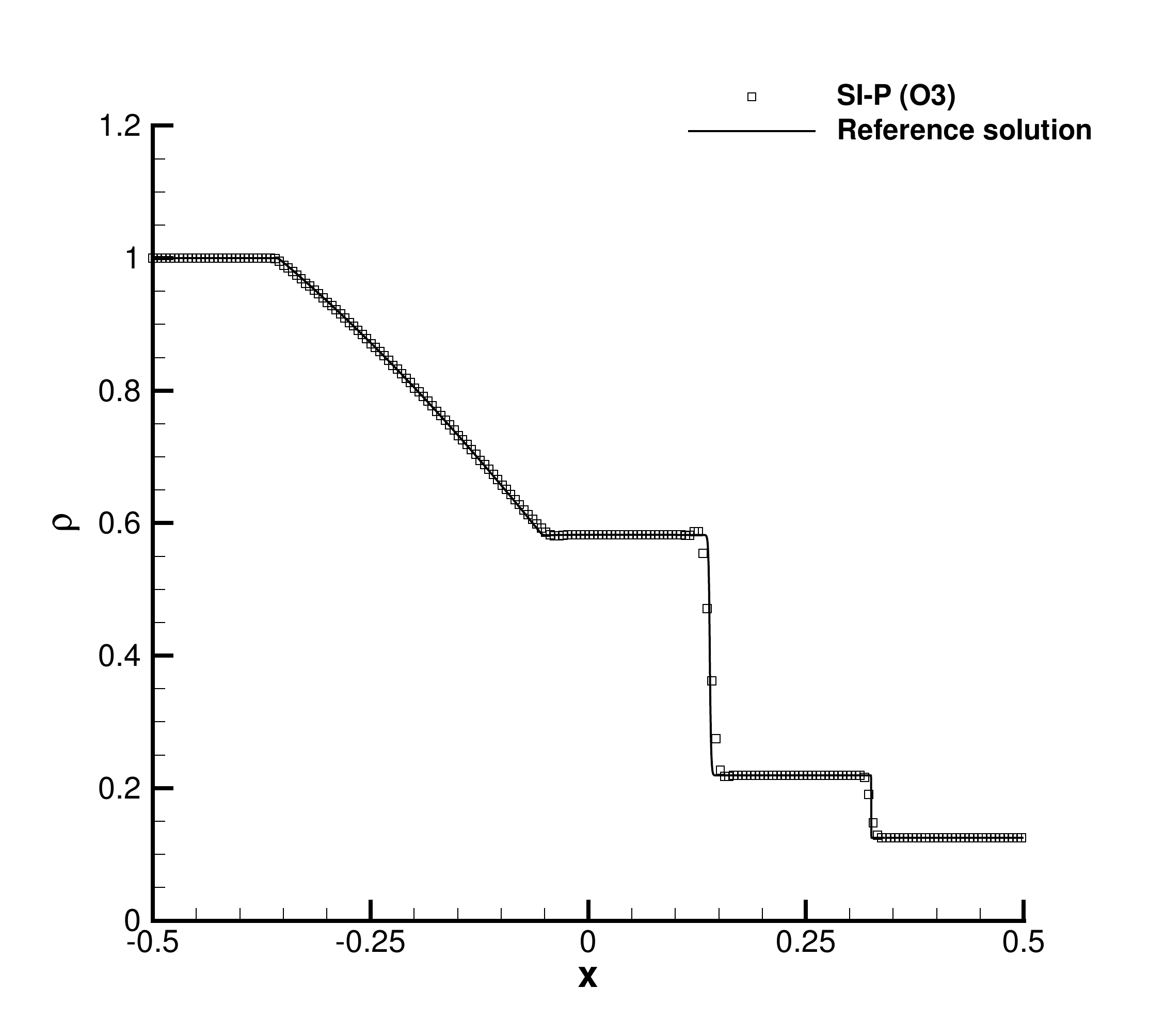} &           
			\includegraphics[width=0.47\textwidth]{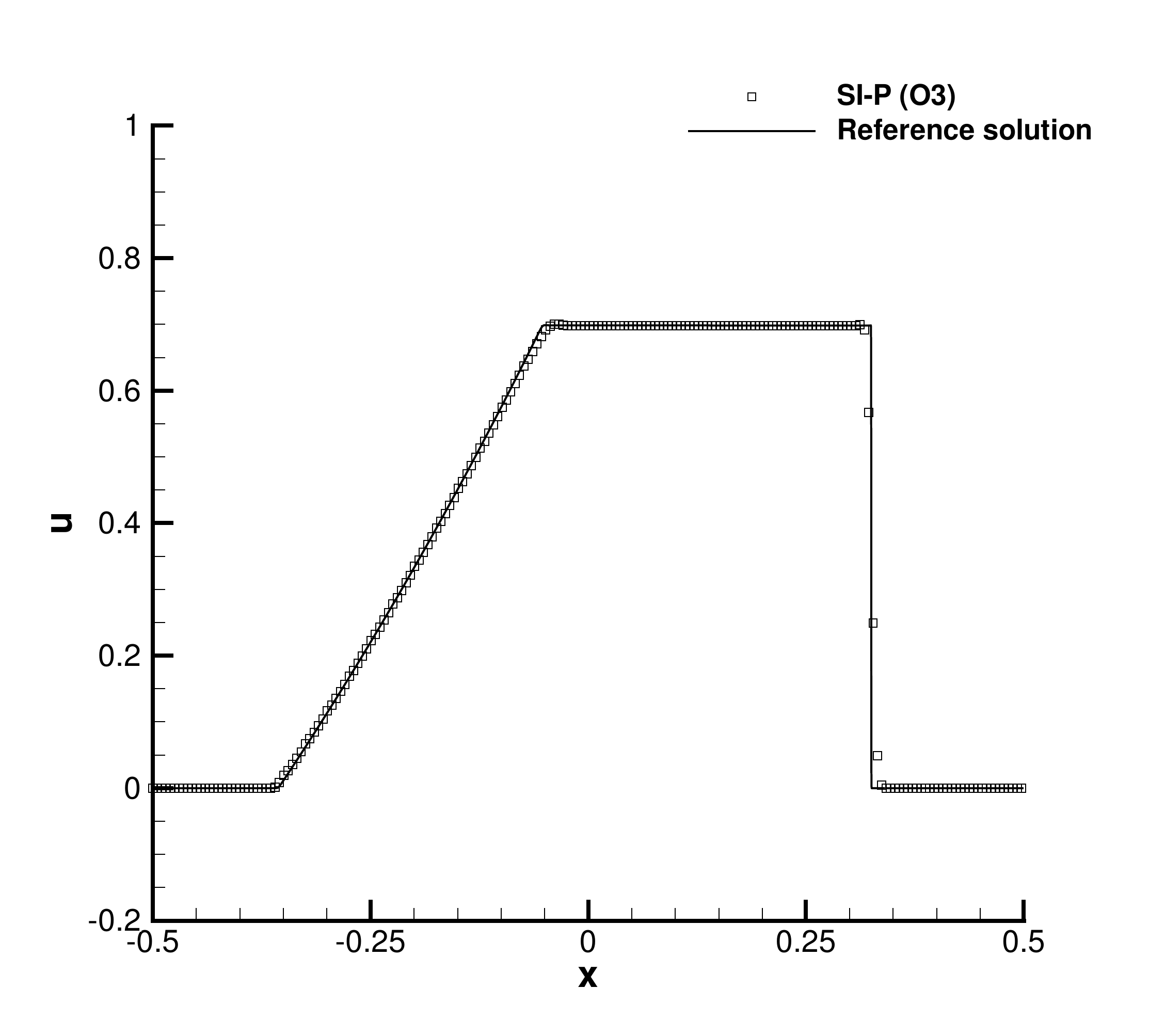} \\
			\includegraphics[width=0.47\textwidth]{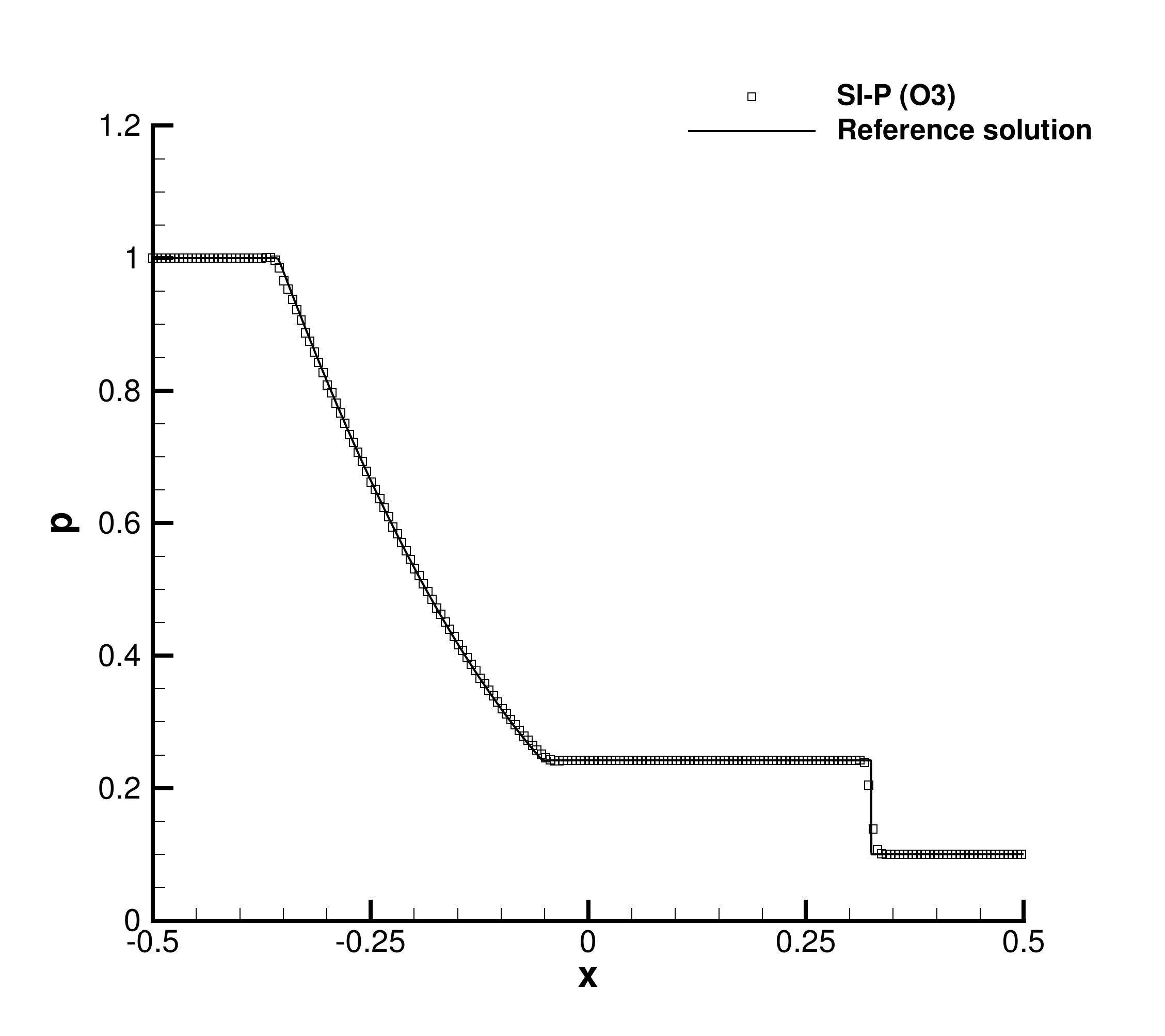} &           
			\includegraphics[width=0.47\textwidth]{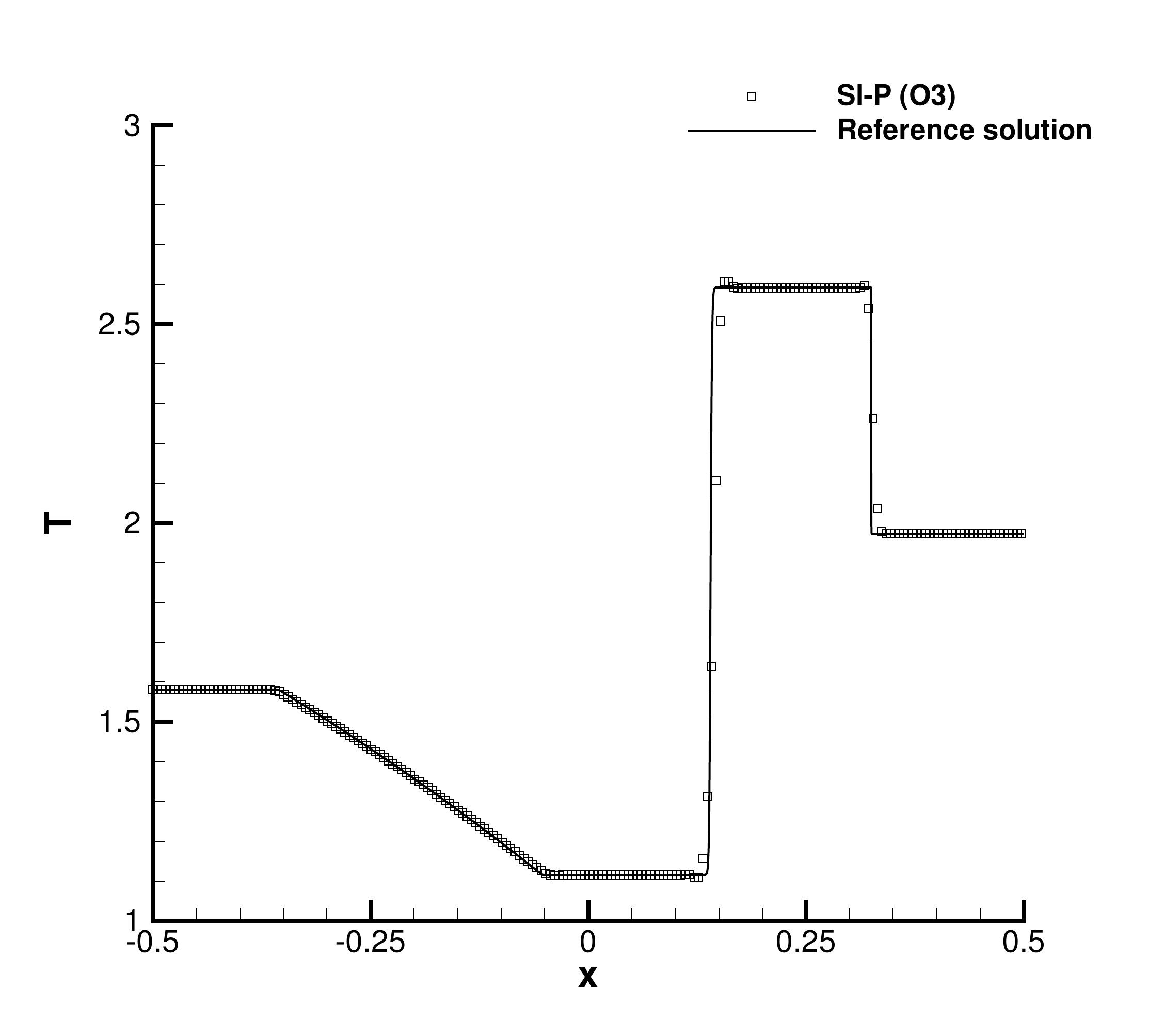} \\
		\end{tabular} 
		\caption{Shock tube problem RK2 at final time $t_f=0.2$ for the Redlich-Kwong EOS. Comparison of density, velocity, pressure and temperature (symbols) versus the reference solution (straight line) for third order SI-P schemes.}
		\label{fig.RK2}
	\end{center}
\end{figure}

A nonlinear equation of state is considered in Riemann problems RK1 and RK2, which results are shown in Figures \ref{fig.RK1} and \ref{fig.RK2}, respectively. In particular, RK2 is the Sod shock tube problem that has been run using the Redlich-Kwong EOS, thus obtaining very different waves in terms of profile and location compared to the same test run with classical ideal gas EOS (see \cite{ToroBook}). The temperature distribution is also shown and the results reasonably match the reference solution. The Newton algorithm \eqref{eqn.pNewton} for the solution of the pressure wave equation \eqref{eqn.p_SI-discr}-\eqref{eqn.b-SI-discr} has always converged to a tolerance $\delta=10^{-10}$ in at most four iterations.

\begin{figure}[!htbp]
	\begin{center}
		\begin{tabular}{cc} 
			\includegraphics[width=0.47\textwidth]{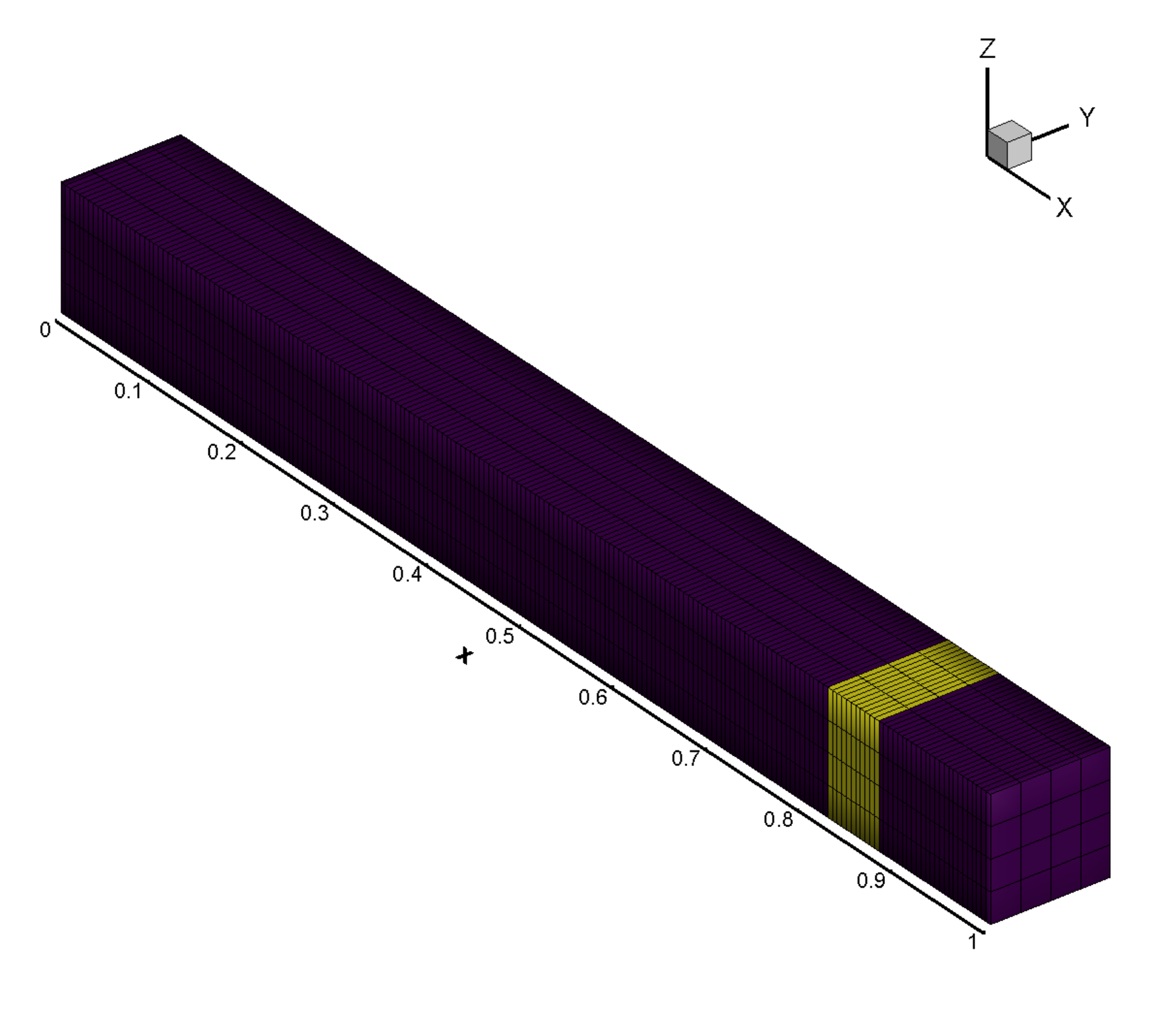}  & 
			\includegraphics[width=0.47\textwidth]{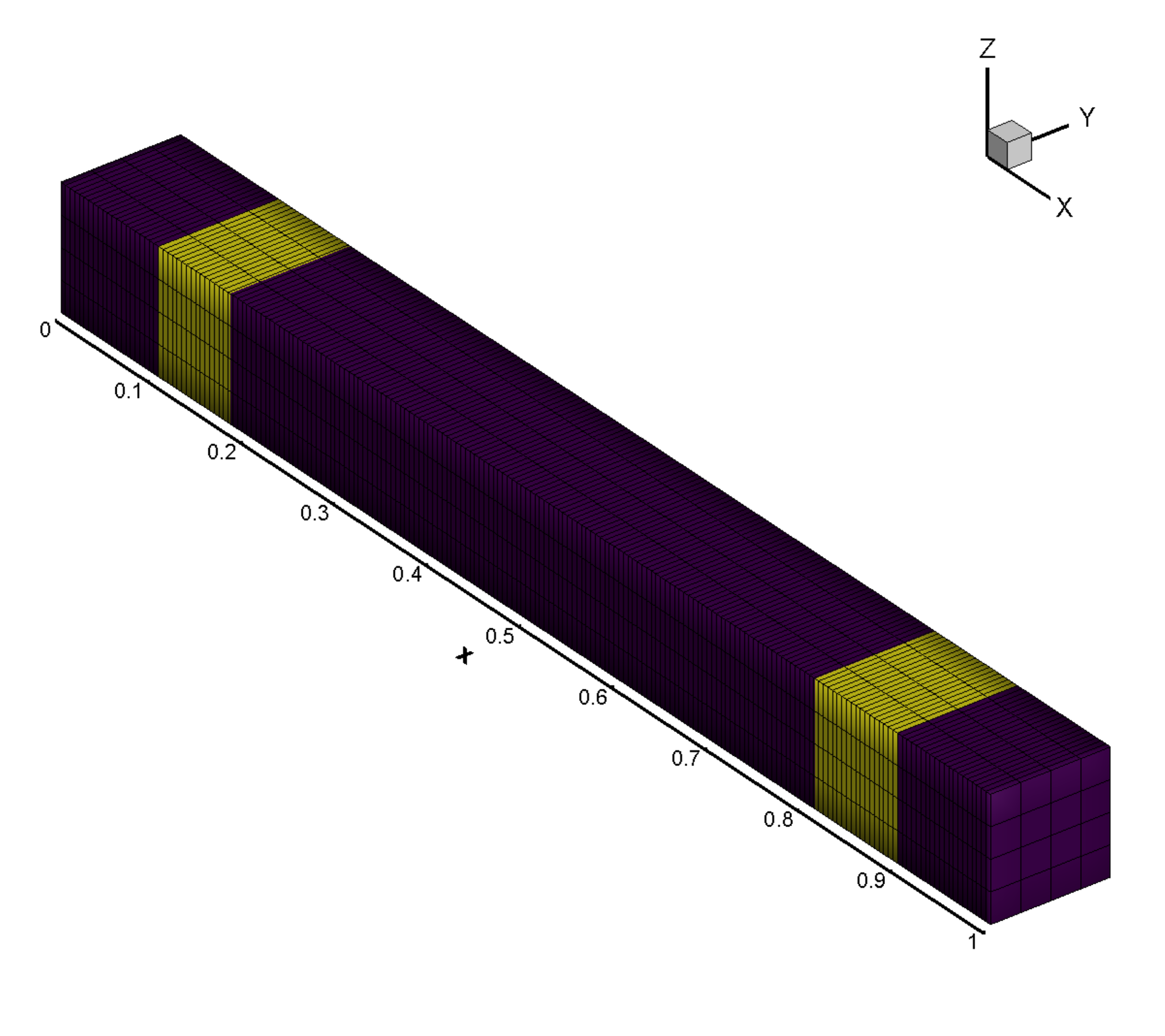} \\  
			\includegraphics[width=0.47\textwidth]{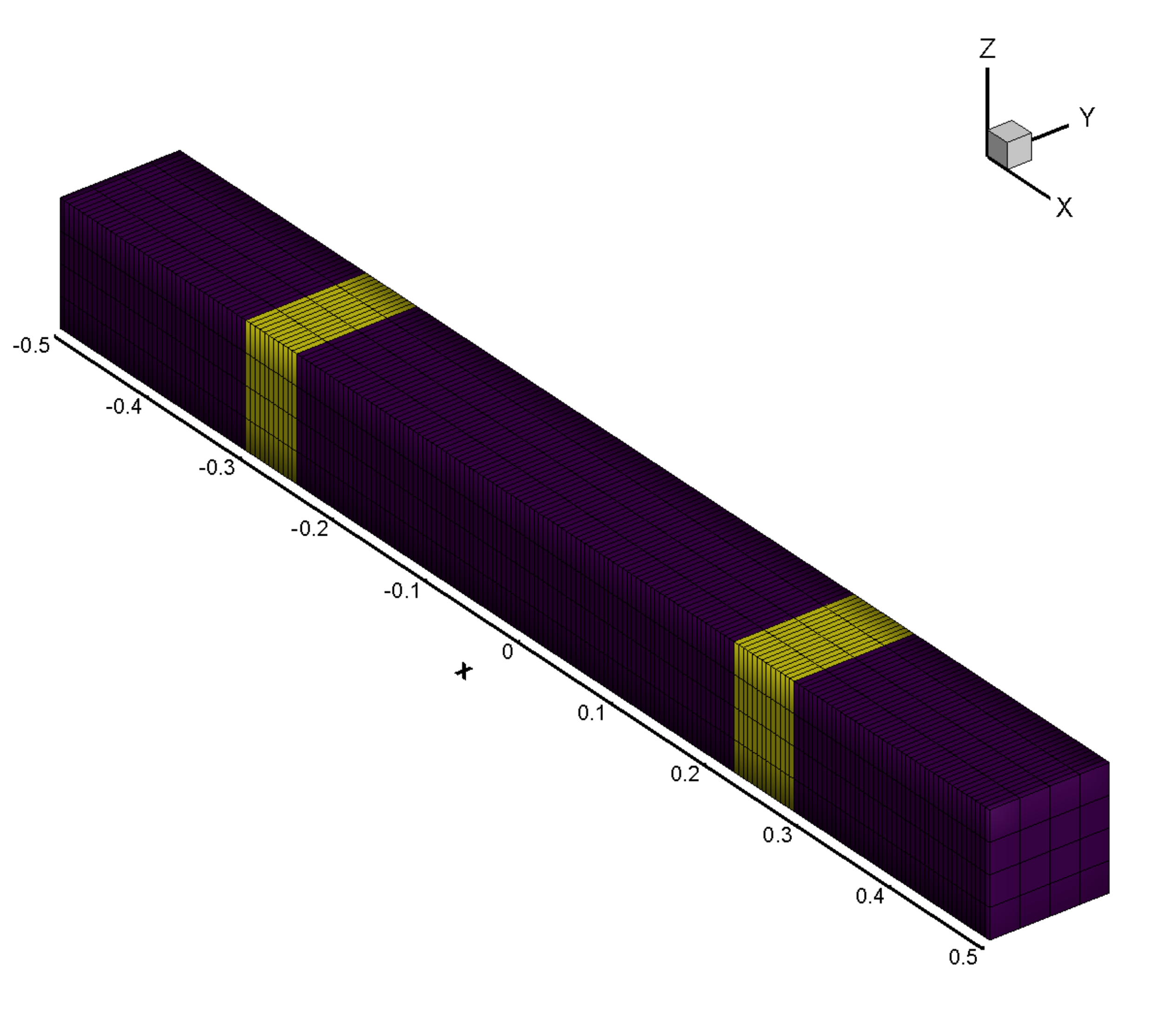}  & 
			\includegraphics[width=0.47\textwidth]{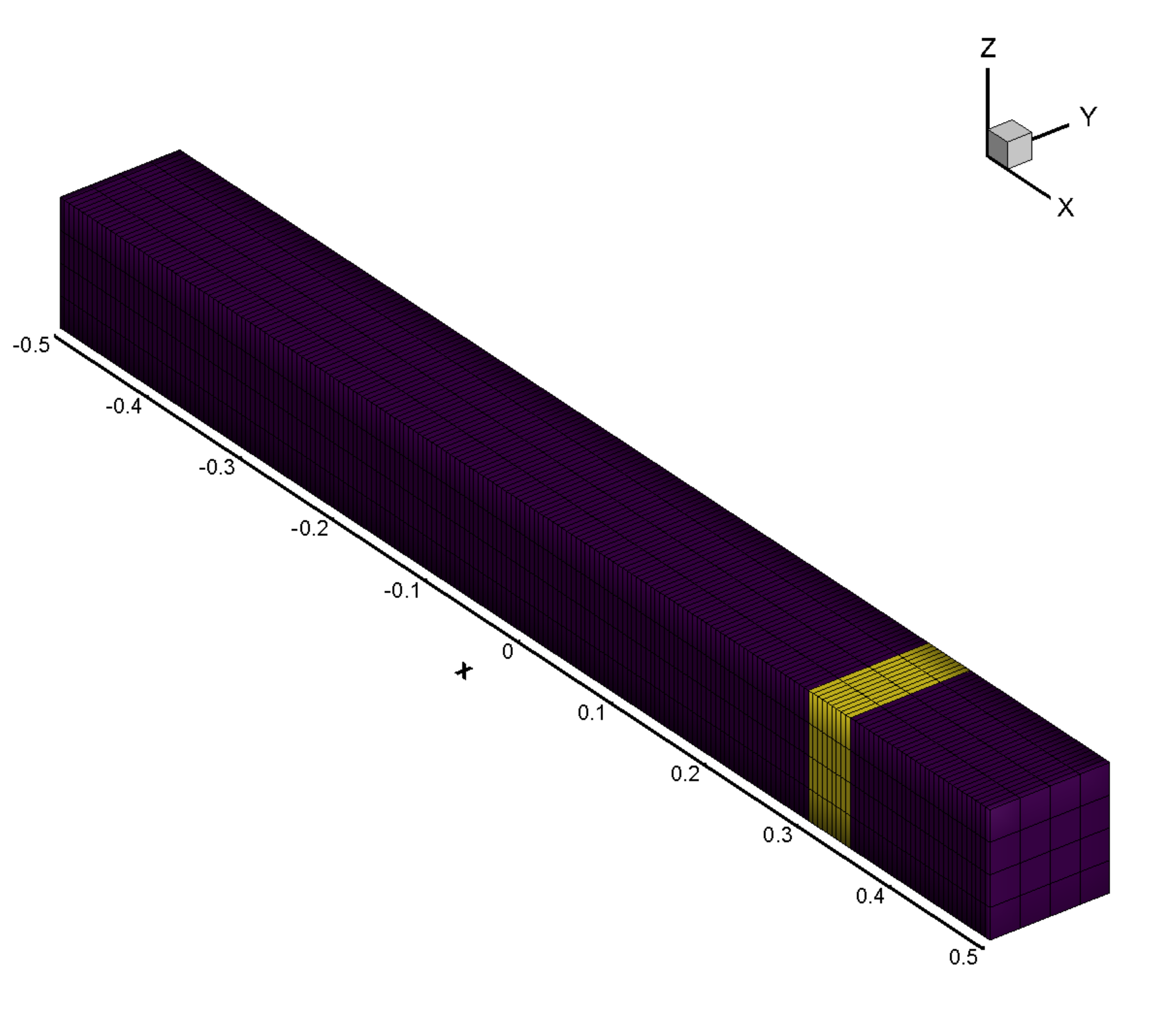} \\      
		\end{tabular} 
		\caption{Flattener indicator for RP1 (top left), RP2 (top right), RK1 (bottom left) and RK2 (bottom right) test at the final time of the simulation.}
		\label{fig.RP-flattenerRP}
	\end{center}
\end{figure}

Finally, Figure \ref{fig.RP-flattenerRP} show the regions of the computational domain where the flattener $\chi^n$ has been activated and the solution has been updated with the convex combination \eqref{eqn.Q1flattener}. The troubled cells are indeed identified only where shocks are located and do not involve a large number of control volumes, thus the majority of the computational cells evolve the solution using the fully high order space and time scheme.

\subsection{Gresho vortex}
\begin{figure}[!htbp]
	\begin{center}
		\begin{tabular}{cc} 
			\includegraphics[width=0.47\textwidth]{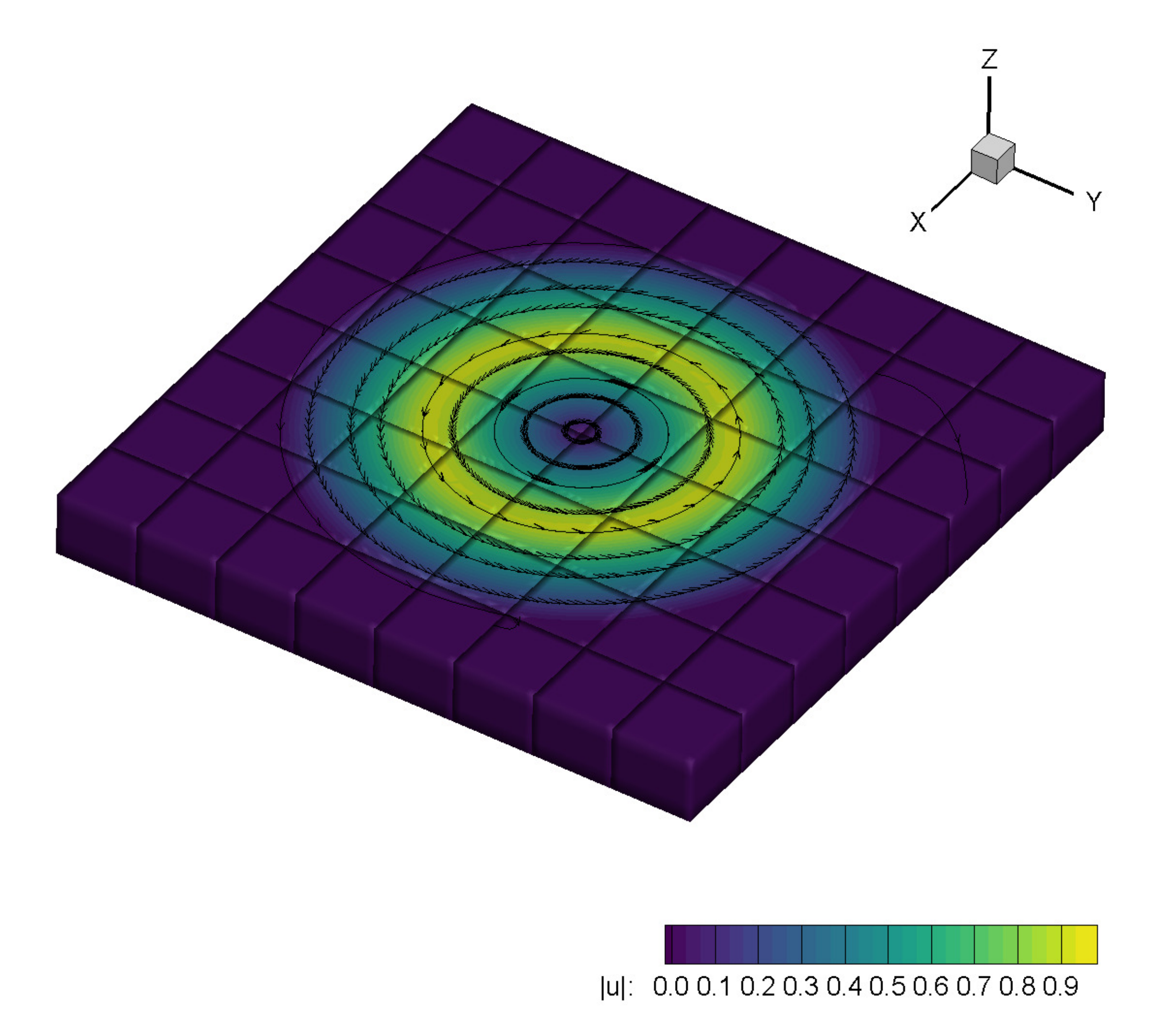} & 
			\includegraphics[width=0.47\textwidth]{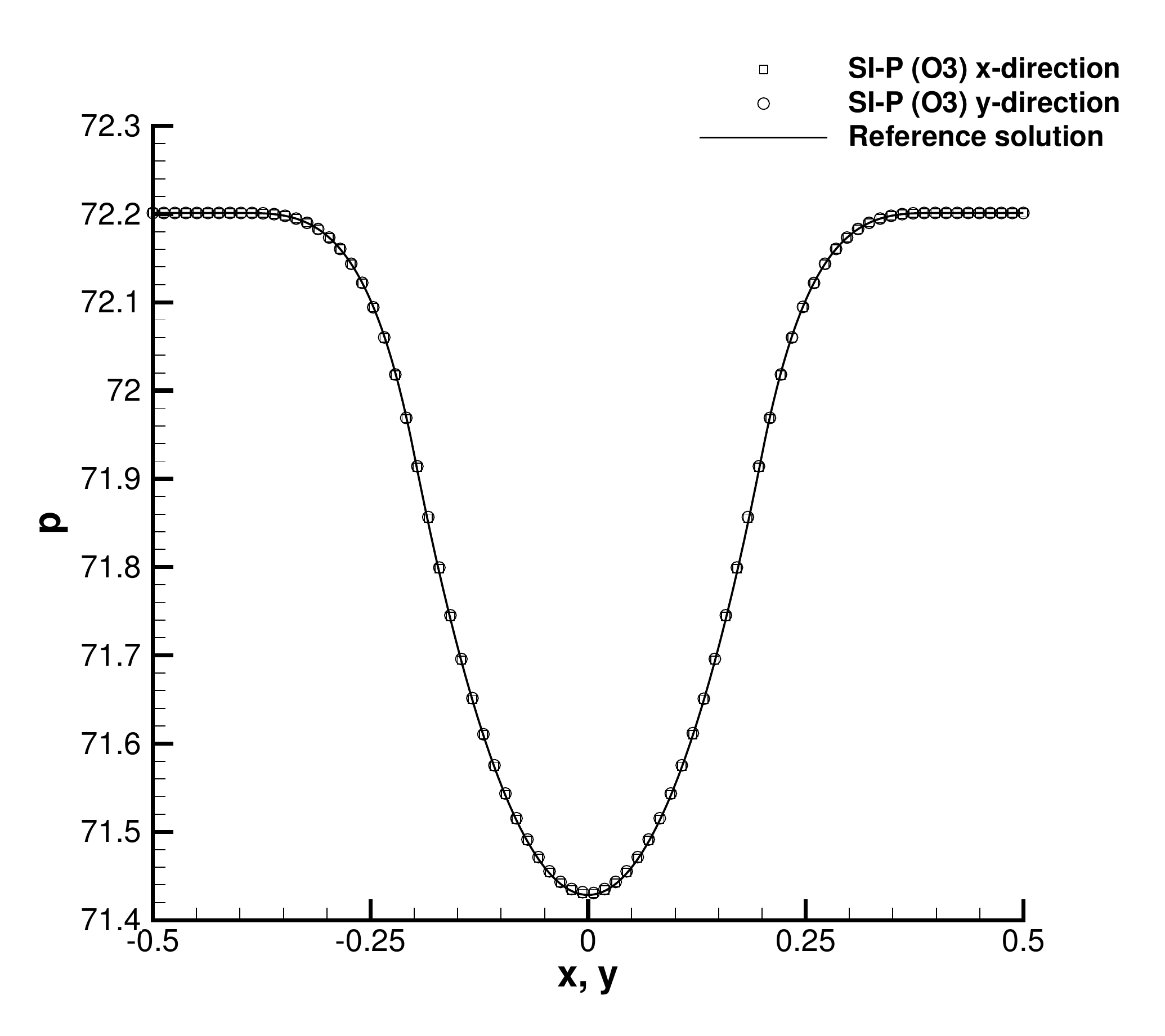} \\          
			\includegraphics[width=0.47\textwidth]{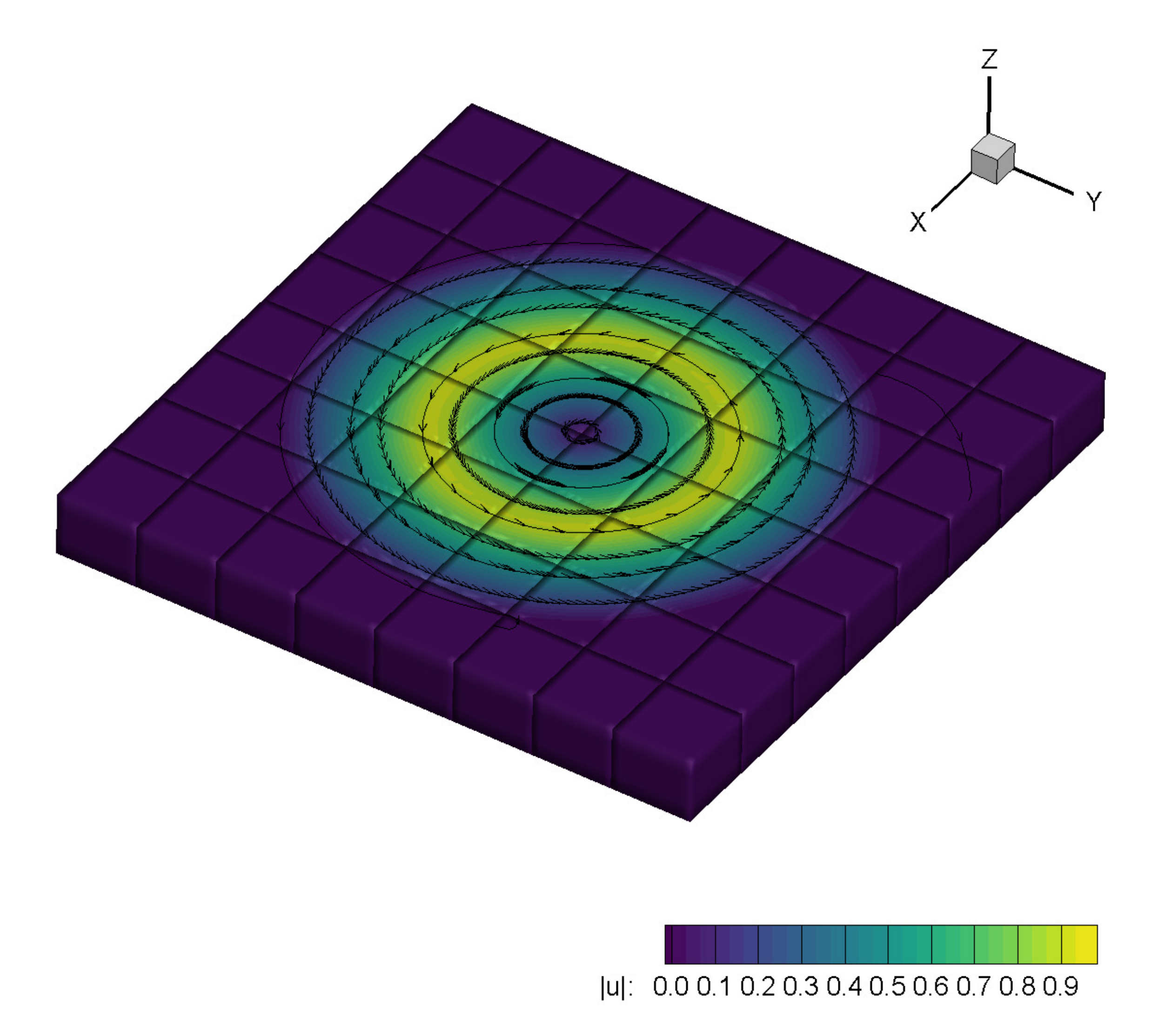} &
			\includegraphics[width=0.47\textwidth]{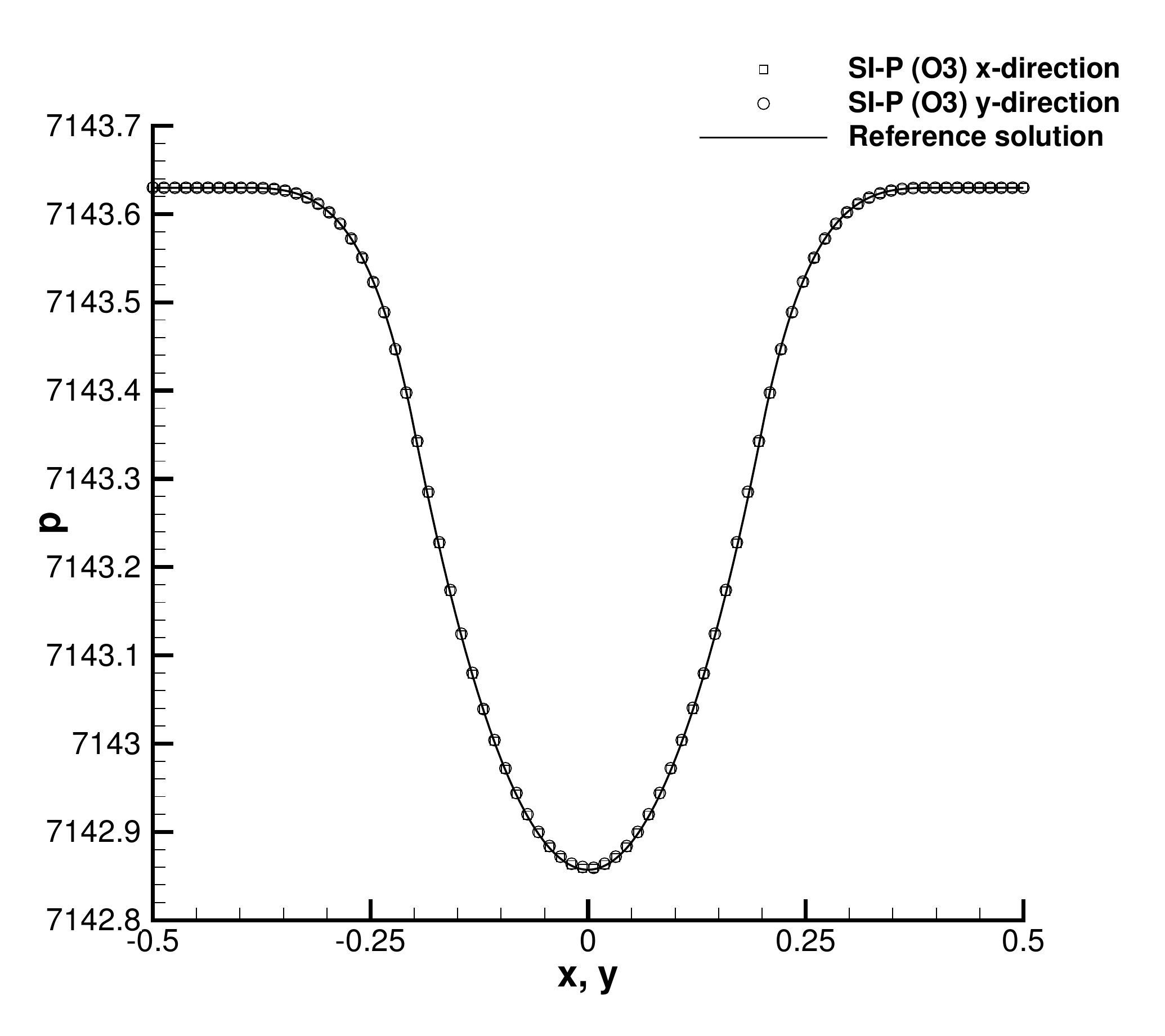} \\
			\includegraphics[width=0.47\textwidth]{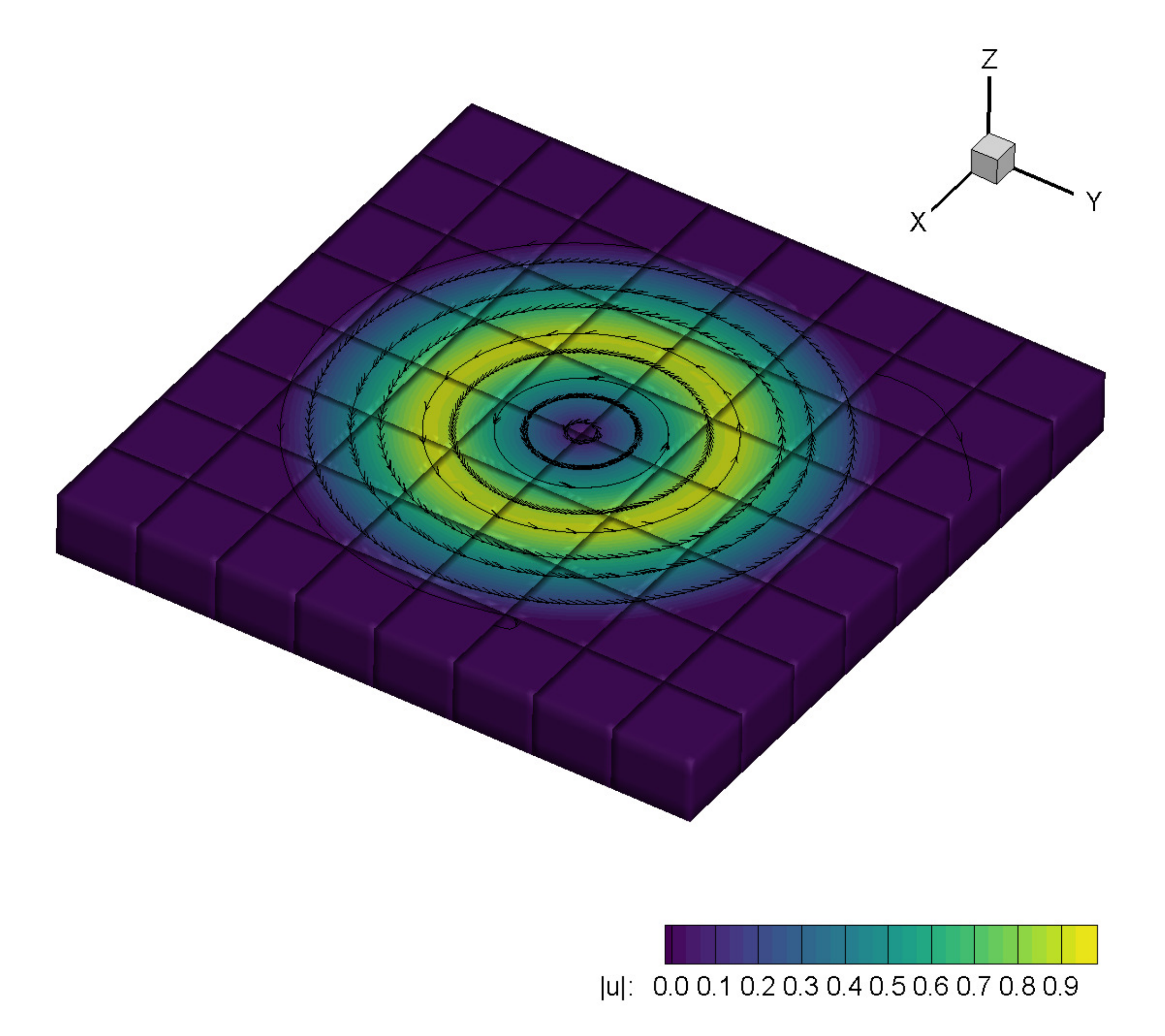} &
			\includegraphics[width=0.47\textwidth]{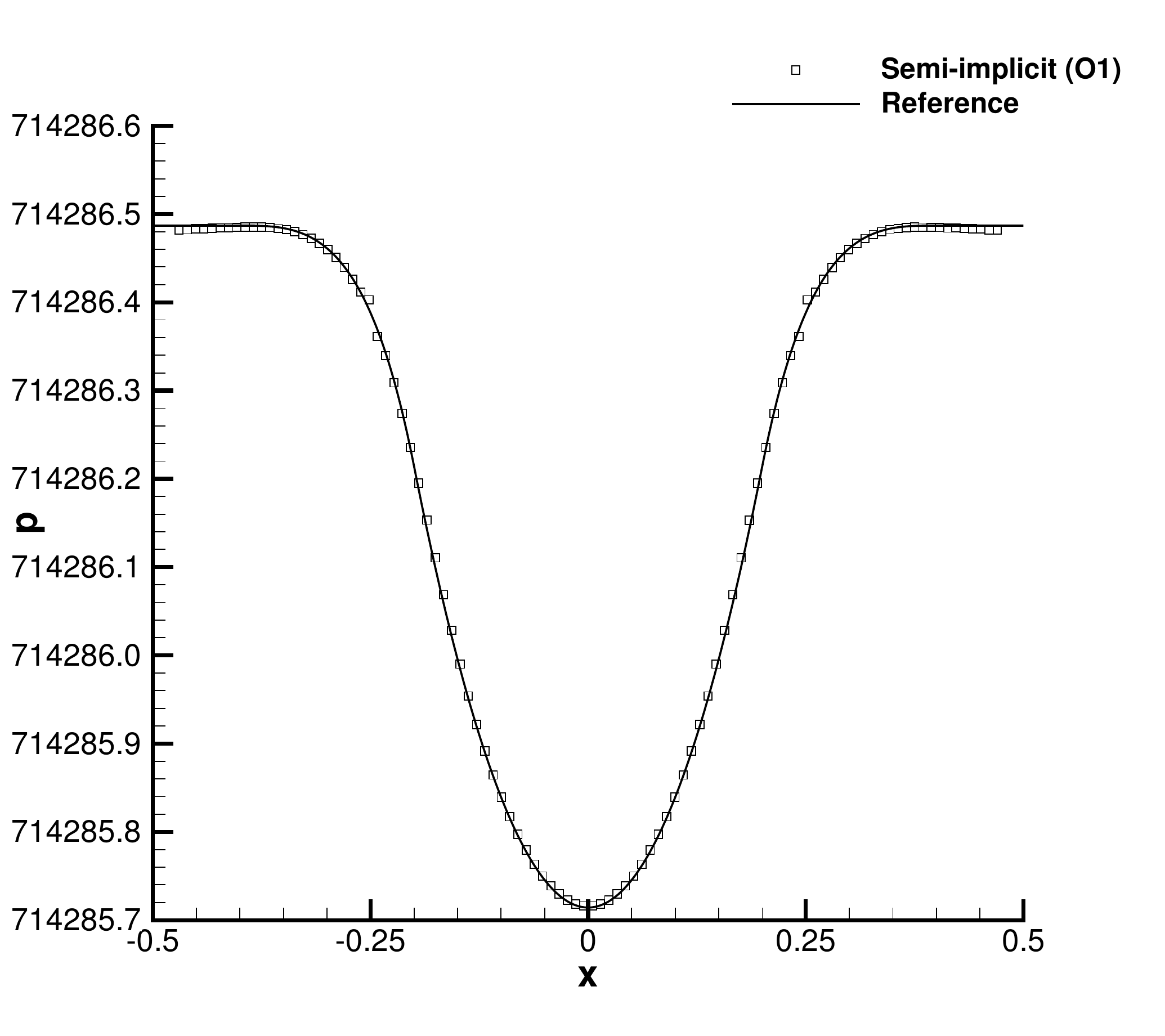} \\          
		\end{tabular} 
		\caption{Gresho vortex problem with third order SI-P scheme at the final time $t_f=0.4 \, \pi$ with Mach number $M=10^{-1}$ (top), $M=10^{-2}$ (middle), $M=10^{-3}$ (bottom). Left: stream-traces of the velocity field with velocity magnitude contours ($30$ levels have been used in the range $[3\cdot 10^{-6};1]$ for all Mach numbers). Right: pressure distribution versus reference solution along a 1D cut in $x-$direction ($y=z=0$) with 80 interpolation points.}
		\label{fig.Gresho}
	\end{center}
\end{figure}

The so-called Gresho vortex problem \cite{Gresho1990} is a known stationary solution of the Euler equation that is typically used for assessing the behavior of numerical methods at different Mach number. The computational domain is defined by $\Omega=[-0.5;0.5]^3$ with Dirichlet boundaries, where the initial condition is imposed. This is given in polar coordinates for density, angular velocity and pressure with $r=\sqrt{x^2+y^2}$ denoting the generic radial position on the $x-y$ plane and $\theta=\arctan(y/x)$ is the corresponding angle:
\begin{eqnarray}
\rho(r)&=&1 \nonumber \\
u_\theta(r) &=& \left\{\begin{array}{l}
\begin{tabular}{ll}
$5r$& $0\leq r <0.2$\\
$2-5r$ & $0.2\leq r <0.4$\\
$0$ & $r \geq 0.4$
\end{tabular}
\end{array}
\right.
\nonumber \\
p(r) &=& \left\{\begin{array}{l}
\begin{tabular}{ll}
$p_0+\frac{25}{2}r^2$	& $0\leq r <0.2$\\
$p_0+\frac{25}{2}r^2+4[1-5r-\ln(0.2)+\ln(r)]$ 		& $0.2\leq r <0.4$\\
$p_0-2+4\ln(2)$ 					& $r \geq 0.4$
\end{tabular}
\end{array}
\right., 
\label{gresho.exact}
\end{eqnarray}
where the background pressure $p_0=\rho / (\gamma M^2)$ is expressed in terms of the Mach number. The velocity field with Cartesian components can be easily obtained from $u_\theta$ with a rotation, that is $(u,v)=u_\theta/r \cdot (-y,x)$. This test is run until the final time $t_f=0.4 \, \pi$ with different magnitudes of the Mach number, namely $M=10^{-1}$, $M=10^{-2}$ and $M=10^{-3}$. The computational mesh is composed of $N_x \times N_y \times N_z = 80 \times 80 \times 4$ control volumes and the time step is evaluated with $\textnormal{CFL}=0.15$ according to \cite{Avgerinos2019}. Figure \ref{fig.Gresho} depicts the velocity magnitude contours together with the stream-traces of the velocity field for each Mach number regime. The pressure profile along the $x$ direction is also shown and compared against the exact solution. An excellent agreement can be appreciated, hence concluding that the novel semi-implicit pressure scheme preserves the stationary solution for a wide range of Mach numbers.   

\begin{figure}[!htbp]
	\begin{center}
		\begin{tabular}{cc} 
			\includegraphics[width=0.47\textwidth]{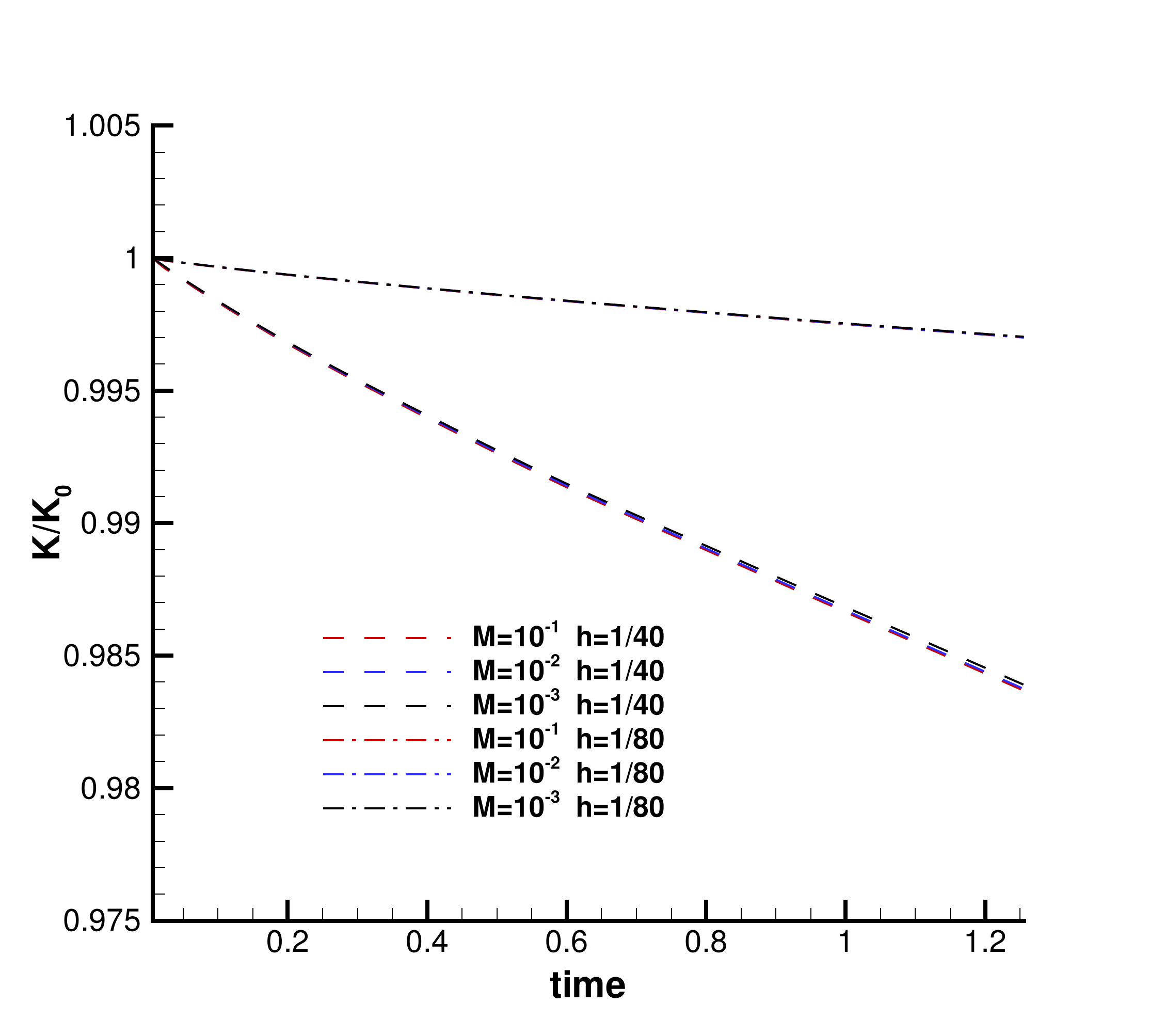} &
			\includegraphics[width=0.47\textwidth]{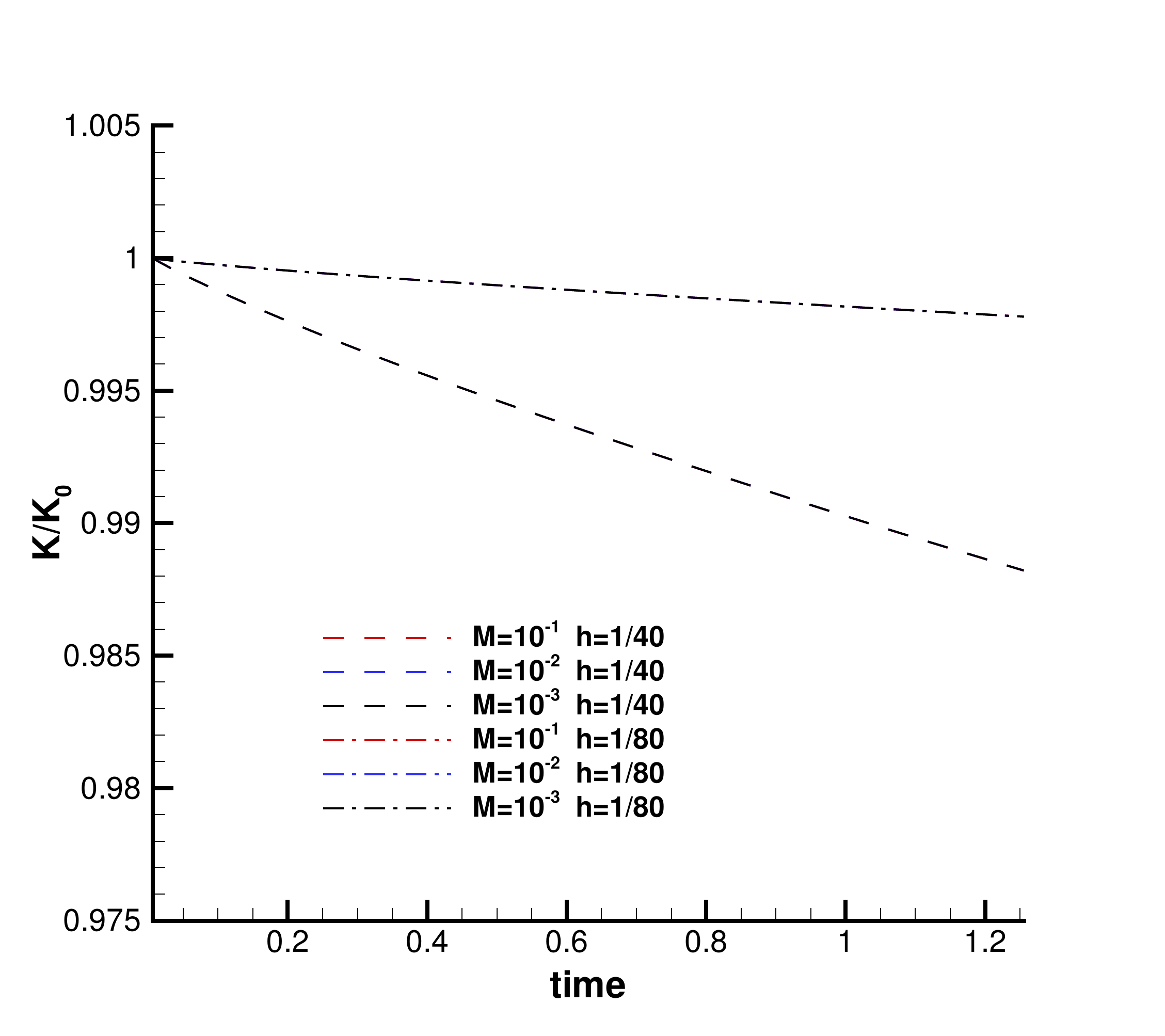}  \\			
		\end{tabular} 
		\caption{Evolution of the total kinetic energy $K$ normalized with respect to the initial kinetic energy $K_0$ of the Gresho vortex problem computed with second order (left) and third order (right) semi-implicit scheme. Dashed lines refer to the mesh size $N_x = N_y = 40$, while dash-dot lines refer to the mesh size $N_x = N_y = 80$. Mach number 0.1 (red), 0.01 (blue), 0.001 (black) are considered.}
		\label{fig.GreshoKENERGY}
	\end{center}
\end{figure}

Figure \ref{fig.GreshoKENERGY} depicts the evolution of the total kinetic energy normalized with with respect to the initial kinetic energy. Two different grids are used with characteristic mesh size of $h_1=1/40$ and $h_2=1/80$ for different values of the Mach number. The time step here is evaluated with $\textnormal{CFL}=0.25$. We consider both second order and third order schemes, in order to give evidences of the less dissipative behavior of the higher order scheme compared to the widespread second order solvers available in the literature for low Mach flows \cite{Avgerinos2019,BDLTV2020,Degond2}. The results do not depend on the stiffness regime because of the asymptotic property of the schemes, that allows all these simulations to be run with the same time step. Indeed, one can notice that the lines are almost overlapping, thus demonstrating that the low Mach regime does not affect neither the stability nor the accuracy of the numerical scheme. Compared to other results in the literature \cite{Avgerinos2019}, the high order accuracy of the SI-P method reduces the numerical dissipation, which is particularly evident on the coarser mesh $h_1$. The kinetic energy dissipation measures $0.997$ and $0.998$ for second and third order scheme, respectively, in the case of the finest mesh. The advantage induced by a higher order discretization is more evident on the coarse mesh, where the ratio $K/K_0$ is $0.984$ for $M=2$ and $0.988$ for $M=3$.

\subsection{Viscous shock} \label{ssec:ShockNS}
\begin{figure}[!htbp]
	\begin{center}
		\begin{tabular}{cc} 
			\includegraphics[width=0.47\textwidth]{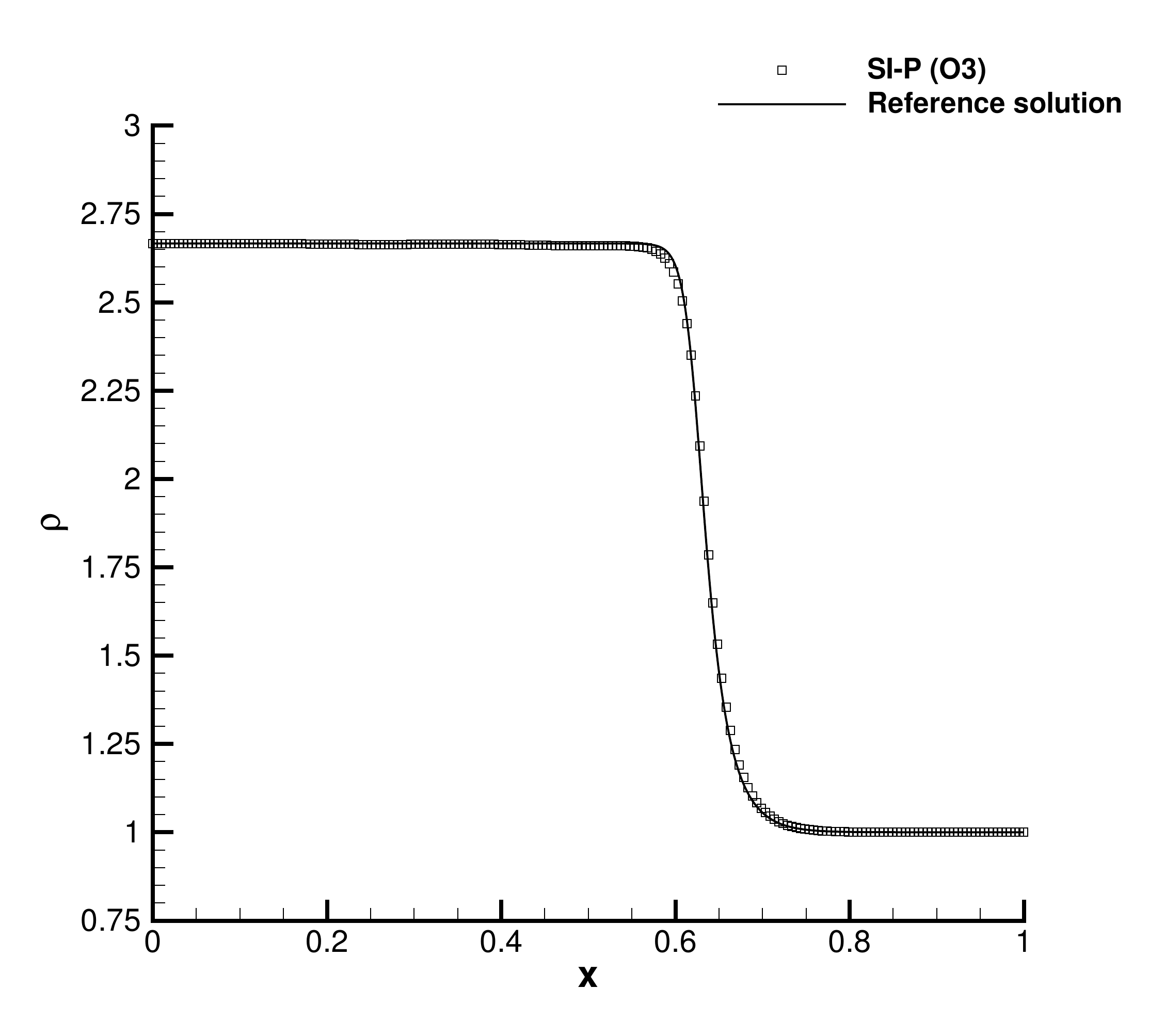} &           
			\includegraphics[width=0.47\textwidth]{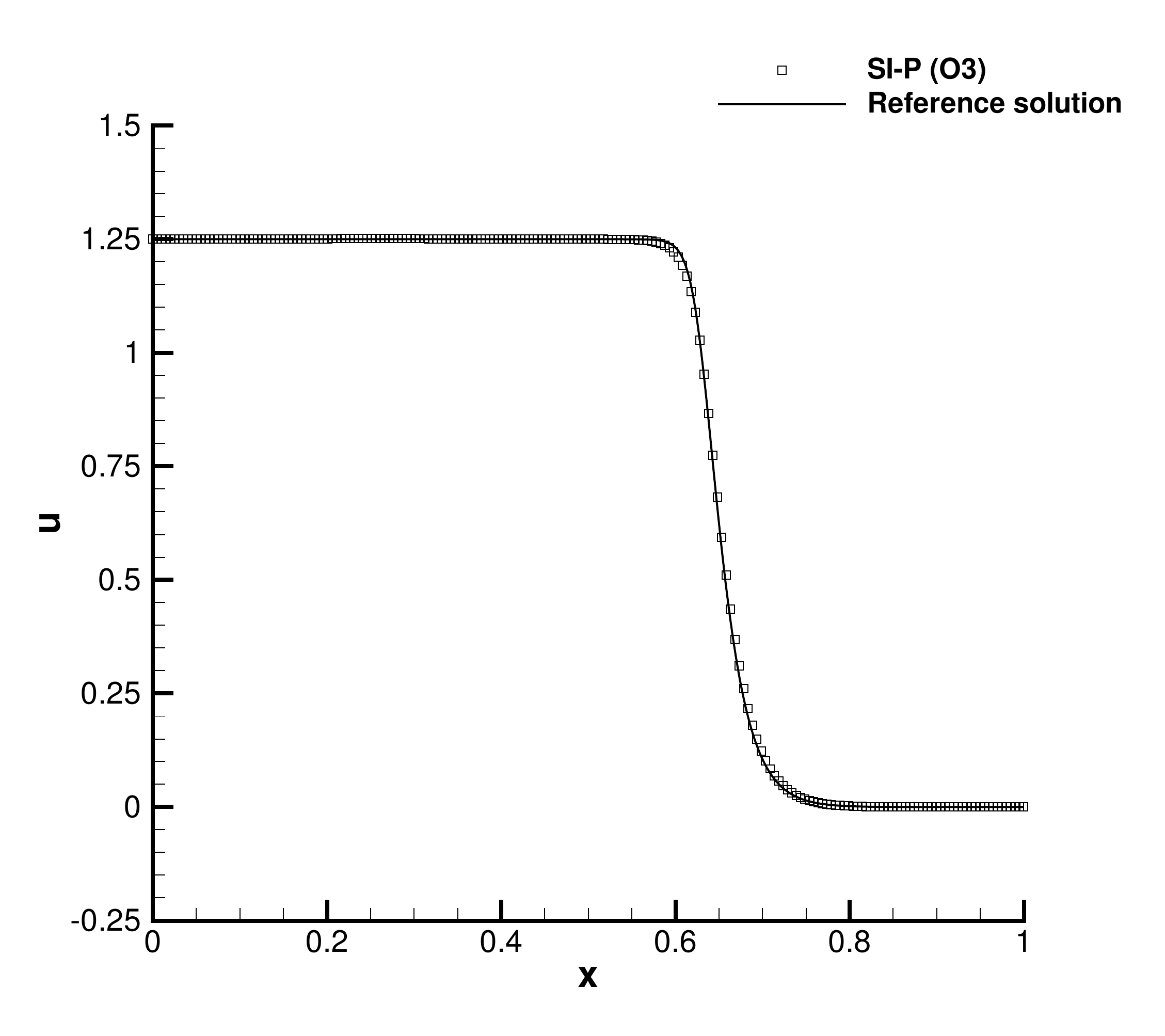} \\
			\includegraphics[width=0.47\textwidth]{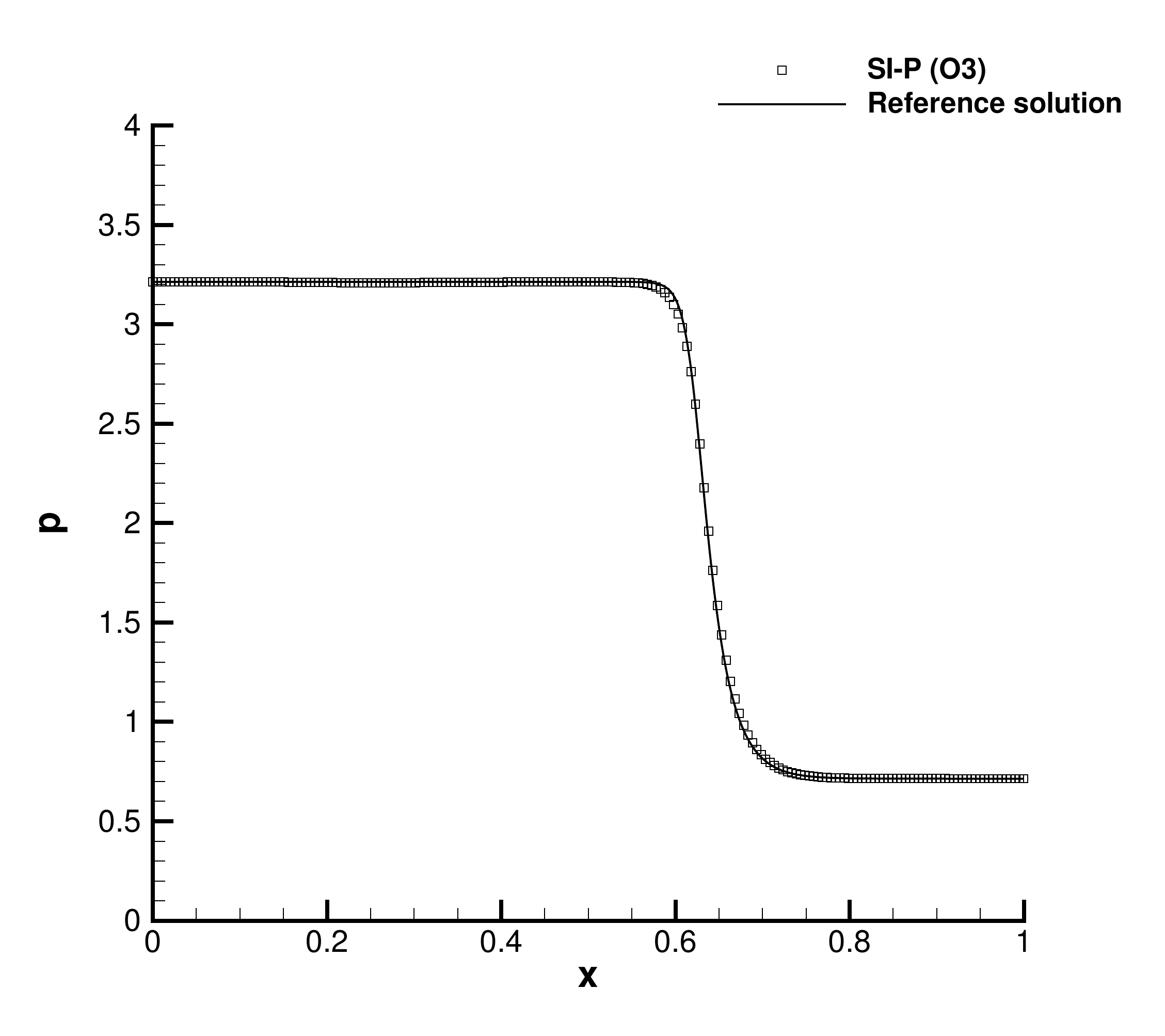} &           
			\includegraphics[width=0.47\textwidth]{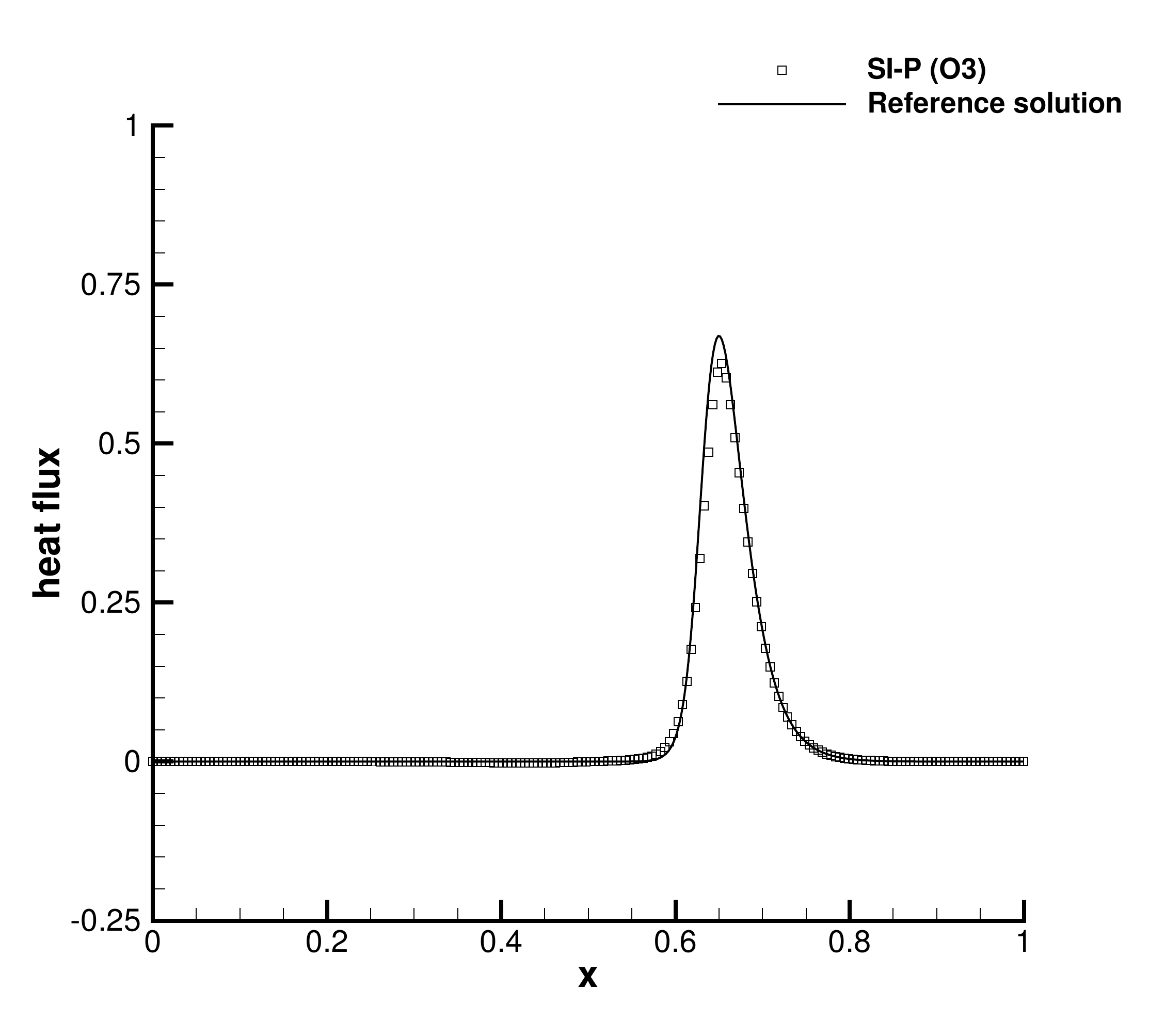} \\
		\end{tabular} 
		\caption{Viscous shock problem with shock Mach number $M_s = 2$ and Prandtl number $Pr = 0.75$. Third order SI-P solution compared against analytical solution for density (top left), horizontal velocity (top right), pressure (bottom left) and heat flux in x-direction (bottom right).	One-dimensional cut of 200 equidistant points along the x-direction at $y = z = 0$.}
		\label{fig.ShockNS}
	\end{center}
\end{figure}

Now we consider the full Navier-Stokes system \eqref{eqn.cns} in the case of supersonic viscous flows. Specifically, we propose to solve the problem of an isolated viscous
shock wave which is traveling into a medium at rest with a shock Mach number of $M_s=2$. The analytical solution of this problem has been obtained in \cite{Becker1923} where the compressible Navier-Stokes equations are solved for the special case of a stationary shock wave at Prandtl number $Pr= 0.75$ with constant viscosity. According to \cite{Becker1923,ADERNSE}, the exact solution is given in terms of dimensionless variables, namely density, pressure and velocity. The dimensionless velocity $\bar u = \frac{u}{M_s \, c_0}$ is related to the stationary shock wave. This can then be computed as the root of the following equation:
\begin{equation} 
\label{eqn.alg.u} 
\frac{|\bar u - 1|}{|\bar u - \lambda^2|^{\lambda^2}} = \left| \frac{1-\lambda^2}{2} \right|^{(1-\lambda^2)} 
\exp{\left( \frac{3}{4} \textnormal{Re}_s \frac{M_s^2 - 1}{\gamma M_s^2} x \right)},
\end{equation}
with
\begin{equation}
\lambda^2 = \frac{1+ \frac{\gamma-1}{2}M_s^2}{\frac{\gamma+1}{2}M_s^2}.
\end{equation}
The solution of equation \eqref{eqn.alg.u} permits to express the dimensionless velocity $\bar u$ as a function of $x$. The form of the viscous profile of the dimensionless pressure $\bar p = \frac{p-p_0}{\rho_0 c_0^2 M_s^2}$ is given by the relation 
\begin{equation}
\label{eqn.alg.p} 
\bar p = 1 - \bar u +  \frac{1}{2 \gamma}
\frac{\gamma+1}{\gamma-1} \frac{(\bar u - 1 )}{\bar u} (\bar u - \lambda^2).  
\end{equation}
Finally, the profile of the dimensionless density $\bar \rho = \frac{\rho}{\rho_0}$ is derived from the integrated continuity equation: $\bar \rho \bar u = 1$. In order to simulate an unsteady shock wave traveling into a medium at rest, one can simply superimpose a constant velocity field $u = M_s c_0$ to the solution of the stationary shock wave found in the previous steps. The computational domain is the rectangular box  $\Omega=[0;1]\times[0;0.2]\times[0;0.2]$ which is discretized with a total number of cells $N_x \times N_y \times N_z =200 \times 4 \times 4$. Periodic boundaries are imposed in $y$ and $z$ direction, while the constant inflow velocity is prescribed for $x=0$ and outflow boundary condition is set at $x=1$. The time step is evaluated with $\textnormal{CFL}=0.5$ and the final time of the simulation is $t_f=0.2$. The initial condition is given by a shock wave centered at $x=0.25$ which is propagating at Mach $M_s=2$ from left to right with a Reynolds number of $Re=100$. The upstream shock state is defined by $\rho_0=1$, $u_0=v_0=0$, $p_0=1/\gamma$ and $c_0=1$, while the fluid viscosity is $\mu=2\times 10^{-2}$ and the ideal gas law is adopted, thus allowing the heat flux in \eqref{eqn.cns} to be treated implicitly in the pressure wave equation \eqref{eqn.p_SI-discr}. The third order SI-P schemes is used to run the simulation and the results are depicted in Figure \ref{fig.ShockNS}, which match very closely the analytical solution for density, horizontal velocity, pressure and heat flux in $x$ direction computed as $q_x=\lambda \, \frac{\partial T}{\partial x}$. Though being a one-dimensional problem, let observe that this test case involves all terms contained in the governing equations, hence including convective and viscous fluxes, pressure gradients as well as temperature gradients and heat fluxes. Furthermore, an analytical solution does exist which permits to compare the numerical results. Looking at the excellent matching between numerical and exact solution we can conclude that the Navier-Stokes system is properly discretized by the novel semi-implicit pressure solver proposed in this article. 

\subsection{3D Taylor-Green vortex} \label{ssec:TGV}
\begin{figure}[!htbp]
	\begin{center}
		\begin{tabular}{cc} 
			\includegraphics[width=0.47\textwidth]{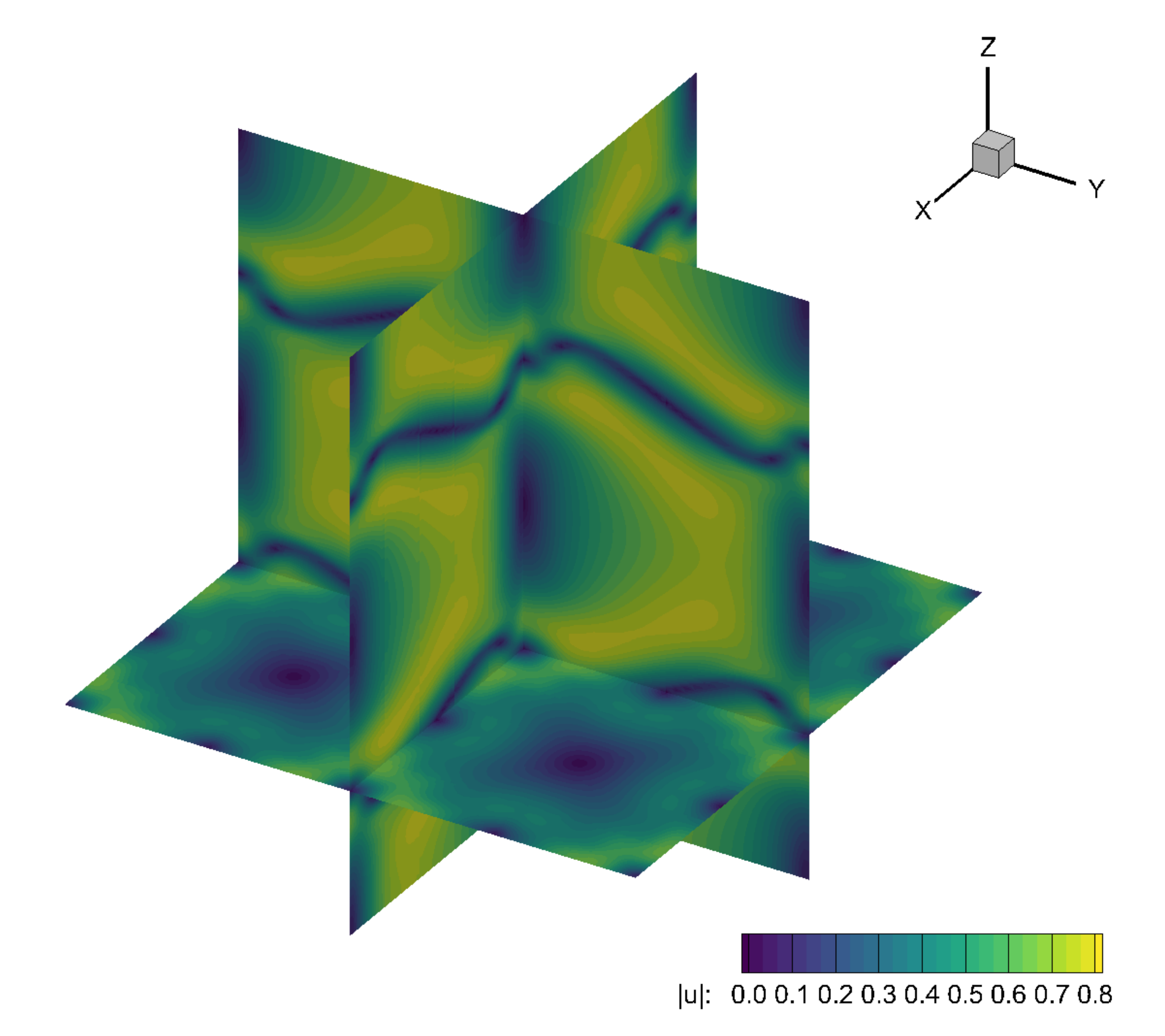} &           
			\includegraphics[width=0.47\textwidth]{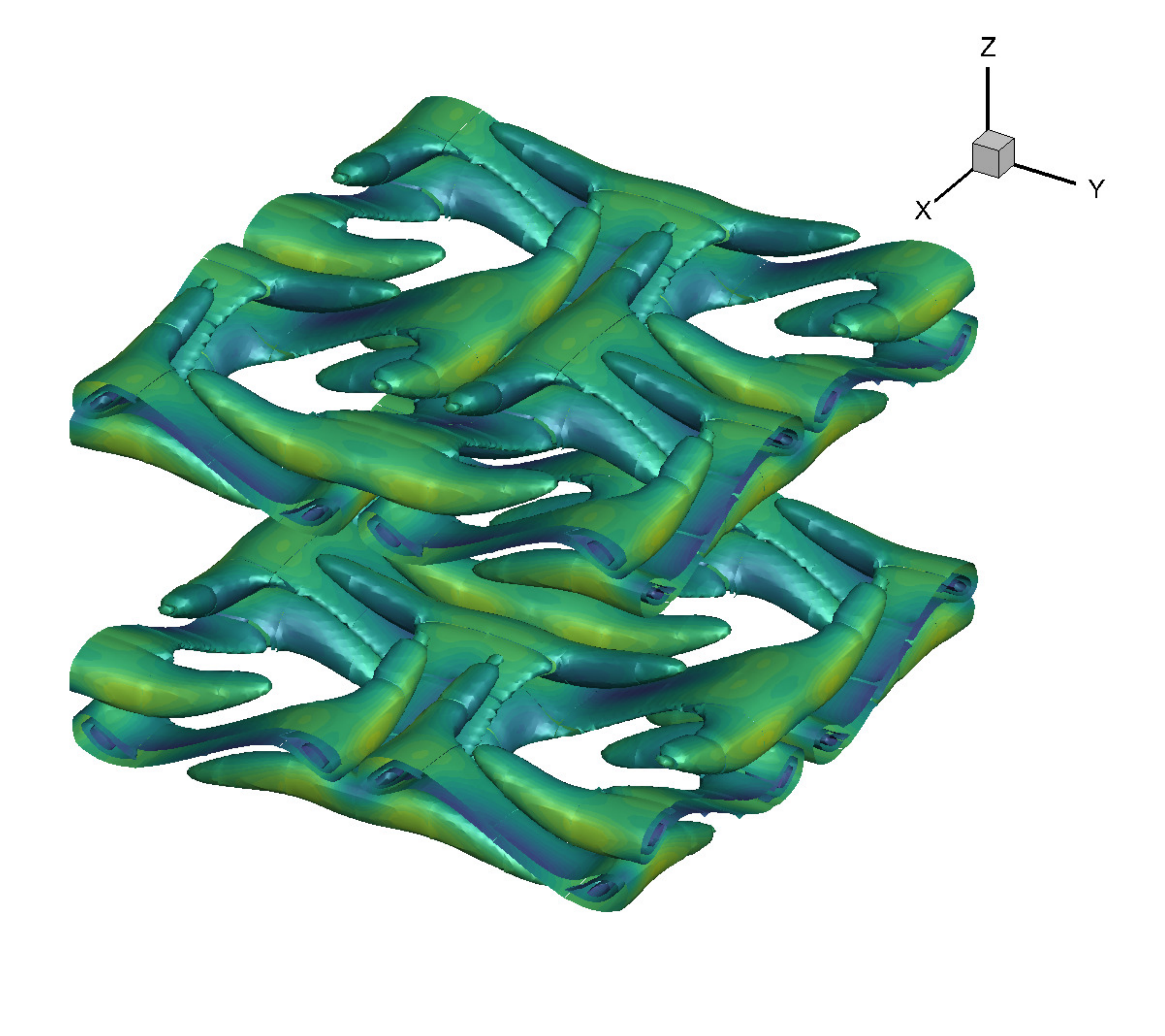} \\
			\includegraphics[width=0.47\textwidth]{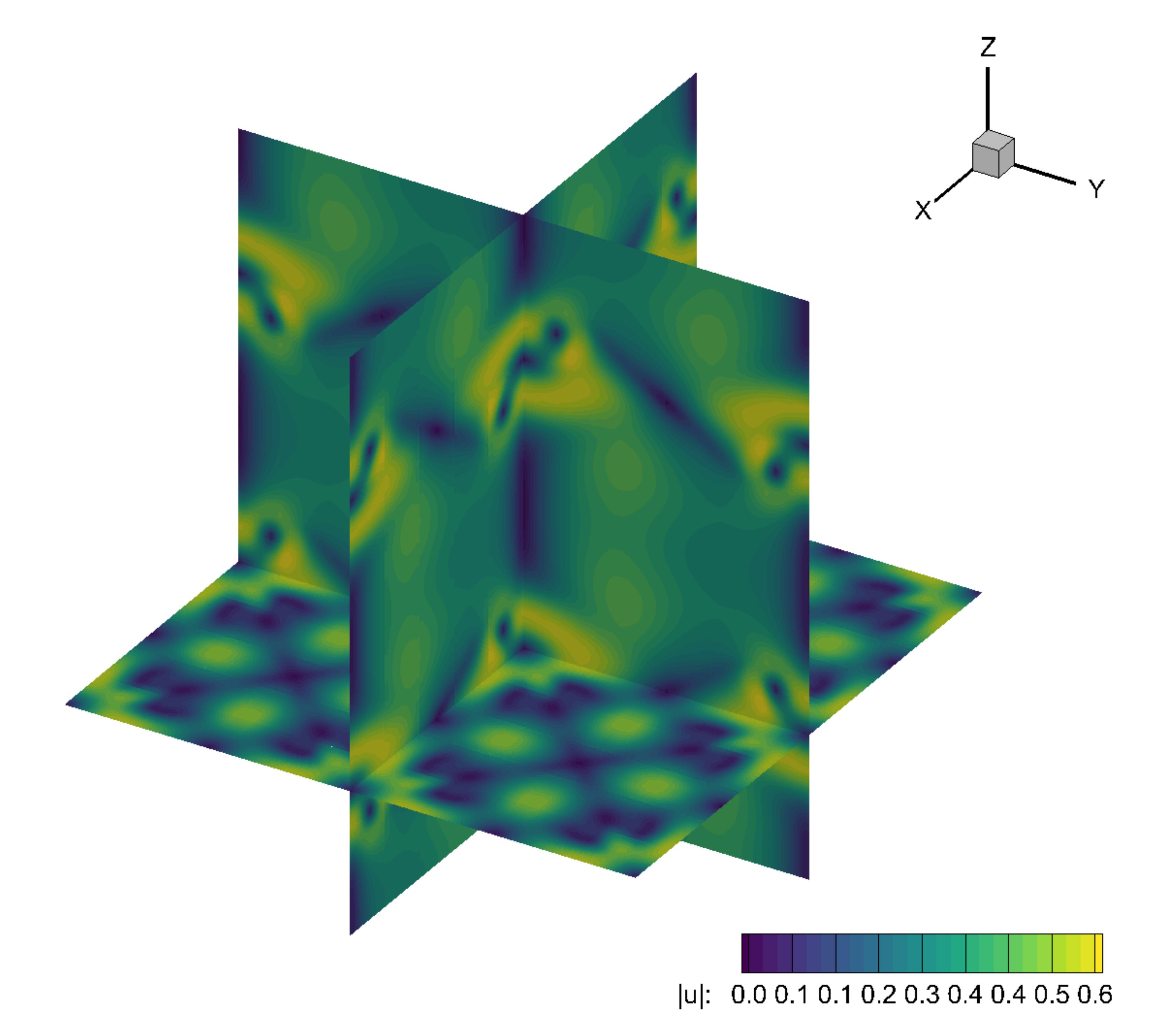} &           
			\includegraphics[width=0.47\textwidth]{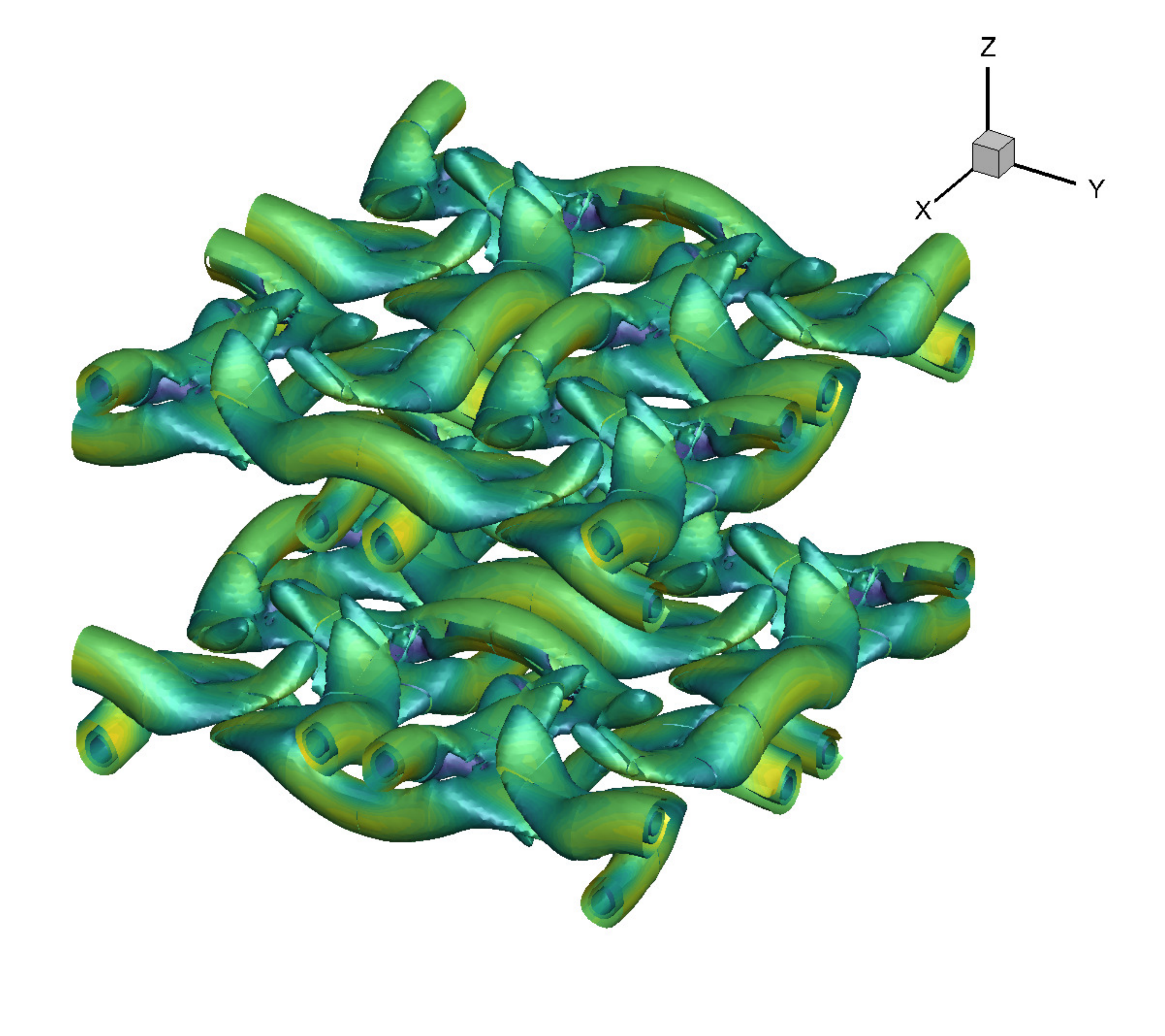} \\
		\end{tabular} 
		\caption{3D Taylor-Green vortex at $Re=200$. Velocity magnitude (left) and vorticity iso-surfaces at levels $\{2,3,5\}$ (right) at time $t=4$ (top) and $t=8$ (bottom).}
		\label{fig.TGV}
	\end{center}
\end{figure}

As last test case we solve the well-known 3D Taylor-Green vortex, that is a widespread test problem used in the context of incompressible flows. The initial condition of the fluid according to \cite{Oriol2015} writes
\begin{eqnarray}
&& \rho(\xx,0) = \rho_0,  \nonumber \\ 
&& \uu(\xx,0) = \left( \, \, \sin(x)\cos(y)\cos(z),  \, \, -\cos(x)\sin(y)\cos(z), \, \, 0 \right),  \nonumber \\ 
&& p(\xx,0) = p_0 + \frac{\rho_0}{16}\left(\cos(2x)+\cos(2y) \right)\left(\cos(2z)+2) \right), 
\label{eq:NT_13_1}
\end{eqnarray}
with $\rho_0=1$. The computational domain is the cube $\Omega=[-\pi,\pi]^3$ and periodic boundary conditions are imposed everywhere. Starting from this smooth initial condition, the flow quickly degenerates into very complex small scale structures, depending on the Reynolds number. Consequently, no analytical solution is available for this highly unsteady flow. Nevertheless, in \cite{BrachetDNS} well-resolved DNS studies for an incompressible fluid are available, hence we consider those results as a very accurate reference solution. In order to mimic the incompressible property of the fluid we set $p_0=10^3$ in our semi-implicit compressible solver. The final time of the simulation is set to $t_f=10$ and the time step is evaluated with $\textnormal{CFL}=0.5$. Two different values of Reynolds number are taken into account, namely $Re=100$ and $Re=200$. The computational domain is discretized with a total number of $1'728'000$ control volumes, obtained by setting $N_x \times N_y \times N_z =120 \times 120 \times 120$. Figure \ref{fig.TGV} depicts the vorticity isosurfaces together with the velocity magnitude for $Re=200$, clearly showing the development of the small-scale structures that arise from the fluid flow.

\begin{figure}[!htbp]
	\begin{center}
		\begin{tabular}{c} 
			\includegraphics[width=0.7\textwidth]{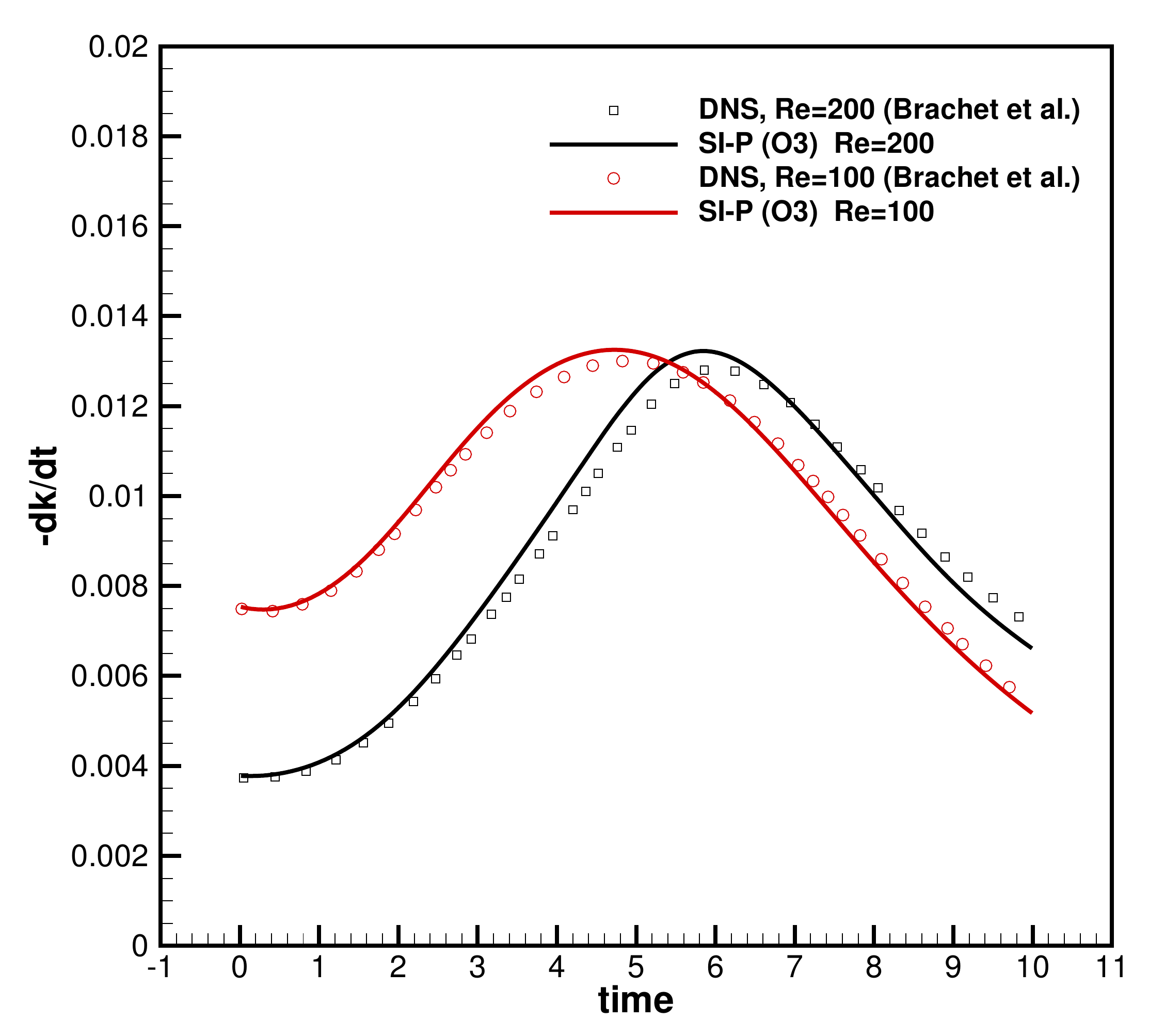}\\      
		\end{tabular} 
		\caption{Time evolution of the kinetic energy dissipation rate $-dK/dt$ for the 3D Taylor-Green vortex compared with available DNS data of brachet et al. \cite{BrachetDNS}.}
		\label{fig.TGV_kinEnergy}
	\end{center}
\end{figure}

Finally, Figure \ref{fig.TGV_kinEnergy} plots the time series of the calculated total kinetic energy dissipation rates $-dK/dt$ compared against the DNS data \cite{BrachetDNS}. Also for this rather complex test case, the numerical results obtained with the third order SI-P schemes fit well with the reference solution for all the considered Reynolds numbers. Notice that the spatial resolution used for running this test is relatively coarse compared to existing second order solvers. However, the results are of good quality because of the high order discretization achieved by the novel semi-implicit schemes. This has also been observed in the case of higher order DG schemes \cite{TavelliIncNS,TavelliCNS} where the spatial resolution could be reduced even further. On the other hand, second order schemes as the method presented in \cite{BDLTV2020} require much more computational cells to carry out the same simulations of the 3D Taylor-Green vortex test case.

\section{Conclusions} \label{sec:conclusion}
A novel high-order semi-implicit numerical method for the solution of the compressible Navier-Stokes equations at all Mach numbers has been derived and discussed. A high order cell-centered quadrature-free finite volume scheme is used for the approximation of the explicit fluxes, whereas finite differences are employed for the discretization of the implicit terms. Collocated Cartesian grids are used in a fully three-dimensional setting and the fluid can be modeled with ideal gas as well as with general equations of state that might lead to a nonlinear relation between internal energy and pressure. The new semi-implicit method splits the kinetic energy and the enthalpy fluxes into an explicit and an implicit part, making the usage of an iterative solver for the pressure unnecessary. Formal analysis of the scheme at the discrete level reveals the asymptotic property of the algorithm, which is capable of retrieving at the discrete level a consistent discretization of the limit model in the zero Mach number regime. High order time stepping is performed relying on IMEX schemes and the resulting stability condition requires a time step limitation based only on the fluid velocity (and eventually the viscous eigenvalues) and not on the acoustic speed, thus making the novel numerical method very efficient and suitable for the simulation of low Mach flows. The implicit sub-system requires the solution of an elliptic equation for the pressure, which allows to easily include nonlinear EOS by adopting a nested Newton solver for the resulting mildly nonlinear system. This would not hold true in the case of the derivation of an implicit equation for the total energy, where the system would become fully nonlinear. A wide set of benchmark problems is proposed to test the accuracy and the robustness of the new algorithm, involving low and high Mach number flows as well as viscous and inviscid fluid simulations.

Future research will concern the solution of more complex systems of hyperbolic equations with nonconservative products and stiff source terms like the GPR model \cite{GPR}, which would simultaneously allow a unified formulation for continuum mechanics, including fluids and solids. Another potential future topic of research may be the application of the new semi-implicit pressure-based IMEX solver to the equations of magnetohydrodynamics with involution constraints, which require the additional property of keeping a zero divergence of the magnetic field at the discrete level.  

%
\section*{Acknowledgments}
The Authors would like to thank the Italian Ministry of Instruction, University and Research (MIUR) to support this research with funds coming from PRIN Project 2017 (No. 2017KKJP4X entitled “Innovative numerical methods for evolutionary partial differential equations and applications”). The Authors also acknowledge the CINECA supercomputing center in Bologna (Italy) for awarding access to the MARCONI100 cluster under the project \textit{IscrC\_HiPPUM\_0}.


\bibliography{biblio}

\begin{thebibliography}{10}

\bibitem{HLL1983}
P.~Lax A.~Harten and B.~van Leer.
\newblock On upstream differencing and godunov-type schemes for hyperbolic
  conservation laws.
\newblock {\em SIAM Rev.}, 25:35--61, 1983.

\bibitem{Ala2}
T.~Alazard.
\newblock Low mach number limit of the full navier-stokes equations.
\newblock {\em Arch. Ration. Mech. Anal.}, 180:1--73, 2006.

\bibitem{Avgerinos2019}
S.~Avgerinos, F.~Bernard, A.~Iollo, and G.~Russo.
\newblock Linearly implicit all mach number shock capturing schemes for the
  euler equations.
\newblock {\em J. Comput. Phys.}, 393:278--312, 2019.

\bibitem{TavelliCNS}
P.~Roe B.~Einfeldt, C.~Munz and B.~Sj{\"o}green.
\newblock On godunov-type methods near low densities.
\newblock {\em J. Comp. Phys.}, 92:273--295, 1991.

\bibitem{BalsaraFlattener}
D.S. Balsara.
\newblock {Self-adjusting, positivity preserving high order schemes for
  hydrodynamics and magnetohydrodynamics}.
\newblock {\em J. Comp. Phys.}, 231:7504 -- 7517, 2012.

\bibitem{Becker1923}
R.~Becker.
\newblock Stosswelle und detonation.
\newblock {\em Physik}, 8:321--1923, 1923.

\bibitem{BilletAbgrall2003}
G.~Billet and R.~Abgrall.
\newblock {A}n adaptive shock-capturing algorithm for solving unsteady reactive
  flows.
\newblock {\em Computers \& Fluids}, 32:1473--1495, 2003.

\bibitem{BosRus}
S.~Boscarino and G.~Russo.
\newblock On a class of uniformly accurate {IMEX} {R}unge-{K}utta schemes and
  applications to hyperbolic systems with relaxation.
\newblock {\em SIAM J. Sci. Comput.}, 31:1926--1945, 2009.

\bibitem{BPR2017}
Sebastiano Boscarino, Lorenzo Pareschi, and Giovanni Russo.
\newblock A unified {IMEX} {R}unge-{K}utta approach for hyperbolic systems with
  multiscale relaxation.
\newblock {\em SIAM J. Numer. Anal.}, 55(4):2085--2109, 2017.

\bibitem{ADERFSE}
W.~Boscheri.
\newblock A space-time semi-lagrangian advection scheme on staggered voronoi
  meshes applied to free surface flows.
\newblock {\em Computers \& Fluids}, 202:104503, 2020.

\bibitem{BoscheriWAO}
W.~Boscheri and D.S. Balsara.
\newblock {High order direct Arbitrary-Lagrangian-Eulerian (ALE) $P_NP_M$
  schemes with WENO Adaptive-Order reconstruction on unstructured meshes}.
\newblock {\em J. Comp. Phys.}, 398:108899, 2019.

\bibitem{BDT_cns}
W.~Boscheri, G.~Dimarco, and M.~Tavelli.
\newblock An efficient all mach second order finite volume solver for
  compressible navier-stokes equations.
\newblock {\em Computer Methods in Applied Mechanics and Engineering}.
\newblock submitted.

\bibitem{NestedNewton}
L.~Brugnano and V.~Casulli.
\newblock Iterative solution of piecewise linear systems.
\newblock {\em SIAM J. Sci. Comput.}, 30:463–472, 2007.

\bibitem{Chalons}
M.~Girardin C.~Chalons and S.~Kokh.
\newblock Large time step and asymptotic preserving numerical schemes for the
  gas dynamics equations with source terms.
\newblock {\em SIAM J. Sci. Comput.}, 35:2874--2902, 2013.

\bibitem{Dou}
S.~Jiang C.~Dou and Y.~Ou.
\newblock Low mach number limit of full navier-stokes equations in a 3d bounded
  domain.
\newblock {\em J. Diff. Equat.}, 258:379--398, 2015.

\bibitem{Casulli1990}
V.~Casulli.
\newblock Semi-implicit finite difference methods for the two--dimensional
  shallow water equations.
\newblock {\em J. Comp. Phys.}, 86:56--74, 1990.

\bibitem{Casulli1999}
V.~Casulli.
\newblock A semi-implicit finite difference method for non-hydrostatic
  free-surface flows.
\newblock {\em Int. J. Num. Meth. in Fluids}, 30:425--440, 1999.

\bibitem{CPSV:2018}
I.~Cravero, G.~Puppo, M.~Semplice, and G.~Visconti.
\newblock {CWENO}: uniformly accurate reconstructions for balance laws.
\newblock {\em Math. Comp.}, 87:1689--1719, 2018.

\bibitem{DegTan}
P.~Degond and M.~Tang.
\newblock All speed scheme for the low {M}ach number limit of the isentropic
  {E}uler equations.
\newblock {\em Commun. Comput. Phys.}, 10:1--31, 2011.

\bibitem{Dellacherie1}
S.~Dellacherie.
\newblock Analysis of godunov type schemes applied to the compressible euler
  system at low mach number.
\newblock {\em J. Comp. Phys.}, 229:978--1016, 2010.

\bibitem{ADERNSE}
M.~Dumbser.
\newblock Arbitrary high order {PNPM} schemes on unstructured meshes for the
  compressible {Navier--Stokes} equations.
\newblock {\em Computers \& Fluids}, 39:60--76, 2010.

\bibitem{ADER_CWENO}
M.~Dumbser, W.~Boscheri, M.~Semplice, and G.~Russo.
\newblock Central weighted {ENO} schemes for hyperbolic conservation laws on
  fixed and moving unstructured meshes.
\newblock {\em {SIAM} J. Sci. Comput.}, 39(6):A2564--A2591, 2017.

\bibitem{Dumbser_Casulli16}
M.~Dumbser and V.~Casulli.
\newblock A conservative, weakly nonlinear semi-implicit finite volume scheme
  for the compressible navier-stokes equations with general equation of state.
\newblock {\em Applied Mathematics and Computation}, 272:479--497, 2016.

\bibitem{DumbserEnauxToro}
M.~Dumbser, C.~Enaux, and E.F. Toro.
\newblock Finite volume schemes of very high order of accuracy for stiff
  hyperbolic balance laws.
\newblock {\em J. Comp. Phys.}, 227:3971--4001, 2008.

\bibitem{Ioriatti2015}
M.~Dumbser, U.~Iben, and M.~Ioriatti.
\newblock An efficient semi-implicit finite volume method for axially symmetric
  compressible flows in compliant tubes.
\newblock {\em Applied Numerical Mathematics}, 89:24--44, 2015.

\bibitem{DumbserKaeser2007}
M.~Dumbser and M.~K{\"a}ser.
\newblock Arbitrary high order non-oscillatory finite volume schemes on
  unstructured meshes for linear hyperbolic systems.
\newblock {\em J. Comp. Phys.}, 221:693--723, 2007.

\bibitem{Toro_Vazquez12}
M.E. V{\'a}zquez-Cend{\'o}n E.F.~Toro.
\newblock Flux splitting schemes for the euler equations.
\newblock {\em Computers \& Fluids}, 70:1--12, 2012.

\bibitem{Degond2}
P.~Degond F.~Cordier and A.~Kumbaro.
\newblock An asymptotic-preserving all-speed scheme for the euler and navier
  stokes equations.
\newblock {\em J. Comp. Phys.}, 231:5685--5704, 2012.

\bibitem{Fambri2017}
F.~Fambri and M.~Dumbser.
\newblock Semi-implicit discontinuous galerkin methods for the incompressible
  navier--stokes equations on adaptive staggered cartesian grids.
\newblock {\em Computer Methods in Applied Mechanics and Engineering},
  324:170--203, 2017.

\bibitem{Fambri2017b}
F.~Fambri, M.~Dumbser, and O.~Zanotti.
\newblock Space--time adaptive ader-dg schemes for dissipative flows:
  Compressible navier--stokes and resistive mhd equations.
\newblock {\em Computer Physics Communications}, 220:297--318, 2017.

\bibitem{DLDV2018}
V.M.~Dansac G.~Dimarco, R.~Loub{\`{e}}re and M.H. Vignal.
\newblock Second-order implicit-explicit total variation diminishing schemes
  for the euler system in the low mach regime.
\newblock {\em J. Comp. Phys.}, 372:178--201, 2018.

\bibitem{GassnerDiffusion}
F.~L{\"o}rcher G.~Gassner and C.D. Munz.
\newblock A contribution to the construction of diffusion fluxes for finite
  volume and discontinuous galerkin schemes.
\newblock {\em J. Comp. Phys.}, 224:1049--1063, 2007.

\bibitem{Godunov1959}
S.~Godunov.
\newblock Finite difference methods for the computation of discontinuous
  solutions of the equations of fluid dynamics.
\newblock {\em Mat. Sb.}, 47:271--306, 1959.

\bibitem{Gresho1990}
P.~M. Gresho and S.T. Chan.
\newblock On the theory of semi-implicit projection methods for viscous
  incompressible flow and its implementation via a finite element method that
  also introduces a nearly consistent mass matrix. part 2: Implementation.

\bibitem{GH}
H.~Guillard and A.~Murrone.
\newblock On the behavior of upwind schemes in the low mach number limit : Ii.
  godunov type schemes.
\newblock {\em Computers \& Fluids}, 33:655--675, 2004.

\bibitem{Guillard}
H.~Guillard and C.~Viozat.
\newblock On the behavior of upwind schemes in the low mach limit.
\newblock {\em Computers \& Fluids}, 28:63--86, 1999.

\bibitem{HuShuTri}
C.~Hu and {C.W.} Shu.
\newblock Weighted essentially non-oscillatory schemes on triangular meshes.
\newblock {\em J. Comp. Phys.}, 150:97--127, 1999.

\bibitem{Ioriatti2018}
M.~Ioriatti and M.~Dumbser.
\newblock Semi-implicit staggered discontinuous galerkin schemes for axially
  symmetric viscous compressible flows in elastic tubes.
\newblock {\em Computers and Fluids}, 167:166--179, 2018.

\bibitem{shu_efficient_weno}
{G.-S.} Jiang and {C.W.} Shu.
\newblock Efficient implementation of weighted {ENO} schemes.
\newblock {\em J. Comp. Phys.}, 126:202--228, 1996.

\bibitem{JINAP1999}
Shi Jin.
\newblock Efficient asymptotic-preserving ({AP}) schemes for some multiscale
  kinetic equations.
\newblock {\em SIAM J. Sci. Comput.}, 21(2):441--454, 1999.

\bibitem{KLARAP1999}
Axel Klar.
\newblock An asymptotic preserving numerical scheme for kinetic equations in
  the low {M}ach number limit.
\newblock {\em SIAM J. Numer. Anal.}, 36(5):1507--1527, 1999.

\bibitem{Klein}
R.~Klein.
\newblock Semi-implicit extension of a godunov-type scheme based on low mach
  number asymptotics i: One-dimensional flow.
\newblock {\em J. Comp. Phys.}, 121:213--237, 1995.

\bibitem{LW1960}
P.~Lax and B.~Wendroff.
\newblock Systems of conservation laws.
\newblock {\em J. Comp. Phys.}, 13:217--237, 1960.

\bibitem{LeVequeBook}
{R.J.} LeVeque.
\newblock {\em Finite Volume Methods for Hyperbolic Problems}.
\newblock Cambridge University Press, 2002.

\bibitem{LPR:99}
D.~Levy, G.~Puppo, and G.~Russo.
\newblock Central {WENO} schemes for hyperbolic systems of conservation laws.
\newblock {\em M2AN Math. Model. Numer. Anal.}, 33(3):547--571, 1999.

\bibitem{Bog}
R.~Klein M.~Boger, F.~Jaegle and C.-D. Munz.
\newblock Coupling of compressible and incompressible flow regions using the
  multiple pressure variables approach.
\newblock {\em Math. Methods Appl. Sci.}, 38:458--477, 2015.

\bibitem{AMR2013}
A.~Hidalgo M.~Dumbser, O.~Zanotti.
\newblock Ader-weno finite volume schemes with space-time adaptive mesh
  refinement.
\newblock {\em J. Comp. Phys.}, 248:257--286, 2013.

\bibitem{BrachetDNS}
D.I.~Meiron M.E.~Brachet and S.A. Orszag.
\newblock Small-scale structure of the taylor-green vortex.
\newblock {\em Journal of Fluid Mechanics}, 130:411--452, 1983.

\bibitem{Munz1994}
C.D. Munz.
\newblock On godunov-type schemes for lagrangian gas dynamics.
\newblock {\em SIAM J. Numer. Anal.}, 31:17--42, 1994.

\bibitem{Oriol2015}
R.~Codina O.~Colomes, S.~Badia and J.~Principe.
\newblock Assessment of variational multiscale models for the large eddy
  simulation of turbulent incompressible flows.
\newblock {\em Comp. Methods in App. Mech. and Eng.}, 285:32--63, 2015.

\bibitem{Hofer}
S.~Osher and F.~Solomon.
\newblock A partially implicit method for large stiff systems of ode's with
  only few equations introducing small time-constants.
\newblock {\em SIAM J. Numer. Anal.}, 13:645--663, 1976.

\bibitem{OS1982}
S.~Osher and F.~Solomon.
\newblock Upwind difference schemes for hyperbolic conservation laws.
\newblock {\em Math. Comput.}, 38:339--374, 1997.

\bibitem{Deg_Del_Sang_Vig}
A.~Sangam P.~Degond, F.~Deluzet and M.-H. Vignal.
\newblock An asymptotic preserving scheme for the euler equations in a strong
  magnetic field.
\newblock {\em J. Comp. Phys.}, 228:3540--3558, 2009.

\bibitem{PR_IMEXHO}
L.~Pareschi and G.~Russo.
\newblock High order asymptotically strong-stability-preserving methods for
  hyperbolic systems with stiff relaxation. in: Hou t.y., tadmor e. (eds)
  hyperbolic problems: Theory, numerics, applications.
\newblock 2003.

\bibitem{PR_IMEX}
L.~Pareschi and G.~Russo.
\newblock Implicit-explicit runge-kutta schemes and applications to hyperbolic
  systems with relaxation.
\newblock {\em J. Sci. Comput.}, 25:129--155, 2005.

\bibitem{Patankar}
S.V. Patankar.
\newblock {\em Numerical heat transfer and fluid flow}.
\newblock New York: McGraw-Hill, 1980.

\bibitem{GPR}
I.~Peshkov and E.~Romenski.
\newblock A hyperbolic model for viscous newtonian flows.
\newblock {\em Continuum Mech Thermodyn}, 28:85–104, 2016.

\bibitem{BosFil2016}
F.~Filbet S.~Boscarino and G.~Russo.
\newblock High order semi-implicit schemes for time dependent partial
  differential equations.
\newblock {\em J. Sci. Comput.}, 68:975--1001, 2016.

\bibitem{BosRus2019}
G.~Russo S.~Boscarino, J.-M.~Qiu and T.~Xiong.
\newblock A high order semi-implicit imex weno scheme for the all-mach
  isentropic euler system.
\newblock {\em J. Comp. Phys.}, 392:594--618, 2019.

\bibitem{Boscarino2019}
G.~Russo S.~Boscarino, J.-M.~Qiu and T.~Xiong.
\newblock A high order semi-implicit imex weno scheme for the all-mach
  isentropic euler system.
\newblock {\em J. Comp. Phys.}, 392:594--618, 2019.

\bibitem{Gottlieb2001}
C.-W.~Shu S.~Gottlieb and E.~Tadmor.
\newblock Strong stability-preserving high-order time discretization methods.
\newblock {\em SIAM Rev.}, 43:89--112, 2001.

\bibitem{stroud}
{A.H.} Stroud.
\newblock {\em Approximate Calculation of Multiple Integrals}.
\newblock Prentice-Hall Inc., Englewood Cliffs, New Jersey, 1971.

\bibitem{TavelliIncNS}
M.~Tavelli and M.~Dumbser.
\newblock A staggered space-time discontinuous galerkin method for the
  three-dimensional incompressible navier-stokes equations on unstructured
  tetrahedral meshes.
\newblock {\em J. Comp. Phys.}, 319:294--323, 2016.

\bibitem{ToroBook}
{E.F.} Toro.
\newblock {\em Riemann Solvers and Numerical Methods for Fluid Dynamics: a
  Practical Introduction.}
\newblock Springer, 2009.

\bibitem{Turkel}
E.~Turkel.
\newblock Preconditioned methods for solving the incompressible and low speed
  compressible equations.
\newblock {\em J. Comp. Phys.}, 72:277--298, 1987.

\bibitem{AscRuuSpi}
S.~J.~Ruuth U.~M.~Ascher and R.~J. Spiteri.
\newblock Implicit-explicit {R}unge-{K}utta methods for time-dependent partial
  differential equations.
\newblock {\em Appl. Numer. Math.}, 25:151--167, 1982.

\bibitem{Vidal}
J.~Vidal.
\newblock {\em Thermodynamics: Applications in Chemical Engineering and the
  Petroleum Industry}.
\newblock Editions Technip, 2003.

\bibitem{BDLTV2020}
R.~Loub{\`{e}}re M.~Tavelli W.~Boscheri, G.~Dimarco and M.H. Vignal.
\newblock A second order all mach number imex finite volume solver for the
  three dimensional euler equations.
\newblock {\em J. Comp. Phys.}, 415:109486, 2020.

\end{thebibliography}
\bibliographystyle{plain}

\end{document}